%% file: main.tex
\documentclass{article}

\usepackage{xcolor}
\definecolor{myred}{RGB}{214, 39, 24}
\definecolor{mygreen}{RGB}{0, 122, 49}
\usepackage[colorlinks,linkcolor=myred,citecolor=blue, pagebackref=true, hypertexnames=false]{hyperref}
\usepackage[margin=1.1in]{geometry}
\usepackage[numbers]{natbib}
\usepackage{booktabs}

\input{preamble}

\usepackage{authblk}

\title{Instance-Log-Optimality of Portfolio-Based\\E-Processes and their Sequential Hypothesis Tests}
\author{
  Ian Waudby-Smith$^\dagger$, Ricardo J.~Sandoval$^\dagger$, and Michael I.~Jordan$^{\dagger\ddagger}$\\
  $^\dagger$University of California, Berkeley\\
  $^\ddagger$Inria, Paris
}
\date{\today}

\begin{document}
\maketitle
\setcounter{tocdepth}{2}
\makeatletter
\renewcommand\tableofcontents{%
  \@starttoc{toc}%
}

\makeatother

\begin{abstract}
\input{abstract}
\end{abstract}

\input{content}

\input{relatedwork}

\input{conclusions}

\input{acks}

\bibliographystyle{plainnat}
\bibliography{references.bib}

\newpage
\appendix
\input{appendix}

\end{document}

%% file: preamble.tex
\usepackage[utf8]{inputenc}
\usepackage{amsmath, amsthm, amssymb, amsfonts, mathtools}
\usepackage{dsfont}
\usepackage{graphicx}
\usepackage[makeroom]{cancel}
\usepackage[shortlabels]{enumitem}
\usepackage{mathabx}
\usepackage{thmtools}
\usepackage{algorithm}
\usepackage{algorithmic}
\usepackage[capitalise,noabbrev]{cleveref}
\usepackage{xcolor}
\usepackage{array}
\usepackage{colortbl}
\usepackage[breakable]{tcolorbox}
\usepackage{accents}
\usepackage{mathrsfs}
\usepackage{tikz}
\usepackage{bm}
\usetikzlibrary{positioning, arrows.meta}

\usepackage[normalem]{ulem}

\usepackage{autonum}

\usepackage{caption}
\usepackage{subcaption}

\declaretheorem[name=Theorem,numberwithin=section]{theorem}
\declaretheorem[name=Proposition,numberwithin=section,sibling=theorem]{proposition}

\declaretheorem[name=Lemma,numberwithin=section,sibling=theorem]{lemma}
\declaretheorem[name=Corollary,numberwithin=section,sibling=theorem]{corollary}
\declaretheorem[name=Remark,numberwithin=section]{remark}
\newtheorem*{theorem*}{Theorem}
\declaretheorem[name=Definition,numberwithin=section,sibling=theorem]{definition}
\declaretheorem[name=Assumption,numberwithin=section,sibling=theorem]{assumption}
\declaretheorem[name=Example,numberwithin=section,sibling=theorem]{example}

\theoremstyle{plain}

\newcommand{\RR}{\mathbb R}
\newcommand{\EE}{\mathbb E}

\newcommand{\NN}{\mathbb N}

\newcommand{\Fcal}{\mathcal F}
\newcommand{\Gcal}{\mathcal G}

\newcommand{\Pcal}{\mathcal P}
\newcommand{\Qcal}{\mathcal Q}
\newcommand{\Rcal}{\mathcal R}

\newcommand{\Wcal}{\mathcal W}

\newcommand{\1}{\mathds{1}}

\newcommand{\seq}[4]{(#1)_{{#2}={#3}}^{#4}}

\newcommand{\infseqn}[1]{\seq{#1}{n}{1}{\infty}}

\newcommand{\dd}{\mathrm{d}}

\DeclareMathOperator*{\argmax}{argmax}

\newcommand{\iid}{\text{i.i.d.}}

\newcommand{\indep}{\perp \!\!\! \perp}
\newcommand{\eps}{\varepsilon}
\newcommand{\Var}{\mathrm{Var}}

\newcommand{\ubar}[1]{\underaccent{\bar}{#1}}

\newcommand{\supP}{{\sup_{P \in \Pcal}}}
\newcommand{\Pin}{{P \in \Pcal}}
\newcommand{\Qcalbar}{{\widebar \Qcal}}

\newcommand{\supQbar}{{\sup_{Q \in \Qcalbar}}}
\newcommand{\Qin}{{Q \in \Qcal}}
\newcommand{\Qinbar}{{Q \in \Qcalbar}}
\newcommand{\lambdain}{{\lambda \in \dsimp}}

\newcommand{\supkm}{{\sup_{k \geq m}}}
\newcommand{\tth}{{\text{th}}}
\newcommand{\mto}{{m \to \infty}}
\newcommand{\nto}{{n \to \infty}}
\newcommand{\alphato}{{\alpha \to 0^+}}

\newcommand{\brackone}{{(1)}}
\newcommand{\brackzero}{{(0)}}
\newcommand{\brackd}{{(d)}}
\newcommand{\bracktwo}{{(2)}}

\newcommand{\brackj}{{(j)}}
\newcommand{\brackalpha}{{(\alpha)}}
\newcommand{\brackP}{{(P)}}
\newcommand{\brackQ}{{(Q)}}
\newcommand{\Dleq}{{(D\leq)}}
\newcommand{\Deq}{{(D=)}}

\newcommand{\ONS}{{\mathrm{ONS}}}

\newcommand{\Qstar}{{(Q\star)}}

\newcommand{\TwoSTOracle}{\mathrm{2ST\star}}

\newcommand{\UP}{{\mathrm{UP}}}
\newcommand{\OJ}{{\mathrm{OJ23}}}
\newcommand{\CO}{{\mathrm{CO96}}}

\newcommand{\eproc}{{W}}
\newcommand{\eprocclass}{{\Wcal}}

\newcommand{\Pe}{\Pcal\text{-}e}
\newcommand{\bE}{{\mathbf{E}}}
\newcommand{\be}{{\mathbf{e}}}

\newcommand{\trans}{{\top}}

\newcommand{\Max}{{\mathrm{max}}}
\newcommand{\DKL}{{D_{\mathrm{KL}}}}

\newcommand{\Regret}{\Rcal}

\newcommand{\Qcalu}{{\Qcal^\circ}}
\newcommand{\Qinu}{{Q \in \Qcalu}}

\newcommand{\thetain}{{\theta \in \Theta}}
\newcommand{\GREE}{\mathrm{GREE}}

\newcommand{\dsimp}{{\triangle_d}}
\newcommand{\onesimp}{{\triangle_1}}

\newif\ifverbose %
\verbosefalse %

\newcommand{\verbose}[1]{\ifverbose\textcolor{gray}{\textit{verbose $\rightarrow$}}\quad #1\textcolor{gray}{\textit{non-verbose $\rightarrow$}}\quad \fi}

\newcommand{\simplex}{\triangle_d}
\renewcommand{\epsilon}{\varepsilon}

\newcommand{\pospart}[1]{\left [ #1 \right ]_+}

\newcommand{\Q}{Q}
\renewcommand{\P}{P}

\newcommand{\dP}{\dd \P}
\newcommand{\dQ}{\dd \Q}

\newcommand{\KL}{\mathrm{KL}}

\newcommand{\UI}{\mathrm{UI}}

%% file: abstract.tex
We consider the problem of sequential hypothesis testing using $e$-processes. For a rich class of composite testing problems---which include bounded mean testing, equal mean testing for bounded random tuples, and some key ingredients of two-sample and independence testing as special cases---we show that any $e$-process satisfying a certain sublinear regret bound is asymptotically and almost surely instance-log-optimal for a composite alternative. This is a strong notion of optimality that has not previously been established for the aforementioned problems, and we provide explicit test supermartingales and $e$-processes satisfying this notion in a more general case. Furthermore, we derive matching lower and upper bounds on the expected rejection time in the high-confidence regime for the resulting sequential tests in all of these cases. The proofs of these results make weak, algorithm-agnostic moment assumptions and rely on a proof technique involving the aforementioned regret and a family of numeraire portfolios. 
Finally, we discuss how all of these theorems hold in a distribution-uniform sense, a notion of log-optimality that is stronger still and seems to be new to the literature.

%% file: content.tex
\section{Introduction}\label{section:introduction}

Let $X\equiv \infseqn{X_n}$ be an \iid{} sequence of random variables on a filtered measurable space $(\Omega, \Fcal)$ and let $\Pcal$ and $\Qcal$ be collections of distributions---the ``null'' and ``alternative'' hypotheses, respectively. Broadly speaking, the goal of sequential hypothesis testing---a research program initiated by \citet{wald1945sequential,wald1947sequential}---is to construct a function (the ``sequential test''), $\phi_n^\brackalpha \equiv \phi^\brackalpha(X_1, \dots, X_n)$ for each $n \in \NN$, with the property that
\begin{equation}
  \forall \alpha \in (0, 1),\quad \supP P \left ( \exists n \in \NN : \phi_n^\brackalpha = 1 \right ) \leq \alpha \quad\text{or equivalently} \quad \supP P \left ( \tau_\alpha < \infty \right ) \leq \alpha,
\end{equation}
where $\tau_\alpha := \inf \{ n \in \NN : \phi_n^\brackalpha = 1\}$ is the stopping time representing the first sample size for which the test $\phi^\brackalpha$ rejects.
By far the most common way to construct such a test (and in some sense the only admissible way to do so \citep{ramdas2020admissible}) is to find a nonnegative stochastic process $\eproc \equiv \infseqn{\eproc_n}$ that is a $P$-supermartingale with mean $\EE_P [\eproc_1] \leq 1$ for every $P \in \Pcal$ and set $\phi_n^\brackalpha := \1 \{ \eproc_n \geq 1/\alpha \}$. Such supermartingales are called ``test supermartingales'' since
\begin{equation}
  \supP P \left ( \tau_\alpha < \infty \right ) = \supP P \left ( \exists n \in \NN : \eproc_n \geq 1/\alpha \right ) \leq \alpha,
\end{equation}
where the final inequality follows from Ville's inequality for nonnegative supermartingales \citep{ville1939etude}. This bound can be thought of as a time-uniform analogue of Markov's inequality and consequently yields a time-uniform notion of type-I error control.\footnote{As is common in the literature on anytime-valid inference, we use the terms ``time'' and ``sample size'' interchangeably.}

Our focus in the current paper is on the alternative rather than the null and in particular on the power for rejecting the null in favor of a (composite) alternative. When considering composite alternatives, the notion of ``power'' for test supermartingales is subtle, with at least two kinds of criteria appearing in the existing literature. For an alternative distribution $\Qin$, the first criterion is the asymptotic growth rate of $\eproc$ given by a $Q$-almost sure lower bound on
\begin{equation}\label{eq:intro-growth rate}
  \liminf_{n \to \infty} \frac{1}{n} \log( \eproc_n ),
\end{equation}
and the second is an upper bound on the expectation of the rejection time $\tau_\alpha := \inf \{n \in \NN : \eproc_n \geq 1/\alpha \}$ of the resulting test under $Q$:
\begin{equation}\label{eq:intro-stopping time}
  \EE_Q \left [ \tau_\alpha \right ].
\end{equation}
In the former case, larger values represent more power against $\Pcal$, while in the latter, larger values represent less power.

While there has been a flurry of work on sequential hypothesis testing and ``testing by betting'' in recent years~\citep[see, e.g.,][]{ramdas2024hypothesis,shafer2021testing},
less is known about useful particular classes of nonparametric test
supermartingales that exhibit optimal growth rates or expected rejection times
under the alternative, or about what conditions give rise to these optimality
properties. Our main contribution is to demonstrate the existence of some such classes, establishing optimality in a strong sense: adaptive,
asymptotic, and almost sure instance-log-optimality. In particular, the setting we study is the following. Fix an integer $d \in \NN$.
Given $(d+1)$-dimensional random vectors $\infseqn{\lambda_n}$ taking values in the $d$-dimensional simplex $\simplex$ and $\infseqn{\bE_n} \equiv \infseqn{(E_n^\brackzero,\dots, E_n^\brackd)}$ taking values in $\RR_+^{d+1}$, the class we study is made up of test supermartingales that take the following form:
\begin{equation}\label{eq:intro-general-test-sm}
  \eproc_n \equiv \eproc_n(\lambda_1, \dots, \lambda_n) := \prod_{i=1}^n \lambda_i^\trans \bE_i,
\end{equation}
where $\infseqn{\lambda_n}$ is a predictable process  (informally, a
``portfolio'' or a ``betting strategy''; see \cref{section:preliminaries}).
Under every distribution considered throughout this paper, the tuples $\infseqn{\bE_n}$ are assumed to be independent and identically distributed (\iid{}) but their entries can be arbitrarily dependent unless otherwise specified.
We consider the null hypothesis given by
\begin{equation}
  \widebar \Pcal := \{ P : \text{$\EE_P[E_1^\brackj] \leq 1$ for each $j \in \{0, \dots, d\}$} \}.
\end{equation}
However, we primarily study the subset $\Pcal \subseteq \widebar \Pcal$ for which $E_1^\brackzero = 1$ deterministically, that is
\begin{equation}\label{eq:null}
 \Pcal := \{ P : \text{$E_1^\brackzero = 1$ and $\EE_\P [E_1^\brackj] \leq 1$ for each $j \in \{1,\dots, d\}$} \}.
\end{equation}
The reason for this distinction is that we study lower and upper bounds for the aforementioned optimality criteria within the smaller class $\Pcal$ but it is nevertheless of interest to derive upper bounds for the class $\widebar \Pcal$ since there exist widely studied test supermartingales in the literature for which the trivial $e$-value of $1$ is not included.
It is routine to check that $\eproc_n$ forms a test supermartingale.
The validity of $W_n$ does not require that the vectors $\infseqn{\bE_n}$ are \iid{}, but it plays an important role in the definition and proof of the optimality guarantees of $W_n$ under the alternative that we state in the sections to come.

Readers familiar with nonparametric sequential tests may  recognize a special case of \eqref{eq:intro-general-test-sm} where $d = 1$ and $E_n^\brackzero = 1$ and $E_n \equiv E_n^\brackone$ for each $n \in \NN$ so that
\begin{equation}\label{eq:intro-general-test-sm-special-case}
  \eproc_n = \prod_{i=1}^n \left ( 1 + \lambda_i\cdot (E_i - 1) \right ).
\end{equation}
See, for example, \citet{wang2022backtesting} and \citet[\S 7]{ramdas2024hypothesis}.

The supermartingales in \eqref{eq:intro-general-test-sm-special-case} (and hence \eqref{eq:intro-general-test-sm}) arise both implicitly and explicitly in a wide range of sequential testing  problems. These include one- and two-sided tests for the mean of a bounded random variable \citep{waudby2020estimating,orabona2021tight,ryu2024confidence}, two-sample testing \citep{shekhar2023nonparametric,podkopaev2023sequentialTwoSample}, marginal independence testing \citep{podkopaev2023sequentialTwoSample},
auditing the fairness of deployed machine learning models \citep{chugg2023auditing}, auditing elections \citep{stark2020sets,waudby2021rilacs,stark2023alpha}, and backtesting to assess risk measures \citep{wang2022backtesting}. We expand on a subset of these problems in \cref{section:corollaries} and make the exact connections explicit. In all of the cases considered in \cref{section:corollaries}, the optimality results that follow for processes of the form \eqref{eq:intro-general-test-sm} lead to new almost sure log-optimality results that had not yet been established for these special cases. The spirit of these results in fact holds for a more abstract class of test supermartingales discussed in \cref{section:generalizations} (encompassing sub-Gaussian mean testing and universal inference \citep{wasserman2020universal}) but we focus on the case of \eqref{eq:intro-general-test-sm} as it balances generality and concreteness while being interpretable as the wealth from a dynamically rebalanced \emph{portfolio}; see \cref{section:preliminaries}.

We can now formulate our primary question: Is there a systematic way to choose $\infseqn{\lambda_n}$ such that the test supermartingales in \eqref{eq:intro-general-test-sm} enjoy growth-rate-optimal and/or expected-rejection-time-optimal properties, and if so, is there a fundamental quantity that characterizes their optimality? We answer this question in the affirmative, where the aforementioned fundamental quantity is the maximum expected logarithmic increment of \eqref{eq:intro-general-test-sm} under some alternative $Q \in \Qcal$; i.e., $\ell_Q^\star := \ell_Q(\lambda_Q)$ where
\begin{equation}
  \lambda_Q \in \argmax_{\lambda \in \dsimp} \ell_Q(\lambda) \quad\text{and}\quad\ell_Q(\lambda) := \EE_Q \left [ \log \left ( \lambda^\trans \bE_1 \right ) \right ].
\end{equation}
Following \citet[\S 15]{cover1999elements} we refer to $\lambda_Q$ as a ``log-optimal'' portfolio or ``strategy'' under $Q$, and as we discuss in \cref{section:generalizations}, this is precisely the numeraire portfolio of \citet{long1990numeraire}; see also \citet{larsson2024numeraire}. The approach of maximizing the expected log-returns is also known as the Kelly criterion \citep{kelly1956new,latane1959criteria,breiman1961optimal}.

With this preamble and notation in mind, we present an informal summary of the main results (the formal versions of which can be found in \cref{section:asymptotic-equivalence,section:stopping time}).

\begin{tcolorbox}[boxrule=0pt, frame empty, breakable]
  \begin{theorem*}[An informal summary of the main results of \cref{section:asymptotic-equivalence,section:stopping time}]
    Let $\eproc \equiv \infseqn{\eproc_n}$ be a test $\Pcal$-supermartingale
    taking the form in \eqref{eq:intro-general-test-sm}. Let $\Qcal$ be any
    alternative hypothesis for which the vector of $e$-values
    $\infseqn{\bE_n}$ is \iid{}, such as
    \begin{equation}
      \Qcal := \left \{ P : \text{$\EE_P[E_1^\brackj] > 1$ for some $j \in [0:d]$} \right \}.
    \end{equation}
    If $\infseqn{\lambda_n}$ is chosen according to any algorithm with a sublinear portfolio regret (e.g., a universal portfolio \citep{cover1991universal}; see \cref{section:preliminaries} for a definition), then
    \begin{enumerate}[label = (\roman*)]
    \item $\eproc$ has an optimal asymptotic growth rate; i.e., for any $\Qin$,
    \begin{equation}
      \liminf_{\nto } \left ( \frac{1}{n} \log (\eproc_n) - \frac{1}{n} \log (\eproc_n') \right ) \geq 0 \quad\text{$Q$-almost surely,}
    \end{equation}
    where $\eproc_n'$ is any other $\Pcal$-supermartingale.
    \item The asymptotic growth rate of $\eproc$ is given by
    \begin{equation}
      \lim_{\nto}\frac{1}{n} \log (\eproc_n ) = \max_{\lambdain} \ell_Q(\lambda)\quad\text{$Q$-almost surely.}
    \end{equation}
    \item If $\ell_Q^\star > 0$ and the log-wealth increments have a finite $(2+\delta)^\tth$ moment for any $\delta > 0$, the rejection time $\tau_\alpha := \inf \{n \in \NN : \eproc_n \geq 1/\alpha \}$ can be bounded in expectation for small $\alpha \in (0, 1)$ as:
    \begin{equation}
      \EE_Q \left [ \tau_\alpha \right ] \lesssim \frac{\log(1/\alpha)}{\max_{\lambdain} \ell_Q(\lambda)}. 
    \end{equation}
    \item The right-hand side of the inequality in $(iii)$ is unimprovable: for any stopping rule $\tau_\alpha'$ satisfying $\sup_\Pin \P(\tau_\alpha' < \infty)\leq \alpha$, we have
    \begin{equation}
      \EE_Q \left [ \tau_\alpha' \right ] \geq \frac{\log(1/\alpha)}{\max_{\lambdain}\ell_Q(\lambda)}.
    \end{equation}
    \end{enumerate}
  \end{theorem*}
\end{tcolorbox}
In several prior works, betting strategies based on regret bounds have been used to show that the corresponding test $\Pcal$-supermartingales diverge to $\infty$ with $Q$-probability one for a collection of alternative hypotheses $Q \in \Qcal$ (and we provide an extended discussion of such work in \cref{section:related work}). As such, the stopping times of the resulting tests are $Q$-almost surely finite:
\begin{equation}\label{eq:intro finite stopping time}
  \Q (\tau_\alpha < \infty) = 1,
\end{equation}
and some of these works have also derived lower bounds on the asymptotic almost sure growth rate; see for example papers by \citet{shekhar2023nonparametric}, \citet{podkopaev2023sequentialTwoSample}, \citet{chugg2023auditing}, and \citet{chen2025optimistic}. However, as we discuss in more detail in \cref{section:corollaries}, those derived growth rates are not optimal in general.
In \citet{wang2022backtesting}, the authors derive optimal lower bounds on asymptotic growth rates in $L_1$ rather than with probability one. We expand on their approach and the guarantees thereof in \cref{remark:gree}.

There is a rich literature deriving upper bounds on expected rejection times, both in the context of sequential hypothesis testing (see, for example, the works of \citet{wald1947sequential,robbins1974expected,shekhar2023nonparametric,chugg2023auditing}) but also in the stochastic multi-armed bandit literature, especially those works studying the ``best-arm identification'' problem (see \citet{kaufmann2016complexity,kaufmann2021mixture,agrawal2020optimal,agrawal2021optimal,honda2010asymptotically,honda2011asymptotically,honda2015non}). However, neither for the setting considered in~\eqref{eq:intro-general-test-sm} nor for the special cases of \cref{section:corollaries} has it been shown that a betting strategy $\infseqn{\lambda_n}$ can attain matching lower and upper bounds in the sense of the above informal theorem. In the case of \citet{honda2010asymptotically}, however, we can use the power of hindsight to justify optimality of their approach for the special case of \cref{section:bounded one sided}. We briefly discuss in \cref{section:testing problems to keep in mind} why one particular but common approach has not attained these matching lower and upper bounds. The contribution of the current paper is to close these gaps and provide a unified perspective on optimality for these problems.  Moreover, our theory is an operational one---we present explicit (and computationally tractable) test supermartingales, $e$-processes, and sequential hypothesis tests for the settings we study. Finally, we also provide distribution-uniform generalizations of our results in \cref{section:uniformity} which may be of independent interest---since a distribution-uniform notion of ``power'' seems to have been missing in the sequential testing literature.
\begin{figure}[h]
  \centering
  \includegraphics[width=0.49\linewidth]{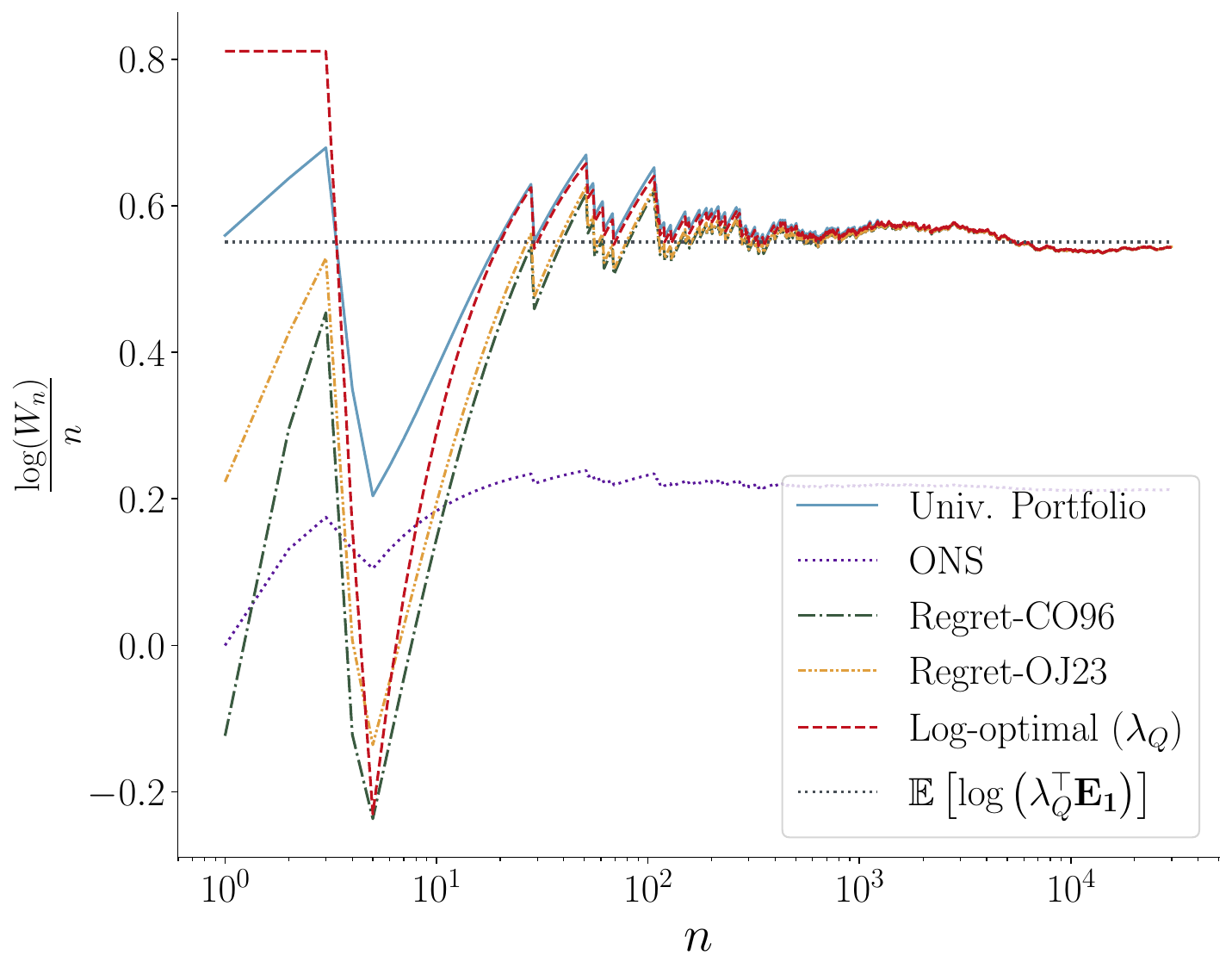}
\includegraphics[width=0.49\linewidth]{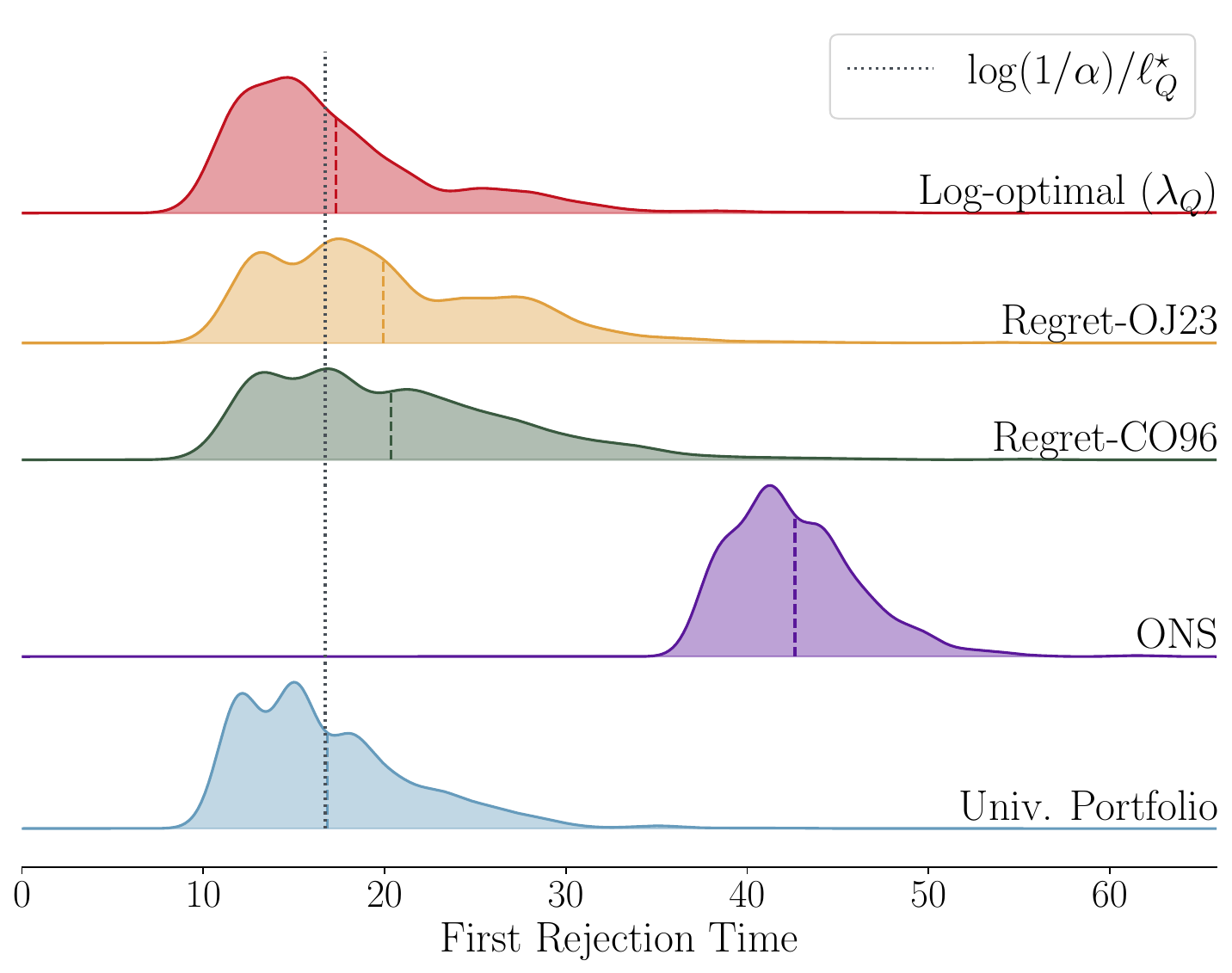}
\caption{Empirical growth rates (left) and distributions of
rejection times (right) for various $e$-processes (see
\cref{section:preliminaries} for a precise definition). As will be discussed in
\cref{corollary:adaptive optimality of UP and OJ}, the $e$-processes labeled
``Univ.~Portfolio,'' ``Regret-CO96,'' and ``Regret-OJ23'' all satisfy a
sublinear portfolio regret bound. For this reason, they all have growth rates
converging to $\ell_Q^\star$ and expected rejection times close to
$\log(1/\alpha) / \ell_Q^\star$. In particular, Online Newton Step (ONS) is a
commonly studied strategy in the literature for deriving growth rate and
rejection time guarantees but it does not satisfy a sublinear portfolio regret bound, and
can consequently be suboptimal from both of the perspectives discussed in the
introduction. More details can be found in
\cref{section:asymptotic-equivalence,section:stopping time}.}
  \label{fig:intro}
\end{figure}

\subsection{Preliminaries: notation and background}\label{section:preliminaries}
In the previous discussions, we glossed over some technical details surrounding filtrations, supermartingales, and stopping times for the sake of exposition. Let us now make these matters more formal for the results to come and introduce a few more crucial concepts.

\paragraph*{(i) Notation and formalities} Throughout, we will work in reference to a filtered measurable space $(\Omega, \Fcal)$ where $\Fcal \equiv (\Fcal_n)_{n=0}^\infty$ will be some discrete time filtration such that $\Fcal_0$ is the trivial sigma-algebra. For concreteness, one can think of $\infseqn{X_n} \equiv (X_1, X_2, \dots)$ as being a sequence of random objects (typically real-valued) and $\Fcal_n = \sigma(X_1, \dots, X_n)$ as being the sigma-algebra generated by those random objects up until some time $n \in \NN$. We will additionally consider two families of distributions, $\Pcal$ and $\Qcal$, which will be used to denote the null and alternative hypotheses respectively. In particular, if we fix some $P \in \Pcal$ in the null, we will be implicitly working with the probability space $(\Omega, \Fcal, P)$ and similarly for some $Q \in \Qcal$ in the alternative. We will write the expectation of a random variable $Y$ under some probability measure $P$ as $\EE_P \left [ Y \right ] := \int Y \dd P$.
Given the centrality of supermartingales and martingales in the discussions to come, we define them formally using the established notation here.
\begin{definition}[Test (super)martingales]\label{definition:supermartingale}
  Let $\Pin$. A stochastic process $\infseqn{\widebar W_n}$ on $(\Omega, \Fcal)$ is said
  to be a $P$-supermartingale if it is adapted to $\Fcal$---meaning that $M_n$
  is $\Fcal_n$-measurable for each $n \in \NN$---and if
\begin{equation}\label{eq:supermart property}
  \forall n \in \NN,\quad\EE_P \left [ \widebar W_n \mid \Fcal_{n-1} \right ]
  \leq \widebar W_{n-1}\quad\text{$P$-almost surely.}
\end{equation}
We say that this process is a $\Pcal$-supermartingale if \eqref{eq:supermart
property} holds for all $P \in \Pcal$. Additionally, the process
$\infseqn{\widebar W_n}$
is said to be a $P$-martingale if \eqref{eq:supermart property} holds but with
the inequality $(\leq)$ replaced by an equality $(=)$ and it is said to be a
$\Pcal$-martingale if it is a $P$-martingale for every $\Pin$.

Finally, we say that $\infseqn{\widebar W_n}$ is a \emph{test $P$-supermartingale}
(resp.~\emph{test $P$-martingale}) if $\widebar W_n \geq 0$ with $P$-probability one and
$\EE_P[ \widebar W_1 ] \leq 1$ (resp.~$\EE_P [ \widebar W_1 ] = 1$).
\end{definition}

\paragraph*{(ii) $E$-processes} All of the results to come can be presented in terms of so-called \emph{$e$-processes} which serve as a generalization of test supermartingales for the purposes of composite sequential testing; see \citep{ramdas2021testing} and \citep[\S 7]{ramdas2024hypothesis}.
\begin{definition}[$E$-processes]\label{definition:e-process}
  Let $\Pcal$ be a collection of distributions (to be thought of as a null
  hypothesis). A nonnegative adapted process $\infseqn{\eproc_n}$ is said to be
  a $\Pcal$-$e$-process if it is $P$-almost surely upper bounded by a test
  $P$-supermartingale $\infseqn{\widebar W_n^\brackP}$ for each $P \in \Pcal$:
\begin{equation}
  \forall n \in \NN, P \in \Pcal\quad \eproc_n\leq \widebar W_n^\brackP
  \quad\text{$P$-almost surely.}
\end{equation} 
\end{definition}
Note that the upper-bounding test supermartingale $\infseqn{\widebar
W_n^\brackP}$ in \cref{definition:e-process} can depend on the distribution $P
\in \Pcal$, while the $e$-process $\infseqn{\eproc_n}$ itself does not; see
\citet{ramdas2021testing} for some work dedicated to this subtlety.
Importantly for the purposes of hypothesis testing, $e$-processes still satisfy
Ville's inequality since for any $\Pin$, $P ( \sup_n \eproc_n \geq 1/\alpha
) \leq P ( \sup_n \widebar W_n^{(P)} \geq 1/\alpha ) \leq \alpha$.
In particular, the process $\infseqn{\phi_n}$ given by $\phi_n := \1 \{ \eproc_n
\geq 1/\alpha \}$ is a level-$\alpha$ sequential test for the null $\Pcal$.
Throughout the paper, we often make statements about the class of
$\Pcal$-$e$-processes given by
\begin{equation}\label{eq:prelim-eprocclass}
  \eprocclass := \{ \eproc : \eproc_n \leq \widebar \eproc_n \text{ and $\widebar \eproc_n$ is a test $\Pcal$-supermartingale of the form \eqref{eq:intro-general-test-sm}} \},
\end{equation}
and this is in an effort to effectively study test supermartingales of the form \eqref{eq:intro-general-test-sm} but also allow their lower-bounding $e$-processes for the sake of generality. Explicit $\Pe$-processes that are not supermartingales will be studied in \cref{corollary:adaptive optimality of UP and OJ}. 

\paragraph*{(iii) Stopping times} A stopping time $\tau$ is a $\NN\cup \{ \infty \}$-valued random variable on $(\Omega, \Fcal)$ such that $\{ \tau = n \} \in \Fcal_n$. We are interested in stopping times $\tau_\alpha$ of the form
\begin{equation}
  \tau_\alpha := \inf \{ n \in \NN : \eproc_n \geq 1/\alpha \},
\end{equation}
for some $\Pcal$-$e$-process $\eproc$, where the above should be interpreted as the first time that the test $\phi_n := \1 \{ \eproc_n \geq 1/\alpha \}$ rejects, and we take $\inf \emptyset = \infty$ as a convention. When using the notation of stopping times, Ville's inequality can be thought of as guaranteeing that the aforementioned test will \emph{never} reject except with small probability:
\begin{equation}\label{eq:prelim-stopping-time}
  \supP P \left ( \tau_\alpha < \infty \right ) \leq \alpha.
\end{equation}
In many cases, it will hold that $\inf_{Q \in \Qcal} \Q(\tau_\alpha < \infty) = 1$---i.e.~that the test $\phi_n$ has eventual power one under the alternative $\Qcal$---but we leave those discussions for specific contexts.
\paragraph*{(iv) Deterministic portfolio regret} The discussion thus far has focused on probabilistic tools from the theory of stochastic processes.  An alternative perspective on sequences of data comes from the fields of online learning,  game theory, and investment theory, where one defines a notion of ``regret'' of an online procedure by comparing its performance to that of an oracle that can see the entire sequence in hindsight.
Regret is typically viewed as a deterministic quantity and there is much work dedicated to deriving algorithms whose regret is bounded for each individual sequence, rather than in expectation or with high probability~\citep{cesa2006prediction,hazan2016introduction}.  Particularly relevant to this paper is a notion of regret that arises in portfolio theory \citep{cover1991universal,cover1996universal}. Since we will refer to this notion frequently by the name ``portfolio regret,'' we provide a definition here.
\begin{definition}[Portfolio regret]\label{definition:portfolio regret}
  Let $\infseqn{\be_n}$ be a deterministic sequence of $(d+1)$-dimensional tuples lying in $[0,\infty)^{d+1}$. Let $\infseqn{W_n}$ be a real sequence. For any $n \in \NN$, we say that $W_n$ has a \emph{portfolio regret} of $\Regret_n$ given by
  \begin{equation}\label{eq:prelim-regret}
    \Rcal_n \equiv \Rcal(W_n) := \max_{\lambda \in \dsimp} \sum_{i=1}^n \log \left ( \lambda^\trans \be_i \right ) - \log(W_n).
  \end{equation}
\end{definition}
The reason for the term ``portfolio regret'' stems from the language used by \citet{cover1991universal} and \citet{cover1996universal} in their work on universal portfolio theory and because the product $\prod_{i=1}^n \lambda_i^\trans \be_i$ can be viewed as the wealth accumulated by an investor using a portfolio of $\lambda_n \in \dsimp$ where they have placed a $\lambda_n^\brackj$-fraction of their wealth on stock $\be_n^\brackj$ for each $j \in [0:d]$ at time $n \in \NN$.
The expression in \eqref{eq:prelim-regret} should be interpreted as the difference between the log-wealth obtained by an oracle (that implements the best-in-hindsight constant rebalanced portfolio for processes of the form \eqref{eq:intro-general-test-sm}) versus the logarithm of some other process $\eproc$. One can think of $W_n$ as being  given by
\begin{equation}\label{eq:prelim-wealth}
  W_n = \prod_{i=1}^n \lambda_i^\trans \be_i
\end{equation}
for a $\dsimp$-valued sequence $\infseqn{\lambda_n}$ where $\lambda_n$ depends only on the stocks up until time $n-1$.\footnote{In the information theory literature, sequences $\infseqn{\lambda_n}$ such that $\lambda_n$ only depends on the first $n-1$ observations are often called ``causal'' or ``nonanticipating'' \citep{cover1991universal,cover1996universal,cover1999elements} but we use the term ``predictable'' as this is more common in probability and statistics.}
As one possible algorithm satisfying a sublinear portfolio regret bound, consider Cover's universal portfolio strategy \citep{cover1991universal}, defined at time $n \in \NN$ as a mixture over the possible wealth values at time $n-1$:
\begin{equation}\label{eq:UP-lambda}
  \lambda_n^\UP := \frac{\int_{\lambda \in \simplex} \lambda \eproc_{n-1}(\lambda) \dd F(\lambda)}{\int_{\lambda\in\simplex} \eproc_{n-1}(\lambda) \dd F(\lambda) }.
\end{equation}
Taking $F(\cdot)$ to be the $(d+1)$-dimensional Dirichlet$(1/2, \dots, 1/2)$ distribution and plugging $\infseqn{\lambda_n^\UP}$ into \eqref{eq:prelim-wealth} to obtain $\infseqn{W_n^\UP}$,
 \citet[Theorem 3]{cover1996universal} establish that for \emph{any} $[0,\infty)^{d+1}$-valued sequence $\infseqn{\be_n}$, for all $n \geq 2$,
\begin{equation}\label{eq:prelim-regret-ub}
  \Regret \left ( W_n^\UP \right ) \leq \frac{d\log ( n + 1)}{2} + \log 2.
\end{equation}
Hence if the sequence $\infseqn{\be_n} \equiv\infseqn{\bE_n}$ is in fact a stochastic process, then so is $\infseqn{\Regret_n}$, and the latter quantity can be bounded for every $\omega \in \Omega$ and hence $P$-almost surely for any probability distribution $P$. Importantly, the \iid{} assumption we impose on $\infseqn{\bE_n}$ is not required to define the notion of portfolio regret nor to attain sublinear portfolio regret. Nonetheless, the \iid{} assumption is required to state and prove the optimality guarantees under the alternative that we present shortly.
Finally, we note that when $d = 1$, a regret bound sharper than \eqref{eq:prelim-regret-ub} was obtained by \citet{orabona2021tight}; we make use of this bound in \cref{corollary:adaptive optimality of UP and OJ}.

\subsection{A concrete testing problem to keep in mind: bounded means}\label{section:testing problems to keep in mind}

Consider the problem of testing whether the mean of a bounded random variable is equal to some $\mu_0 \in [0, 1]$.
That is, let $\infseqn{X_n}$ be \iid{} random variables supported on the unit interval $[0, 1]$ with mean $\mu_P = \EE_P[X_1]$ under the distribution $P$. We are considering the unit interval $[0, 1]$ without loss of generality since one can always take a random variable bounded on $[a,b]$ for $a < b$ and transform it to the unit interval via $x \mapsto (x -a)/(b-a)$. Suppose that we are interested in testing the following null hypothesis $\Pcal^=$ versus the alternative $\Qcal^\neq$:
\begin{equation}\label{eq:bounded-means-null-one-sided}
  \Pcal^= = \{ P : \EE_P[X_1] = \mu_0 \} \quad\text{versus}\quad \Qcal^\neq = \{ P : \EE_P[X_1] \neq \mu_0\}
\end{equation}
for some $\mu_0 \in [0, 1]$.
Following some prior works that study this problem such as \citet{hendriks2021test,waudby2020estimating,orabona2021tight,ryu2024confidence}, as well as others in more implicit forms and in different contexts \citep{shekhar2023nonparametric,chugg2023auditing,podkopaev2023sequentialKernelized,podkopaev2023sequentialTwoSample,ryu2025improved}, consider the $\Pcal^=$-test martingale given by
\begin{equation}\label{eq:bounded-means-martingale-two-sided}
  \eproc_n^= := \prod_{i=1}^n \left (1 + \gamma_i \cdot (X_i - \mu_0) \right ),
\end{equation}
where $\infseqn{\gamma_n}$ is any $[-1/(1-\mu_0), 1/\mu_0]$-valued predictable sequence. Indeed, \citet[Proposition 2]{waudby2020estimating} show that every test $\Pcal^=$-martingale must take the form in \eqref{eq:bounded-means-martingale-two-sided}.
Several strategies for choosing $\infseqn{\gamma_n}$ have been proposed in prior works \citep{waudby2020estimating,orabona2021tight,chugg2023auditing,ryu2024confidence,shekhar2023near,shekhar2023nonparametric,podkopaev2023sequentialKernelized,podkopaev2023sequentialTwoSample,cho2024peeking} and while many of them exhibit excellent empirical performance, they have not yet been shown to explicitly satisfy the optimality desiderata appearing in the introduction. Moreover, a commonly studied strategy for the purposes of obtaining growth-rate- and expected-rejection-time-optimality guarantees in bounded mean testing is given by the so-called Online Newton Step (ONS) algorithm \citep{hazan2007logarithmic,cutkosky2018black} defined with $\gamma_1^\ONS := 0$ and then recursively as
\begin{align}
  \gamma_n^\ONS := \left ( \left ( \gamma_{n-1}^\ONS - \frac{2}{2-\log(3)} \cdot \frac{Z_{n-1}}{1 + \sum_{i=1}^{n-1} Z_i^2} \right ) \land \frac{1}{2} \right ) \lor -\frac{1}{2},\\
  \text{where}~~Z_i := \frac{-(X_i - \mu_0)}{1+\gamma_i^\ONS \cdot (X_i - \mu_0)}.
\end{align}
We will not make use of the exact form of $\gamma_n^\ONS$ in this paper, but it is worth noting that $\gamma_n^\ONS \in [-1/2, 1/2]$ for every $n \in \NN$ and hence if the maximizer of the expected log-wealth increment under $Q$ lies outside of this range, ONS cannot adapt to this fact, leading to suboptimal growth rates and expected rejection times; see \cref{fig:intro}.
Comparing to processes of the form \eqref{eq:intro-general-test-sm} more explicitly, it can be verified that \eqref{eq:bounded-means-martingale-two-sided} is a special case of \eqref{eq:intro-general-test-sm} with $d = 1$ for the $e$-values $E_n^\brackzero := (1-X_n) / (1-\mu_0)$ and $E_n^\brackone := X_n / \mu_0$ and by defining
\begin{equation}
  \lambda_n := \mu_0 + \gamma_n \mu_0(1-\mu_0)
\end{equation}
for each $n$.
In particular, the range $\gamma_n \in [-1/2, 1/2]$ constrains $\lambda_n$ to lie in $\mu_0 \pm \mu_0(1-\mu_0) /2$ which is always a strict subset of $[0, 1]$. Meanwhile, any strategy satisfying sublinear portfolio regret \eqref{eq:prelim-regret} for arbitrary sequences must have access to all of $[0, 1]$; hence the aforementioned suboptimality of ONS when compared to those with sublinear portfolio regret.

Note that \citet{waudby2020estimating} used the symbol $\lambda$ in place of $\gamma$ in \eqref{eq:bounded-means-martingale-two-sided}, but to keep notation consistent with the textbook chapter of \citet[\S 7]{ramdas2024hypothesis}, we reserve $\lambda$ for betting strategies used in test supermartingales of the form \eqref{eq:intro-general-test-sm}. While all of the results to follow will be stated in the form of \eqref{eq:intro-general-test-sm}, the bounded mean testing problem can be kept in mind as a concrete and nontrivial special case that has been of interest in the literature.

\section{Asymptotic and almost sure instance-log-optimality}\label{section:asymptotic-equivalence}

As alluded to in the introduction, one often thinks of ``powerful''
$e$-processes as those that have a fast growth rate under alternative hypotheses
$\Qcal$, and we pursue this perspective throughout this section. We are
interested in making statements about $e$-processes and test supermartingales
that we can explicitly construct, especially insofar as they compare to an
``oracle'' $e$-process that we cannot construct but that is optimal in a certain
sense. The following two definitions provide some of the requisite language to
make such comparisons formally. 
\begin{definition}[Asymptotic and almost sure instance-log-optimality and equivalence]
  \label{definition:log-optimality}
  Let $\eprocclass$ be a collection of $\Pcal$-$e$-processes. We say that $\eproc^\star = \infseqn{\eproc_n^\star} \in \eprocclass$ is $Q$-almost surely \emph{log-optimal} within $\eprocclass$ if for any other $\eproc \in \eprocclass$ it holds that
  \begin{equation}\label{eq:log-optimal-def}
    \liminf_{n \to \infty} \frac{1}{n} \left (  \log \eproc_n^\star  - \log \eproc_n  \right ) \geq 0\quad\text{$Q$-almost surely,}
  \end{equation}
  and we say that $\eproc^\star$ is \emph{$\Qcal$-instance-log-optimal} if \eqref{eq:log-optimal-def} holds for every $Q \in \Qcal$.
  Furthermore, we say that $\eproc, \eproc' \in \eprocclass$ are \emph{$Q$-almost surely asymptotically equivalent} with a rate of $\infseqn{a_n/n}$ if
  \begin{equation}\label{eq:equivalent-def}
      \frac{1}{n} \left ( \log \eproc_n -  \log \eproc'_n \right ) 
      = o(a_n/n)\quad\text{$Q$-almost surely,}
  \end{equation}
  where $\infseqn{a_n}$ is a nonnegative, monotonically nondecreasing sequence that diverges to infinity, a condition that we write as $a_n \nearrow \infty$.
\end{definition}
In words, a process $\eproc^\star$ satisfying \eqref{eq:log-optimal-def} is one that diverges to $\infty$ no slower than any other in $\eprocclass$ with $Q$-probability one. Moreover, if $\eproc^\star$ consists of \iid{} multiplicands---i.e., $\eproc_n^\star := \prod_{i=1}^n E_i$ for \iid{} $\infseqn{E_n}$---then we have that $\frac{1}{n} \log (\eproc_n^\star) \to \EE_Q [\log E_1]$ $Q$-almost surely by the strong law of large numbers. Consequently, the process $\eproc_n^\star$ can be written in the following exponential form:
\begin{equation}
  \eproc_n^\star = \exp \left \{ n \EE_Q[\log E_1] + o(n) \right \}\quad\text{$Q$-almost surely.}
\end{equation}
Notice that if $\eproc$ is any $e$-process in $\eprocclass$ that is $Q$-almost surely asymptotically equivalent to $\eproc^\star$ in the sense of \eqref{eq:equivalent-def}, then $W$ is itself $Q$-almost surely log-optimal. The notions of asymptotic equivalence and optimality defined above are similar to those of \citet[Definition 3]{wang2022backtesting} in the context of the test supermartingale in \eqref{eq:intro-general-test-sm-special-case} but those authors consider convergence in $L_1(Q)$ rather than with $Q$-probability one. The authors do show that some explicit betting strategies satisfy this notion of $L_1(Q)$-log-optimality and we discuss one of them further in \cref{section:related work}. While we focus entirely on almost sure notions of log-optimality for the remainder of the main text, \cref{section:l1} considers analogues of our results under $L_1$ convergence. These are relegated to the appendix for brevity.
\begin{remark}
  In the non-sequential context with a composite null $\Pcal$ but a \emph{point} alternative $Q$, the phrase ``log-optimality'' is often used to refer to a single $\Pcal$-$e$-value $E^\star := \argmax_{E \in \mathcal E} \EE_Q[\log E]$ where $\mathcal E$ consists of all possible $\Pcal$-$e$-values. Such an $e$-value is referred to as the \emph{numeraire $e$-value} and this is given a detailed treatment in \citet{larsson2024numeraire}; see also \citet[\S 6]{ramdas2024hypothesis}, as well as \citet{lardy2024reverse} and \citet{grunwald2019safe}. The criterion of growth-rate optimality was considered for point alternatives (GRO) and in a worst-case sense (GROW) by \citet{grunwald2019safe}. To distinguish the aforementioned use of the phrase from that in \cref{definition:log-optimality}, the latter articulates a notion of log-optimality that is (a) sequential (through the use of $e$-processes rather than $e$-values), (b) almost sure (in the form of its asymptotic convergence or dominance over other processes), (c) instance-dependent (since optimality is defined with respect to any unknown $\Qin$ that must be adapted to within a composite alternative $\Qcal$), and (d) within a class $\eprocclass$ which is typically taken to be a strict subset of all possible $\Pcal$-$e$-processes.
  Finally, the work of \citet{larsson2024numeraire} focuses on showing the \emph{existence and uniqueness} of such a numeraire $e$-value among $\mathcal E$ as well as demonstrating optimality properties thereof, while we are interested in providing sufficient conditions for explicit algorithms to attain $\Qcal$-instance-log-optimality as in \cref{definition:log-optimality}.
\end{remark}
Let us now center our attention on the null hypothesis $\Pcal$ describing a $(d+1)$-vector of $e$-values.

\subsection{Log-optimality via sublinear portfolio regret}
We first write a fundamental upper bound on the almost sure asymptotic growth rate of any $\Pcal$-$e$-process, providing a benchmark against which other $\Pcal$-$e$-processes shall be compared.

\begin{proposition}\label{proposition:growth rate upper bound}
  Let $\infseqn{W_n}$ be any $\Pe$-process. Let $Q$ be an alternative distribution. Then, the almost sure asymptotic growth rate of $W$ cannot exceed $\ell_Q^\star \equiv \max_{\lambda \in \dsimp} \EE_Q \left[\log\left(\lambda^\trans \bE_1 \right)\right]$. That is,
  \begin{equation}
    \limsup_{n \to \infty} \frac{1}{n} \log(W_n) \leq \ell_Q^\star 
  \end{equation}
  with $\Q$-probability one.
\end{proposition}
The proof of \cref{proposition:growth rate upper bound} is illustrative so we provide it here, leaving technical lemmas to the supplementary material.
\begin{proof}[Proof of \cref{proposition:growth rate upper bound}]
  Let $\P\in\Pcal$ and $\Q \in \Qcal$. We begin by showing that
  \begin{equation}
    \Q \left ( \limsup_{n \to \infty} \frac{1}{n} \log(W_n) \leq \KL(\Q \| \P) \right ) = 1.
  \end{equation}
  Indeed, assume that $\KL(\Q \| \P) < \infty$ otherwise the result is trivial.
  Fix
  $\eps > 0$ and $m \in \NN$. For each $n \in \NN$, define $L_n := \prod_{i=1}^n \dQ(\bE_i) / \dP(\bE_i)$. We have that
  \begin{align}
    \Q \left ( \supkm \left \{ \frac{1}{k}\log \left ( W_k \right ) - \frac{1}{k} \log(L_k) \right \} \geq \eps \right ) &= \Q \left ( \supkm \left \{ W_k / L_k \right \} \geq \exp \{ k\eps \} \right )\\
                                                                                                                         &\leq \sum_{k=m}^\infty \EE_\Q \left [ W_k / L_k \right ] \exp \left \{ -k\eps \right \}\\
                                                                                                                         &\leq \sum_{k=m}^\infty \exp \left \{ -k\eps \right \},
  \end{align}
  where the final inequality follows from \cref{lemma:LR is log-opt}. Clearly, the final series vanishes as $m \to \infty$. Since $\eps> 0$ was arbitrary,
    $\limsup_{n\to\infty} \frac{1}{n} \log(W_n / L_n) \leq 0$
  with $\Q$-probability one. Since $\Q(\frac{1}{n} \log(L_n) \to \KL(\Q \| \P)) = 1$ we have that
  \begin{equation}
    \limsup_{n\to \infty} \frac{1}{n} \log \left ( W_n \right ) \leq \KL(\Q \| \P)
  \end{equation}
  with $Q$-probability one.
  Since $P \in \Pcal$ was arbitrary and $W$ does not depend on $\P$, we have that
  \begin{equation}
    \limsup_{n\to \infty} \frac{1}{n} \log \left ( W_n \right ) \leq \inf_\Pin \KL(\Q \| \P) = \ell_Q^\star,
  \end{equation}
  and the equality follows from \cref{lemma:klinf identity}. This completes the proof.
\end{proof}
\cref{proposition:growth rate upper bound} begs the question of whether there exist processes with almost sure growth rates attaining $\ell_Q^\star$. 
It is known from \citet[Theorem 15.3.1]{cover1999elements} that for a fixed $Q \in \Qcal$,
\begin{equation}\label{eq:q-oracle}
  \eproc_n(\lambda_Q) := \prod_{i=1}^n \lambda_Q^\trans \bE_i
\end{equation}
is $Q$-almost surely log-optimal and by the strong law of large numbers, $\frac{1}{n} \log (\eproc_n(\lambda_Q)) \to \ell_Q^\star = \ell_Q(\lambda_Q)$ with $Q$-probability one. However, \eqref{eq:q-oracle} is not $\Qcal$-instance-log-optimal in the composite sense of \cref{definition:log-optimality} in general since its construction depends on the particular distribution $Q$. Instead, as we show in the following theorem, $e$-processes with sublinear portfolio regret are $\Qcal$-instance-log-optimal and that their asymptotic growth rates are given by $\ell_Q^\star$ for every $\Qin$.

\begin{theorem}[Instance log-optimality of $e$-processes with sublinear portfolio regret]\label{theorem:adaptive optimality}
  Let $\eprocclass$ be the class of all $\Pe$-processes.
  If $\eproc \in \eprocclass$ is a $\Pcal$-$e$-process satisfying a portfolio
  regret bound of $\Regret_n \leq r_n$ for some sublinear $r_n = o(n)$, then
  \begin{enumerate}[label = (\roman*)]
  \item $\eproc$ is $\Qcal$-instance-log-optimal in $\eprocclass$.
  \item Fix any $Q \in \Qcal$ and let $\eproc(\lambda_Q)$ be the test
      supermartingale given in \eqref{eq:q-oracle}. $\eproc$ and 
      $\eproc(\lambda_Q)$ are
      $Q$-almost surely asymptotically equivalent, meaning that:
      \begin{equation}
        \lim_{n \to \infty} \left(\frac{1}{n} \log(\eproc_n) -
        \frac{1}{n} \log(\eproc_n(\lambda_Q))\right) = 0 \quad Q\text{-almost surely}\,.
      \end{equation}
    Moreover, if $\infseqn{a_n}$ is a sequence such that $a_n \nearrow \infty$ for which $\sum_{k=1}^\infty \exp \left \{ -a_k \eps / 2 \right \}
    < \infty$ for any $\eps > 0$ and $a_n^{-1} r_n \to 0$, then the
    aforementioned asymptotic equivalence holds with a rate of
    $\infseqn{a_n/n}$.
    \item For every $\Qin$, $\eproc$ has an asymptotic growth rate of
    $\ell_Q^\star \coloneqq \EE_Q\left[\log\left(\lambda_Q^\top \bE_1 \right)\right]$, 
    meaning that
    \begin{equation}
        \lim_{\nto}\frac{1}{n} \log W_n = \ell_Q^\star
        \quad\text{$Q$-almost surely.}
    \end{equation}
  \end{enumerate}
\end{theorem}

The proof of \cref{theorem:adaptive optimality} can be found in
\cref{proof:adaptive optimality-specific} and proceeds by first deriving a
nonasymptotic concentration inequality that yields $(i)$--$(iii)$ as consequences.

As an application of \cref{theorem:adaptive optimality}, we can rely on the
sharp regret bounds for universal portfolios provided by
\citet{cover1999elements} (as well as improvements in the one-dimensional case
by \citet{orabona2021tight}) to derive explicit $e$-processes by taking the
exponential of the empirical maximum log wealth less the regret as suggested by
\citet{orabona2021tight}. The details are provided in the following corollary of
\cref{theorem:adaptive optimality}. 
\begin{corollary}\label{corollary:adaptive optimality of UP and OJ}
  Consider the same setup as \cref{theorem:adaptive optimality}.
  and let $\infseqn{\lambda_n^\UP}$ denote the outputs of the universal portfolio algorithm as described in \eqref{eq:UP-lambda}. Then the process $\infseqn{W_n^\UP}$ given by
  \begin{equation}
    W_n^\UP := \prod_{i=1}^n (\lambda_i^\UP)^\trans \bE_i
  \end{equation}
  is a test $\Pcal$-supermartingale with logarithmic portfolio regret. Furthermore, the process $\infseqn{W_n^\CO}$ given by
  \begin{equation}
    W_n^\CO := \exp \left \{ \max_{\lambda \in \dsimp} \log \eproc_n(\lambda) - d\log (n + 1)/2 - \log (2) \right \}
  \end{equation}
  is a $\Pe$-process with logarithmic portfolio regret. Finally, in the case of $d = 1$, the process $\infseqn{W_n^\OJ}$ given by
  \begin{equation}
    W_n^\OJ := \exp \left \{ \log \eproc_n(\lambda^\Max_n) - \max_{j = 0, 1, \dots, n} \log \left ( \frac{\pi (\lambda_n^\Max)^j (1-\lambda_n^\Max)^{n-j} \Gamma(n+1)}{\Gamma(j+1/2) \Gamma(n-j + 1/2)} \right ) \right \}
  \end{equation}
  is a $\Pe$-process with logarithmic portfolio regret where $\lambda_n^\Max := \argmax_{\lambda \in \onesimp} \sum_{i=1}^n \log ( (1 - \lambda^\brackzero) E_i^\brackzero + \lambda^\brackone E_i^\brackone )$. Consequently, $W^\UP$, $W^\CO$, and $W^\OJ$ all satisfy the conditions of \cref{theorem:adaptive optimality} and are hence $\Qcal$-instance-log-optimal.
\end{corollary}
The fact that $\eproc^\UP_n$ forms a $\Pcal$-$e$-process (indeed, a test $\Pcal$-supermartingale) follows from the fact that the process in \eqref{eq:intro-general-test-sm} forms a test $\Pcal$-supermartingale combined with the fact that $\infseqn{\lambda_n^\UP}$ is predictable and $\simplex$-valued as can be deduced from its definition in \eqref{eq:UP-lambda}. Furthermore, we have that $\eproc^\CO_n$ and $\eproc^\OJ_n$ both form $\Pcal$-$e$-processes as they are pathwise upper-bounded by $\eproc^\UP_n$ for every $n \in \NN$.

\begin{figure}[!h]
    \centering
    \begin{subfigure}{0.495\linewidth}
        \centering
        \includegraphics[width=\linewidth]{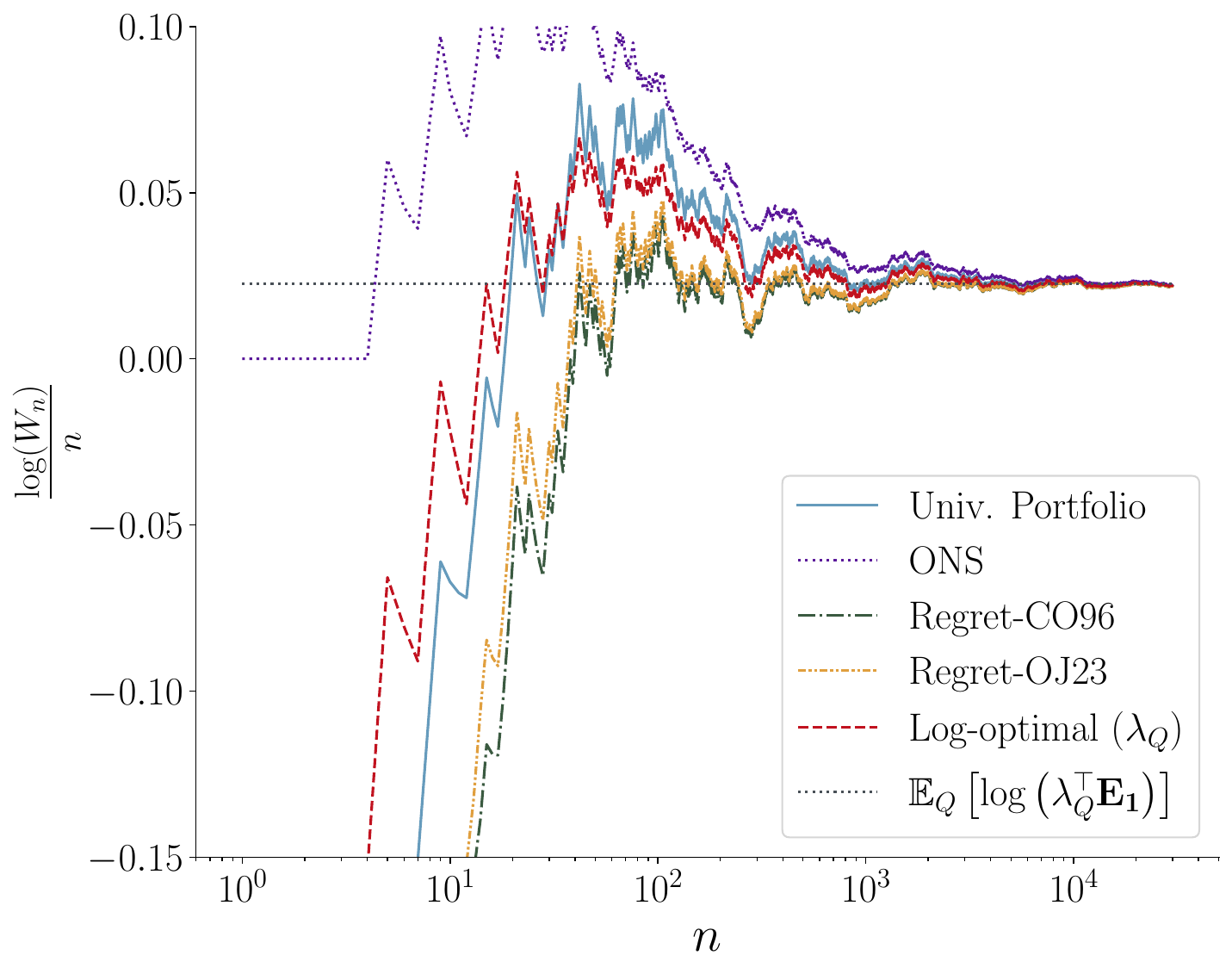}
        \caption{$\Pcal = \{ P : \EE_P[X_1] \leq 0.3\}$, $Q = \text{Bernoulli}(0.4)$}
    \end{subfigure}
    \hfill
    \begin{subfigure}{0.495\linewidth}
       \centering 
        \includegraphics[width=\linewidth]{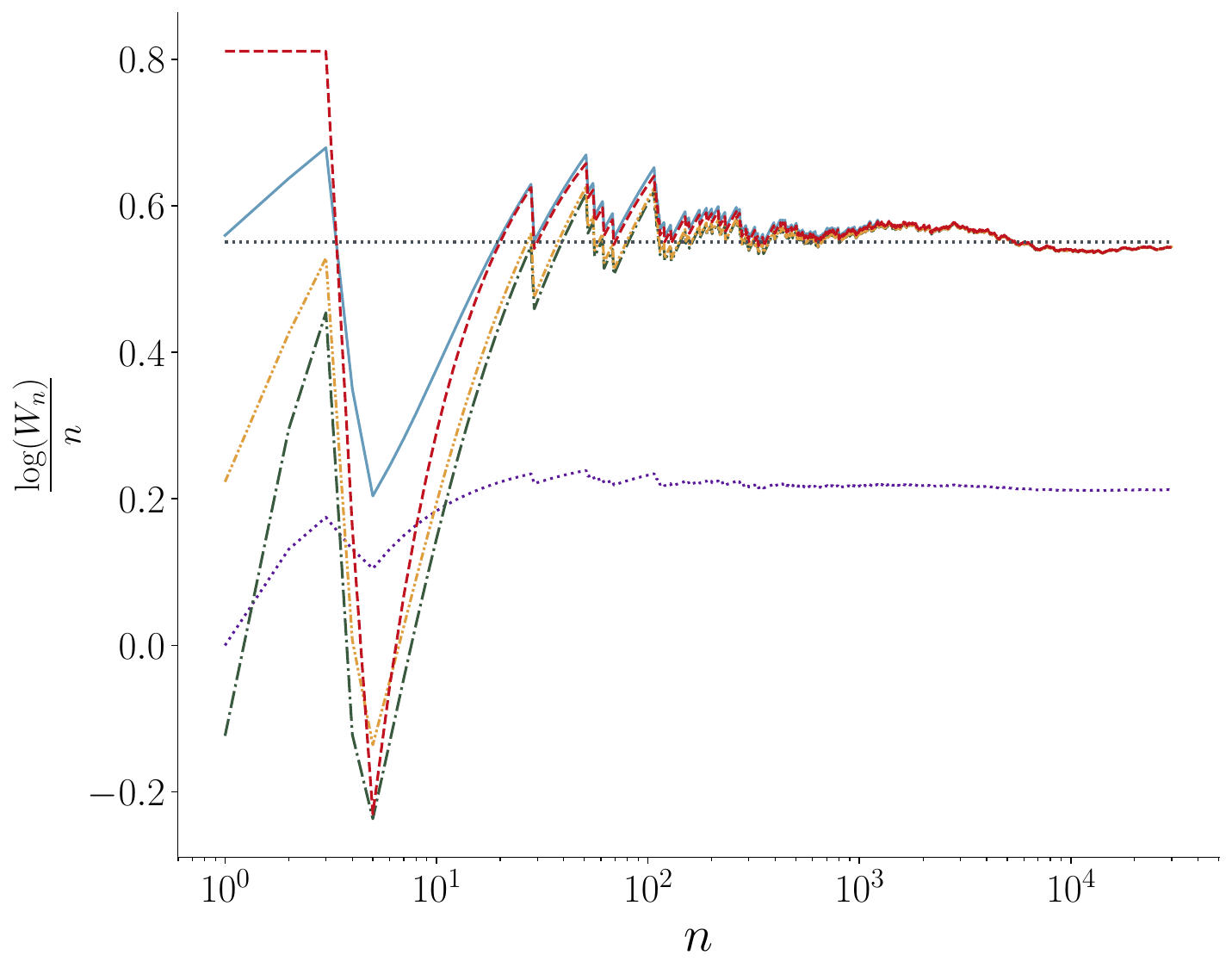}
        \caption{$\Pcal = \{ P : \EE_P[X_1] \leq 0.4\}$, $Q = \text{Bernoulli}(0.9)$}
    \end{subfigure}
    \caption{Empirical growth rates $\log (\eproc_n) / n$ for various $e$-processes that test whether the mean of a bounded random variable is at most $0.3$ or $0.4$---for scenarios (a) and (b), respectively (the exact $e$-processes are discussed in \cref{corollary:adaptive optimality of UP and OJ} and \cref{section:bounded one sided}). In the left-hand and right-hand side plots, the true distributions of the bounded observations are Bernoulli(0.4) and Bernoulli(0.9) and hence the optimal log-wealth increments $\ell_Q^\star$ are 0.023 and 0.55, respectively. In particular, the maximizer $\lambda_Q$ of the growth rate is within the implicit range $[1/4, 3/4]$ that is available in the regret bound for ONS in the first scenario while it is outside of that range in the second.  Hence, in the latter case, only those $e$-processes with sublinear portfolio regret (Univ.~Portfolio, Regret-CO96, Regret-OJ23) are asymptotically equivalent to the wealth of the log-optimal strategy $\lambda_Q$. 
}
    \label{fig:log-optimality}
\end{figure}
\begin{remark}[On the use of $W^\CO$ and $W^\OJ$ over $W^\UP$.]
One may wonder what is the advantage of using the lower-bounding processes $W^\CO$ and $W^\OJ$ over $W^\UP$ since the latter is strictly larger in general. The universal portfolio algorithm of \citet{cover1991universal} can display computationally unstable behavior especially for large values of $n$ since its implementation involves taking the ratio of integrals of degree $n$ polynomials. Under the alternative hypothesis, these processes tend to become large and hence carrying out Cover's algorithm involves computing and storing large coefficients of large-degree polynomials. On the other hand, the regret-based $e$-processes $W^\CO$ and $W^\OJ$ only require the solution to a straightforward convex optimization problem and are simpler to compute as a result. Asymptotically, they are all equivalent and enjoy the same $\Qcal$-instance-log-optimality guarantees.

As seen in \cref{fig:log-optimality}, the difference in empirical growth rates tends to be negligible for large $n$, but the universal portfolio algorithm does see non-negligible benefits at smaller sample sizes. Consequently, as a broad recommendation, we suggest implementing Cover's algorithm up until the sample size $n'$ where computational instability appears, at which point the algorithm can switch to $W_{n'}^\CO, W_{n'+1}^\CO, \dots$ (or $W^\OJ$ if $d = 1$). 
\end{remark}

\begin{remark}[Relationship to the GREE strategy of \citet{wang2022backtesting}]\label{remark:gree}
    In a related work, \citet{wang2022backtesting} focus on the test supermartingales given in \eqref{eq:intro-general-test-sm-special-case}:
    \begin{equation}
      \eproc_n := \prod_{i=1}^n (1 + \lambda_i^\brackone (E_i - 1)),
    \end{equation}
    and they study growth optimality properties when $\infseqn{E_n}$ are \iid{} $\Pcal$-$e$-values.
    In particular, they show that if $\infseqn{\lambda_n^\brackone}$ are chosen according to a follow-the-leader-type strategy---referred to as \emph{growth-rate for empirical $e$-statistics} (GREE)---defined by
    \begin{equation}
      \lambda_n^\GREE := \argmax_{\lambda \in [0, 1]}\frac{1}{n-1}\sum_{i=1}^{n-1} (1 + \lambda (E_i - 1)),
    \end{equation}
    then the resulting $e$-process $\eproc_n^\GREE := \prod_{i=1}^n (1 + \lambda_i^\GREE (E_i - 1))$ will be asymptotically log-optimal in $L_1(Q)$ \citep[Definition 3]{wang2022backtesting}. Since $L_1(Q)$ convergence neither implies nor is implied by $Q$-almost sure convergence, $L_1(Q)$-log-optimality is neither stronger nor weaker than \cref{definition:log-optimality}. Nevertheless, in \cref{sec:l1-log-optimality} we show that sublinear portfolio regret is also a sufficient condition for satisfying $L_1(Q)$-log-optimality.
    We end this remark by noting how the GREE strategy from \citep{wang2022backtesting} does not enjoy sublinear portfolio regret, which stems from the fact that the follow-the-leader strategy \citep[\S 5, p. 71]{hazan2016introduction} incurs linear regret in the worst case even under exp-concave and bounded utilities. Whether or not GREE satisfies $Q$-almost sure instance log-optimality under certain distributional assumptions is an open question.
\end{remark}

In \cref{section:uniformity,section:generalizations}, we
show that the growth-rate log-optimality and numeraire properties hold 
in a stronger distribution-uniform sense and for a more abstract class of
$e$-processes. Nonetheless, growth-rate log-optimality is the first of the two
notions of power for $e$-processes that were discussed in the introduction. We
now turn our attention to the second: bounds on the expected rejection
time of the level-$\alpha$ sequential test obtained by thresholding an
$e$-process at $1/\alpha$ for some $\alpha \in (0, 1)$. As we will see, lower
and upper bounds are both characterized by the reciprocal of the same
log-optimal growth rate $\ell_Q^\star$ appearing in \cref{theorem:adaptive
optimality}.

\section{Sharp bounds on the expected rejection time}\label{section:stopping time}
In the classical (non-sequential) regime, one is often interested in finding a test with a small ``sample complexity,'' meaning that the number of samples required to reject some null hypothesis is small. In the context of sequential testing, the sample complexity is synonymous with the stopping time $\tau_\alpha := \inf \{ n \in \NN : \eproc_n \geq 1/\alpha \}$, which is now a random variable. As such, it is typically of interest to find bounds on its expectation, especially one that scales with the ``difficulty'' of the test, both in terms of the desired type-I error level $\alpha \in (0, 1)$ and certain properties of the alternative $Q \in \Qcal$ \citep{breiman1961optimal,chugg2023auditing,shekhar2023nonparametric,garivier2016optimal,agrawal2021optimal,kaufmann2016complexity,kaufmann2021mixture}. Informally speaking, upper bounds (and corresponding lower bounds, if they exist) of expected rejection times typically take the form
\begin{equation}\label{eq:stopping-time-bound-prose-general}
  \frac{\EE_Q \left [ \tau_\alpha \right ]}{\log(1/\alpha)} \lesssim \frac{1}{\text{easiness under $Q$}}.
\end{equation}
The ``easiness'' of rejecting $\Pcal$ under an alternative $Q$ in \eqref{eq:stopping-time-bound-prose-general} is problem-specific; for example, in the context of bounded mean testing as in \cref{section:bounded one sided,section:bounded two sided,section:difference-in-means}, \citet{chugg2023auditing} derive a bound of the form in \eqref{eq:stopping-time-bound-prose-general} where the ``easiness under $Q$'' is measured through $\Delta^2 = (\mu_Q - \mu_0)^2$, so that the problem of distinguishing between the true mean $\mu_Q$ under $Q$ and the null mean $\mu_0$ becomes easier when they are farther apart. As another example, when testing a simple null, $H_0: X_1 \sim P$ versus a simple alternative $H_1 : X_1 \sim Q$, the rejection time of any sequential test must satisfy a lower bound $\EE_Q[\tau] \geq \log(1/\alpha) / \DKL(Q \| P)$, a classical result due to \citet[Eq.~(4.81)]{wald1945sequential} (see also \citet{garivier2016optimal}, \citet[Lemma 4]{shekhar2023nonparametric}, and \citet{agrawal2025stopping}).
In what follows, we derive information-theoretic lower bounds and matching upper bounds on the expected rejection time in the $\alpha \to 0^+$ regime. 

We remark that while growth-rate log-optimality has been a central criterion appearing in information theory \citep{kelly1956new,krichevsky1981performance}, mathematical finance \citep{breiman1961optimal,long1990numeraire,karatzas2007numeraire}, and more recently sequential inference, bounds on expected stopping times are ubiquitous in the stochastic multi-armed bandit literature \citep{burnetas1996optimal,garivier2016optimal,kaufmann2016complexity,agrawal2021optimal,agrawal2021optimal,honda2010asymptotically,honda2011asymptotically,honda2015non}. Specifically, it is one of the key criteria in evaluating algorithms for the fixed-confidence best-arm identification problem. Similar to the guarantees that we provide in the results to follow, those works often provide nonasymptotic lower bounds on expected stopping times that are then attained by specific algorithms in the so-called ``high-confidence'' regime, i.e., as $\alpha \to 0^+$.
\subsection{Optimal expected rejection times via sublinear portfolio regret}\label{section:rejection time from sublinear portfolio regret}

Let us begin the present section by discussing lower and upper bounds for constant rebalanced portfolios (\cref{proposition:expected stopping time of constant betting strategy}) and later discuss lower bounds for arbitrary portfolios (\cref{proposition:expected stopping time oracle}) and corresponding upper bounds for those with sublinear portfolio regret (\cref{theorem:stopping time}).
The following result shows that in the setting presented in \eqref{eq:intro-general-test-sm}, any constant rebalanced portfolio $\lambda \in [0, 1]$ has a $\log(1/\alpha)^{-1}$-scaled expected rejection time given by $1/\ell_Q(\lambda)$ in the $\alpha \to 0^+$ regime.
\begin{proposition}[Expected rejection times of constant rebalanced portfolios]\label{proposition:expected stopping time of constant betting strategy}
  Let $\lambda \in \dsimp$ be any constant portfolio on $d$ $e$-values. Define $\eproc_n(\lambda) := \prod_{i=1}^n \lambda^\trans \bE_i$. Then for any $\alpha \in (0, 1)$, the expectation of the first rejection time $\tau_\alpha(\lambda) := \inf \{ n \in \NN : \eproc_n(\lambda) \geq 1/\alpha \}$ has the following nonasymptotic lower bound:
  \begin{equation}\label{eq:constant betting strategy lower bound}
    \frac{\EE_Q[\tau_\alpha(\lambda)]}{\log (1/\alpha)} \geq \frac{1}{\ell_Q(\lambda)}.
  \end{equation}
  Furthermore, if the logarithmic increments have a finite $s^\tth$ moment for $s > 2$,
  \begin{equation}
    \rho_s(\lambda) := \EE_Q \left \lvert \log \left ( \lambda^\trans \bE_1 \right ) - \ell_Q(\lambda) \right \rvert^s < \infty,
  \end{equation}
  Then, the expected rejection time matches the above lower bound in the $\alpha \to 0^+$ regime:
  \begin{equation}
    \lim_{\alpha \to 0^+} \frac{\EE_Q \left [ \tau_\alpha(\lambda) \right ]}{\log(1/\alpha)} = \frac{1}{\ell_Q(\lambda)}.
  \end{equation}
\end{proposition}
A proof can be found in \cref{proof:expected stopping time of constant betting strategy}.
In the case of $\ell_Q(\lambda) = 0$, the left-hand and right-hand sides of the (in)equalities should be interpreted as taking the value $\infty$.
The results that follow are stated for adaptively selected $\infseqn{\lambda_n}$ but we present \cref{proposition:expected stopping time of constant betting
strategy} due to its relationship to some discussions of \citet[Theorem 1]{breiman1961optimal}. Specifically, \cref{proposition:expected stopping time of constant betting strategy}
provides a more detailed study of the differences in expected rejection times. Indeed, it can be deduced that for any $\lambda \in \simplex$,
\begin{equation}
  \lim_{\alpha \to 0^+} \left ( \frac{\EE_Q[\tau_\alpha(\lambda_Q)]}{\log(1/\alpha)} - \frac{\EE_Q[\tau_\alpha(\lambda)]}{\log(1/\alpha)} \right ) = \frac{1}{\ell_Q^\star} - \frac{1}{\ell_Q(\lambda)} \leq 0,
\end{equation}
which holds with equality if and only if $\lambda = \lambda_Q$. By contrast, \citet[Theorem 1]{breiman1961optimal} did not include a rescaling by $\log(1/\alpha)$, and hence the corresponding right-hand side would be either be $0$ or $-\infty$ if $\lambda = \lambda_Q$ or $\lambda \neq \lambda_Q$, respectively. In the remainder of the section, we assume that $\ell_Q^\star := \max_{\lambda \in \dsimp} \ell_Q(\lambda) > 0$, meaning that \emph{some} betting strategy (and in particular $\lambda_Q$) has a positive growth rate. 

Thus far, we have only derived bounds on expected rejection times for constant rebalanced portfolios. Let us now show that $1/\ell_Q^\star$ is in fact a lower bound on the $\log(1/\alpha)$-rescaled expected stopping time for \emph{any} stopping rule, not just for ones obtained by thresholding $W_n(\lambda)$ for a fixed $\lambda$.
\begin{proposition}[Lower bounds on the expected rejection time]\label{proposition:expected stopping time oracle}
  Let $\tau_\alpha$ be any stopping rule satisfying
  \begin{equation}\label{eq:alphacorrect}
    \sup_\Pin \P \left ( \tau_\alpha < \infty \right ) \leq \alpha.
  \end{equation}
  Then for any alternative distribution $Q \in \Qcal$,
  \begin{equation}\label{eq:stopping time lower bound general}
    \frac{\EE_\Q[\tau_\alpha]}{\log(1/\alpha)} \geq \frac{1}{\ell_Q^\star}.
  \end{equation}
\end{proposition}
The proof can be found in \cref{proof:expected stopping time oracle}, and it relies on an inequality of \citet{wald1945sequential}, alongside an application of Wald's identity, and an identity relating KL divergences to $\ell_Q^\star$ in \cref{lemma:klinf identity}. Note that without the application of \cref{lemma:klinf identity}, the lower bound is written in terms of the reciprocal of an infimum over KL divergences; such results for arbitrary composite nulls can be found in \citet{agrawal2025stopping} (not just for the problem studied here). The condition \eqref{eq:alphacorrect} is referred to as ``$\alpha$-correct'' in the best-arm identification literature.

\begin{figure}[h]
    \centering
    \begin{subfigure}{0.49\linewidth}
        \centering
        \includegraphics[width=\linewidth]{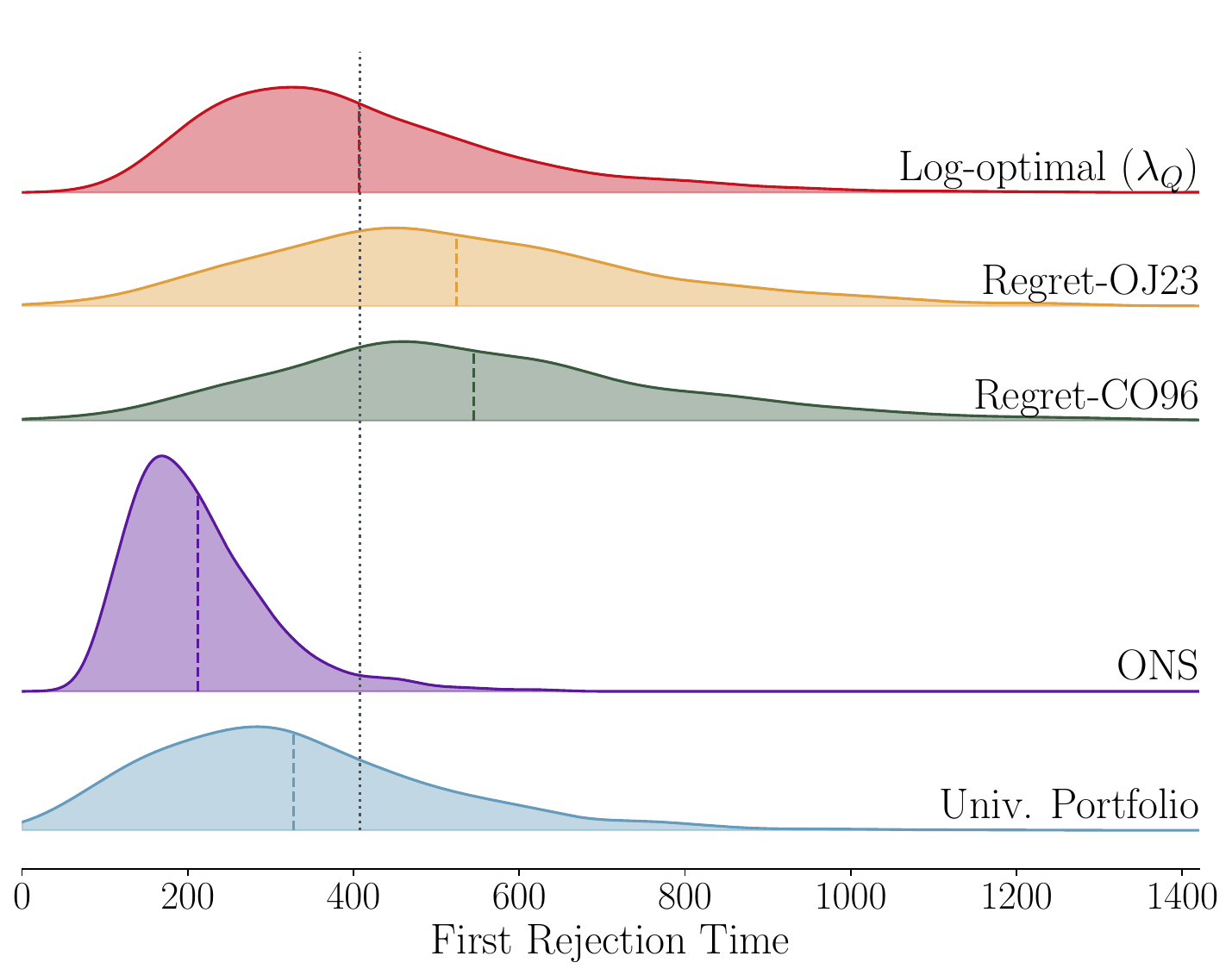}
        \caption{$\Pcal = \{ P : \EE_P[X_1] \leq 0.3\}$, $Q = \text{Bernoulli}(0.4)$}
    \end{subfigure}
    \hfill
    \begin{subfigure}{0.49\linewidth}
       \centering 
        \includegraphics[width=\linewidth]{figures/new_figures/rejection_time_density_alpha_0001_bern_09.pdf}
        \caption{$\Pcal = \{ P : \EE_P[X_1] \leq 0.4\}$, $Q = \text{Bernoulli}(0.9)$}
    \end{subfigure}
   \caption{Distributions of the first rejection time for various $e$-processes under two scenarios for $\alpha = 0.01$.
   In particular, the log-optimal strategy lies inside the allowable range available to ONS in the scenario considered in the left-hand plot whereas it lies outside this range for that of the right-hand plot. Consequently, we see that the distribution of rejection times for $e$-processes with sublinear portfolio regret lie near the optimum $\log(1/\alpha) / \ell_Q^\star$ in both cases, whereas those of ONS only do so in the left-hand plot.}
    \label{fig:expected rejection times}
\end{figure}

With \cref{proposition:expected stopping time oracle} in mind, it is natural to wonder whether there exists a betting strategy achieving a matching \emph{upper} bound of $1/\ell_Q^\star$ in the $\alpha \to 0^+$ regime as in \cref{proposition:expected stopping time of constant betting strategy}. The following theorem provides an answer to this question. As in \cref{theorem:adaptive optimality}, sublinear portfolio regret is a sufficient condition to attain optimality. Furthermore, the following only makes moment assumptions under the log-optimal strategy rather than invoking almost sure boundedness imposed by the data and/or the betting strategy, assumptions which are often used in bounding the expected rejection time in nonparametric settings.
\begin{theorem}[The expected rejection time of sublinear portfolio regret $e$-processes]
    \label{theorem:stopping time}
  Let $\infseqn{\eproc_n}$ be any $\widebar \Pcal$-$e$-process satisfying the portfolio regret bound $\Regret(\eproc_n) \leq r_n$ for some sublinear $r_n = o(n)$. Let $Q \in \Qcal$ be an element of the alternative hypothesis for which
  \begin{equation}
   \rho_s(\lambda_Q) := \EE_Q \left \lvert \log\left ( \lambda_Q^\trans \bE_1 \right ) - \ell_Q^\star \right \rvert^{s} < \infty,  
  \end{equation}
  for some $s > 2$. 
  Then the expectation of the first rejection time $\tau_\alpha := \inf \{ n \in \NN : \eproc_n \geq 1/\alpha \}$ has the property that
  \begin{equation}\label{eq:stopping time asymptotic}
    \lim_{\alpha \to 0^+} \frac{\EE_Q [\tau_\alpha]}{\log (1/\alpha)} = \frac{1}{\ell_Q^\star},
  \end{equation}
  matching the lower bound provided in \cref{proposition:expected stopping time oracle} for $\alpha$ tending to zero.
\end{theorem}
The proof of \cref{theorem:stopping time} can be found in \cref{proof:stopping time}. 
When taken together, \cref{proposition:expected stopping time oracle} and
\cref{theorem:stopping time} allow us to conclude that the bound
$1/\ell_Q^\star$ on $\EE_Q[\tau_\alpha] / \log(1/\alpha)$ as $\alpha \to 0^+$ is
both attainable and unimprovable among any $1/\alpha$-thresholded $e$-process, or more generally among any stopping rule satisfying \eqref{eq:alphacorrect}.
$\eprocclass$. In \cref{section:uniformity} we show that such
a conclusion holds in a stronger distribution-uniform sense.
See \cref{fig:expected rejection times} for empirical distributions of rejection times for the $\Pcal$-$e$-processes with logarithmic portfolio regret that were discussed in \cref{corollary:adaptive optimality of UP and OJ}.

Let us now make some remarks on the proof techniques used to arrive at \cref{theorem:stopping time} and how they differ from those that can be found in the literature. While we are not aware of bounds on rejection times that have been derived for test supermartingales of the exact form \eqref{eq:intro-general-test-sm}, there do exist some explicit results in specific nonparametric contexts of differences-in-bounded-means testing \citep{chugg2023auditing,chen2025optimistic} and two-sample testing \citep{shekhar2023nonparametric}. Even in these cases, some comparisons are in order.
 Existing proofs in the literature sometimes rely in some way on almost sure boundedness, either through conditions on the betting strategy $\infseqn{\lambda_n}$ or on the observed random variables themselves, or both. Boundedness is sometimes exploited through an application of a multiplicative Chernoff method which yields sub-Gaussian concentration; see \citep{chugg2023auditing,chen2025optimistic}. By contrast, in our case, we neither assume that the log-wealth increments nor the $e$-values are bounded and we make no use of multiplicative Chernoff bounds nor sub-Gaussianity. Instead, we derive certain concentration results through finiteness of the $s^\tth$ moment for some $s > 2$ along with a Chebyshev-like inequality that can be deduced from an inequality of \citet{nemirovski2000topics}; see \cref{lemma:concentration-via-nemirovski} for details. Exceptions to this reliance on boundedness exist, with expected stopping time bounds that rely on some structure other than convex combinations of $e$-values. For example, existing works may rely on moment bounds \citep{agrawal2021optimal} or exponential families \citep{garivier2016optimal}.

 With the main results from \cref{section:asymptotic-equivalence,section:stopping time} in mind, we now turn our attention to some testing problems that are special cases of the $e$-processes considered thus far.

\section{Implications for some sequential testing problems}\label{section:corollaries}

The purpose of this section is to discuss a few nonparametric sequential testing problems for which our main results can be instantiated to yield optimal bounds on growth rates and expected rejection times. \cref{section:bounded one sided,section:bounded two sided} describe test (super)martingales that can be used to test whether the mean of a bounded random variable is at most or exactly equal to some prespecified null value, and \cref{section:difference-in-means} discusses a related problem of testing whether the difference in means of two bounded random variables is different from zero. We additionally discuss a simplified version of the two-sample and marginal independence testing problems in \cref{section:two-sample} and discuss how some of their key aspects can be reduced to the difference-in-means testing problem of \cref{section:difference-in-means}. However, this final discussion is left to the supplementary material as it requires additional context on integral probability metrics and distance measures with variational representations which are outside the main scope of this paper.
\subsection{One-sided tests for the mean of a bounded random variable}\label{section:bounded one sided}

Returning to (a slightly simplified version of) the bounded mean testing problem given in \cref{section:testing problems to keep in mind}, consider the \emph{one-sided} null and alternative hypotheses given by
\begin{equation}
  \Pcal^\leq := \{ P : \EE_P [X_1] \leq \mu_0 \} \quad\text{versus}\quad \Qcal^> := \{ P : \EE_P [X_1] > \mu_0 \}.
\end{equation}
and the $\Pe$-process
\begin{equation}\label{eq:one-sided test supermartingale}
  \eproc_n^\leq := \prod_{i=1}^n (1 + \gamma_i (X_i - \mu_0))
\end{equation}
for some predictable sequence $[0, 1/\mu_0]$-valued sequence $\infseqn{\gamma_n}$. It is easy to check that $W_n^\leq$ is a special case of the test supermartingale in \eqref{eq:intro-general-test-sm} when taking
\begin{equation}
  \left (d, E_n^\brackzero, E_n^\brackone, \lambda_n^\brackone \right ) = \left (1, 1, \frac{X_n}{\mu_0}, \mu_0 \gamma_n \right )\quad \text{for each $n \in \NN$.}
\end{equation}
Invoking both \cref{theorem:adaptive optimality,theorem:stopping time}, we have the following corollary for one-sided tests for the mean of a bounded random variable.
\begin{corollary}[Growth rates and expected rejection times for one-sided bounded mean testing]\label{corollary:one-sided bounded means}
  If $\infseqn{\lambda_n}$ is chosen in such a way so that $\eproc_n^\leq$ has a sublinear portfolio regret (e.g., through the constructions $\eproc_n^\UP, \eproc_n^\CO$, and $\eproc_n^\OJ$ as described in \cref{corollary:adaptive optimality of UP and OJ}), then we have that for any $Q \in \Qcal^>$,
  \begin{equation}
    \lim_{n \to \infty} \frac{1}{n} \log (\eproc_n^\leq) = \max_{\gamma \in [0,1/\mu_0]} \EE_Q \left [ \log (1 + \gamma (X_1 - \mu_0)) \right ] \quad\text{$Q$-almost surely,}
  \end{equation}
  and this is unimprovable in the sense that for any other $\Pcal^\leq$-$e$-process $\eproc_n'$, and any alternative distribution $Q \in \Qcal^>$ it must hold that
  \begin{equation}
    \limsup_{n\to \infty} \frac{1}{n} \log (\eproc_n') \leq \max_{\gamma \in [0, 1/\mu_0]} \EE_Q \left [ \log (1 + \gamma (X_1 - \mu_0)) \right ] \quad\text{$Q$-almost surely.}
  \end{equation}
  Furthermore, for any $Q \in \Qcal^>$, letting $\gamma_Q^\star \in \argmax_{\gamma \in [0,1/\mu_0]} \EE_Q [\log(1 + \gamma(X_1 - \mu_0))]$ and $\tau_\alpha := \inf \{ n \in \NN : W_n^\leq \geq 1/\alpha \}$, we have that
  \begin{equation}
    \lim_{\alpha \to 0} \frac{\EE_Q \left [ \tau_\alpha \right ] }{\log(1/\alpha)} = \frac{1}{\EE_Q [\log( 1 + \gamma_Q^\star (X_1 - \mu_0))]},
  \end{equation}
  and this is once again unimprovable in the sense that for any stopping rule $\tau_\alpha'$ satisfying \eqref{eq:alphacorrect}, the right-hand side is a lower bound on $\EE_\Q[\tau_\alpha']$.
\end{corollary}

We remark that when taken together and translated from regret minimization to the hypothesis testing problem underlying best-arm identification, the works of \citet{burnetas1996optimal} and \citet{honda2010asymptotically,honda2011asymptotically} obtain the same lower and upper bounds on the expected rejection time for this one-sided null.
In the following section, we consider the slightly more complicated setting of equality nulls and two-sided alternatives.

\subsection{Two-sided tests for the mean of a bounded random variable}\label{section:bounded two sided}
Consider the example problem presented in \cref{section:testing problems to keep in mind} focused on testing whether the mean of a bounded random variable is equal to $\mu_0 \in [0, 1]$. That is, we will consider the equality null and two-sided alternative analogues of \cref{section:bounded one sided}:
\begin{equation}
  \Pcal^= := \{ P :\EE_P [X_1] = \mu_0 \}\quad\text{versus}\quad \Qcal^{\neq} := \{ P : \EE_P[X_1] \neq \mu_0 \},
\end{equation}
and define the test $\Pcal^=$-martingale as in \citet{waudby2020estimating}:
\begin{equation}
  \eproc_n^= := \prod_{i=1}^n \left ( 1 + \gamma_i (X_i - \mu_0) \right ),
\end{equation}
where, unlike in \cref{section:bounded one sided}, $\infseqn{\gamma_n}$ may now take values in $[-1/(1-\mu_0), 1/\mu_0]$. As mentioned in \cref{section:testing problems to keep in mind}, the connection to the test supermartingale given in \eqref{eq:intro-general-test-sm} is that $\eproc_n^=$ can be equivalently written as
\begin{equation}\label{eq:bounded two sided martingale}
  \eproc_n^= = \prod_{i=1}^n \left ( (1-\lambda_i^\brackone) \frac{1-X_i}{1-\mu_0} + \lambda_i^\brackone \frac{X_i}{\mu_0} \right ),
\end{equation}
where we have taken
\begin{equation}
  \left ( d, E_n^\brackzero, E_n^\brackone, \lambda_n^\brackone \right ) = \left (1, \frac{1-X_n}{1-\mu_0}, \frac{X_n}{\mu_0}, \mu_0 + \gamma_n \mu_0 (1 - \mu_0) \right ).
\end{equation}
This is also the form of a test martingale that can be found in \citet[Algorithm 1]{orabona2021tight}, \citet{ryu2024confidence}; see also \citep{ryu2024gambling,ryu2025improved}.
By \citet[Proposition 2]{waudby2020estimating}, processes given by $\eproc_n^=$ are not just \emph{some} of potentially many test $\Pcal^=$-martingales, they in fact consist of \emph{all} such martingales. That is, if a process is a test $\Pcal^=$-martingale, it can be written in the form \eqref{eq:bounded two sided martingale} for some predictable $\onesimp$-valued $\infseqn{\lambda_n}$. Going further, \citet{clerico2025optimality} shows that $W_n^=$ consists of all admissible $\Pcal^=$-$e$-processes, i.e., those that cannot be pathwise dominated by another.

Invoking \cref{theorem:adaptive optimality,theorem:stopping time} in this setting, we have the following corollary for bounded mean testing.
\begin{corollary}[Log-optimality and expected rejection times for two-sided bounded mean testing]\label{corollary:two-sided bounded means}
  If $\infseqn{\lambda_n}$ is chosen in such a way so that $\eproc_n^=$ has a sublinear portfolio regret, such as those described in \cref{corollary:adaptive optimality of UP and OJ}, then we have that for any $Q \in \Qcal^\neq$,
  \begin{equation}
    \lim_{n \to \infty} \frac{1}{n} \log \eproc_n^= = \max_{\gamma \in [-1/(1-\mu_0), 1/\mu_0]} \EE_Q \left [ \log (1 + \gamma (X_1 - \mu_0)) \right ] \quad\text{$Q$-almost surely,}
  \end{equation}
  and this is unimprovable in the sense that for any $\Pcal^=$-$e$-process $\eproc_n'$, it must hold that
  \begin{equation}\label{eq:two sided growth rate upper bound}
    \limsup_{n\to \infty} \frac{1}{n} \log ( \eproc_n') \leq \max_{\gamma \in [-1/(1-\mu_0), 1/\mu_0]} \EE_Q \left [ \log (1 + \gamma (X_1 - \mu_0)) \right ] \quad\text{$Q$-almost surely.}
  \end{equation}
  Furthermore, letting $\gamma_Q^\star$ be the maximizer of $\EE_Q[\log(1 + \gamma(X_1 - \mu_0))]$ over the range $[-1/(1-\mu_0), 1/\mu_0]$, we have that
  \begin{equation}
    \lim_{\alpha \to 0} \frac{\EE_Q \left [ \tau_\alpha \right ] }{\log(1/\alpha)} = \frac{1}{\EE_Q [\log (1 + \gamma_Q^\star (X_i - \mu_0))]},
  \end{equation}
  and this is once again unimprovable over the class of all $\Pcal^=$-$e$-processes.
\end{corollary}
The lower bound on the $\Q$-almost sure limit superior of $\frac{1}{n}\log W_n^=$ and the asymptotic upper bound on the expected rejection time are both immediate from \cref{theorem:adaptive optimality} and \cref{theorem:stopping time}, respectively. It takes a few extra steps to justify why the upper bound in \eqref{eq:two sided growth rate upper bound} and the lower bound on $\EE_\Q[\tau_\alpha] / \log(1/\alpha)$ are also consequences since in this case, $\bE_1$ does not contain a trivial $e$-value of 1 and thus cannot be described by the null hypothesis in \eqref{eq:null}. Nevertheless, by \cref{theorem:adaptive optimality}, it holds that
\begin{equation}\label{eq:discussion a.s. upper bound}
  \limsup_{n\to\infty} \frac{1}{n} \log( W_n') \leq \max_{\lambda \in \triangle_2} \EE_\Q \left [ \log\left (\lambda^\brackzero + \lambda^\brackone \frac{1-X_1}{1-\mu_0} + \lambda^\bracktwo \frac{X_1}{\mu_0} \right ) \right ],
\end{equation}
and similarly for the lower bound on the expected rejection time. One can assume without loss of generality that $\lambda^\brackzero = 0$ since if $(\lambda_\Q^\brackzero, \lambda_\Q^\brackone, \lambda_\Q^\bracktwo)$ is a maximizer to the program in \eqref{eq:discussion a.s. upper bound}, then so is $(0, \lambda_\Q^\brackone + (1 - \mu_0)\lambda_\Q^\brackzero, \lambda_\Q^\bracktwo + \mu_0 \lambda_\Q^\brackzero)$ due to the identity
\begin{equation}
    \lambda_\Q^\brackzero = \left ( (1-\mu_0)\frac{1-X_1}{1-\mu_0} + \mu_0\frac{X_1}{\mu_0} \right )\lambda_\Q^\brackzero.
\end{equation}

Let us now demonstrate how some deterministic inequalities for the function $y \mapsto \log (1 + y)$ can be used to derive lower bounds on growth rates and upper bounds on expected rejection times that qualitatively resemble existing bounds in the literature that invoke the ONS strategy. Consider the difference of alternative and null means under $Q \in \Qcal^\neq$:
\begin{equation}
  \Delta_Q := \EE_Q[X_1] - \mu_0,
\end{equation}
and let $\sigma_Q^2 := \Var_Q(X_1)$ denote the variance of the random variable $X_1$ under $Q$.
Using the inequality $\log(1 + y) \geq y - y^2$ for all $y \in [-1/2,1/2]$, we have that $\log(1+\gamma(X_1 - \mu_0)) \geq \gamma(X_1 - \mu_0) - \gamma^2(X_1 - \mu_0)^2$ whenever $\gamma \in [-1/2, 1/2]$, and it is not hard to check that the maximizer of the expectation of this lower bound is given (and can be further lower-bounded) by
\begin{equation}
  \argmax_{\gamma \in [-1/2,1/2]} \EE_Q \left [ \gamma(X_1 - \mu_0) - \gamma^2(X_1 - \mu_0) \right ] = \frac{\Delta_Q}{2 (\Var_Q(X_1) + \Delta_Q^2)} \geq \frac{\Delta_Q}{2(1/4 + \Delta_Q^2)},
\end{equation}
where we note that the right-hand side always lies in $[-1/2, 1/2]$. Plugging this back into aforementioned lower bound, we can conservatively lower bound the asymptotic growth rate of $\eproc_n^\neq$ and upper bound the expected rejection time as
\begin{equation}\label{eq:conservative bounds without variance}
  \lim_{\nto} \frac{1}{n} \log \eproc_n^\neq \geq \frac{\Delta^2_Q}{1 + 4 \Delta_Q^2}\quad\text{and}\quad \lim_{\alphato} \frac{\EE_Q[\tau_\alpha]}{\log(1/\alpha)} \leq 4 + \frac{1}{\Delta_Q^2}.
\end{equation}
Furthermore, in the regime where $\sigma_Q^2$ is not too small relative to $\Delta_Q$, i.e., if $|\Delta_Q| \leq \frac{1}{2} [ 1-(1-4\sigma_Q^2)^{1/2} ]$,
then the former lower and upper bounds can be written in a way that depend on $\sigma_Q^2$:
\begin{equation}\label{eq:conservative bounds with variance}
  \lim_{\nto} \frac{1}{n} \log \eproc_n^\neq \geq \frac{\Delta^2_Q}{4 (\sigma_Q^2 + \Delta_Q^2)}\quad\text{and}\quad \lim_{\alphato} \frac{\EE_Q[\tau_\alpha]}{\log(1/\alpha)} \leq 4 + \frac{4\sigma_Q^2}{\Delta_Q^2}.
\end{equation}
Qualitatively, the expected rejection time bounds in \eqref{eq:conservative bounds without variance} and \eqref{eq:conservative bounds with variance} resemble those found in \citet[Proposition 1]{chugg2023auditing} and implicitly in \citet[Proposition 1]{shekhar2023nonparametric}, \citet[Theorem 2]{podkopaev2023sequentialKernelized}, and \citet[Theorem 1]{podkopaev2023sequentialTwoSample} in the contexts of difference-in-means, two-sample, and independence testing, and we make a few more remarks on the former in the following section. Nevertheless, the growth rate and expected rejection time bounds in \cref{corollary:two-sided bounded means} are always sharper than those of \eqref{eq:conservative bounds without variance} and \eqref{eq:conservative bounds with variance}. Moreover, the analyses of the aforementioned prior work rely on almost sure bounds on the log-wealth increments---e.g., $\log(1 + \gamma(X_1 - \mu_0)) \in [a,b]$ uniformly in the parameter $\gamma$ for some $-\infty < a \leq b < \infty$---which is achieved due to the fact that ONS restricts $\gamma$ to lie in $[-1/2,1/2]$ but which cannot be guaranteed to hold when employing Cover's universal portfolio strategy.

\subsection{Difference-in-means testing for bounded random tuples}\label{section:difference-in-means}
Suppose $\infseqn{X_n, Y_n}$ is a sequence of \iid{} tuples taking values in $[0, 1]^2$ and define $D_n := X_n - Y_n$ for each $n$. One may be interested in the one-sided null $\Pcal^{\Dleq} := \{ P  : \EE_P D_1 \leq 0 \}$ or the equality null $\Pcal^{\Deq} := \{ P : \EE_P D_1 = 0 \}$ with the alternatives defined analogously to \cref{section:bounded one sided,section:bounded two sided}. After taking the transformation
\begin{equation}
 Z_n := (D_n + 1)/2 \in [0, 1] \quad\text{for each $n \in \NN$},
\end{equation}
we note that $\Pcal^{\Dleq}$ and $\Pcal^{\Deq}$ are equivalent to $\Pcal^\leq$ and $\Pcal^=$ from \cref{section:bounded one sided} and \cref{section:bounded two sided}, respectively, but with $X_n$ replaced by $Z_n$ and in the special case of $\mu_0 = 1/2$ (with their corresponding alternatives coinciding as well). While this problem can be reduced to a special case of the aforementioned bounded mean testing problems, we still highlight it because it is a nontrivial one for which the literature contains some results about asymptotic growth rates and bounds on expected rejection times \citep{chugg2023auditing,chen2025optimistic} (and less directly, \citep{shekhar2023nonparametric,podkopaev2023sequentialKernelized,podkopaev2023sequentialTwoSample}). Hence, when applied to this problem, our main results (such as \cref{theorem:adaptive optimality,theorem:stopping time}) yield some new insights and improvements on those results in the literature. Indeed, consider the test $\Pcal^\Deq$-martingale $\eproc_n^\Deq$ that is essentially the same as that described in \citet[Algorithm 1]{chugg2023auditing} but with a generic $[-1, 1]$-valued predictable sequence $\infseqn{\gamma_n}$ rather than the $[-1/2, 1/2]$-valued sequence implied by their use of ONS:
\begin{equation}\label{eq:chugg-supermartingale}
  \eproc_n^\Deq := \prod_{i=1}^n \left ( 1 + \gamma_i D_i \right ) = \prod_{i=1}^n \left ( (1-\lambda_i^\brackone) \frac{(1 - Z_i)}{1/2} + \lambda_i^\brackone \frac{Z_i}{1/2} \right ),
\end{equation}
where $\lambda_n^\brackone = (1+\gamma_n)/2$ yields a portfolio $\lambda_n$ taking values in $\onesimp$. Note that if ONS were employed, then $\lambda_n$ would have been restricted to the range $[1/4,3/4]$. Clearly, $\eproc_n^\Deq$ is an instantiation of the test supermartingale in \eqref{eq:intro-general-test-sm} for the $e$-values
\begin{equation}
  E_n^\brackzero:= 2(1-Z_n) \equiv 1/2 - D_n \quad\text{and}\quad E_n^\brackone := 2 Z_n \equiv D_n - 1/2\quad\text{for each $n \in \NN$.}
\end{equation}
Defining $\Delta_Q := \EE_Q [D_1]$,\footnote{Note that \citet{chugg2023auditing} define $\Delta$ as the absolute value of $\EE_Q[D_1]$ whereas we consider its signed version but this will not matter for the expected rejection time discussions that follow since $\Delta_Q$ will only appear when squared, and the use of a signed difference is only for the sake of consistency with \cref{section:bounded two sided}.} the proof of Proposition 1 found in \citet{chugg2023auditing} states that when employing ONS, the resulting test has a rejection time $\tau_\alpha^\ONS := \inf \{ n \in \NN  : \eproc_n^\Deq \geq 1/\alpha\}$ that is bounded in expectation as
\begin{equation}
  \EE_Q \left [ \tau_\alpha^\ONS \right ] \leq \frac{81}{\Delta_Q^2} \log \left ( \frac{162}{\Delta_Q^2 \alpha} \right ) + \frac{\pi^2}{2},
\end{equation}
so that in the $\alphato$ regime considered in \cref{section:stopping time}, it holds that
\begin{equation}
  \limsup_\alphato \frac{\EE_Q \left [ \tau_\alpha^\ONS \right ]}{\log(1/\alpha)} \leq \frac{81}{\Delta_Q^2}.
\end{equation}
Nevertheless, when $\gamma_n := 2\lambda_n - 1$ is chosen in such a way so that $\eproc_n^\Deq$ satisfies a sublinear portfolio regret bound, the following series of inequalities hold for the stopping time $\tau_\alpha := \inf\{ n \in \NN: \eproc_n^\Deq \geq 1/\alpha \}$ as an immediate consequence of \cref{corollary:two-sided bounded means} combined with the discussion thereafter:
\begin{equation}
  \lim_{\alphato} \frac{\EE_Q \left [ \tau_\alpha \right ]}{\log(1/\alpha)} = \frac{1}{\EE_Q [ \log (1 + \gamma_Q^\star D_1) ]} \leq 4 + \frac{1}{\Delta_Q^2} \leq \frac{81}{\Delta_Q^2},
\end{equation}
where the final inequality follows from the fact that $\Delta_Q^2 \leq 1$ by construction.

\section{Distribution-uniform instance-log-optimality}\label{section:uniformity}

In \cref{section:asymptotic-equivalence,section:stopping
time,section:corollaries}, we derived matching lower and
upper bounds on $Q$-almost sure asymptotic growth rates and expected rejection
times under $Q$ for every $Q \in \Qcal$ where $\Qcal$ was a rich family of
alternative distributions. However, these results were distribution-\emph{pointwise}, in the sense that
certain limiting statements were given for every $Q \in \Qcal$ once the
distribution $Q$ was fixed, but they were not shown to hold \emph{uniformly} in
some class $\Qcalu \subseteq \Qcal$. In this section, we will show that all of
the results of this paper can be generalized to hold $\Qcalu$-uniformly for a
rich alternative $\Qcalu$, obtaining the former pointwise results as
corollaries.
Before stating these results, we provide a requisite definition.

\begin{definition}[Distribution-uniform asymptotic almost sure log-optimality
and equivalence]\label{definition:uniform log optimality}
Let $\eprocclass$ be a collection of $\Pcal$-$e$-processes. We say that
$\eproc^\star \equiv \infseqn{\eproc_n^\star}$ is \emph{$\Qcalu$-uniformly and
instance-log-optimal} in $\eprocclass$ if for any other $\eproc \in
\eprocclass$, it holds that
\begin{equation}
  \forall \eps > 0,\quad \lim_{\mto } \sup_\Qinu\Q \left ( \supkm \frac{1}{k} \log( \eproc_k / \eproc_k^\star) \geq \eps \right ) = 0.
\end{equation}
Furthermore, we say that $\eproc^\brackone$ and $\eproc^\bracktwo$ are \emph{$\Qcalu$-uniformly asymptotically equivalent} at a rate of $a_n / n$, for $a_n \nearrow \infty$, if
\begin{equation}
  \forall \eps > 0, \quad \lim_\mto \sup_\Qinu \Q \left ( \supkm a_k^{-1} \left \lvert \log(\eproc_n^\brackone) - \log(\eproc_n^\bracktwo) \right \rvert \geq \eps \right ) = 0.
\end{equation}
\end{definition}
It may not be obvious why \cref{definition:uniform log optimality} serves as a distribution-uniform generalization of \cref{definition:log-optimality}, but the relationship between the two directly mirrors the relationship between notions pointwise and uniform almost sure convergence that are used to describe distribution-uniform strong laws of large numbers \citep{chung_strong_1951,waudby2024distribution}. Indeed, when $\Qcalu = \{Q\}$ is taken to be a singleton, \cref{definition:log-optimality} and \cref{definition:uniform log optimality} are equivalent since for any $\infseqn{Z_n}$, it is the case that $\limsup_n Z_n \leq 0$ with $Q$-probability one \emph{if and only if} $\forall \eps > 0$, $\lim_m \Q \left ( \supkm Z_k \geq \eps \right ) = 0$. We will show the somewhat surprising fact that $\Qcalu$-uniform log-optimality and equivalence will hold for the \emph{entire} alternative $\Qcalu = \Qcal$ while convergence of $n^{-1} \log \eproc_n$ to the optimal rate of growth $\ell_Q^\star$ will hold $\Qcalu$-uniformly within a restricted but rich class $\Qcalu \subseteq \Qcal$ satisfying a certain uniform integrability condition. We summarize these two points as follows.
\begin{theorem}[Uniformly instance-log-optimal $e$-processes]\label{theorem:uniform log optimality}

Let $\eprocclass$ be the class of all $\Pcal$-$e$-processes.
If $\eproc \in \eprocclass$ satisfies a portfolio
regret bound of $\Regret_n \leq r_n$ for some sublinear $r_n = o(n)$, then
\begin{enumerate}[label = (\roman*)]
\item $\eproc$ is $\Qcal$-uniformly and instance-log-optimal in $\eprocclass$.
\item If $a_n \nearrow \infty$, $\sum_{k=1}^\infty \exp \left \{ -a_k \eps / 2 \right \} < \infty$ and $a_n^{-1} r_n \to 0$, then $\eproc$ and $\eproc(\lambda_Q)$ given in \eqref{eq:q-oracle} are $\Qcal$-uniformly asymptotically equivalent at a rate of $\infseqn{a_n / n}$. 
\item Let $\Qcalu \subseteq \Qcal$ be any subset of alternative distributions
    for which the logarithmic wealth increment under $\lambda_Q$ given by
    $L_Q^\star := \log  ( \lambda_Q^\trans \bE_1 )$ has a $\Qcalu$-uniformly integrable
    $p^\tth$ moment for some $p \in [1, 2)$:
  \begin{equation}
\lim_\mto \sup_\Qinu \EE_Q \left [ |L_Q^\star - \ell_Q^\star |^p \1 \{ |L_Q^\star - \ell_Q^\star|^p \geq m \} \right ]=0.
  \end{equation}
  If the portfolio regret bound additionally satisfies $r_n = o(n^{1/p})$, then $\frac{1}{n} \log \eproc_n$ converges $\Qcalu$-uniformly to $\ell_Q^\star$ almost surely at a rate of $o(n^{1/p - 1})$:
  \begin{equation}
    \forall \eps > 0,\quad  \lim_\mto \sup_\Qinu \Q \left ( \supkm  k^{1-1/p} \left \lvert  \frac{1}{k}\log \eproc_k - \ell_Q^\star \right \rvert \geq \eps \right ) = 0.
  \end{equation}
\end{enumerate}
\end{theorem}
The proof of \cref{theorem:uniform log optimality} can be found in a more general form in \cref{proof:uniform log optimality}. Parts $(i)$ and $(ii)$ follow from nonasymptotic probability inequalities that yield the desired results once suprema over $\Qin$ and limits as $\mto$ are taken, while part $(iii)$ additionally relies on uniform strong laws of large numbers due to \citet{chung_strong_1951} and \citet{waudby2024distribution}. In particular, note that if $L_Q^\star$ has a finite $p^\tth$ moment under some $\Qin$, it holds that
\begin{equation}
n^{-1} \log \eproc_n - \ell_Q^\star = o(n^{1/p - 1})\quad\text{$Q$-almost surely.}
\end{equation}
Nevertheless, note that \cref{theorem:uniform log optimality} is a strict generalization of \cref{theorem:adaptive optimality} and the latter is recovered when $\Qcalu = \Qcal = \{Q\}$ is a singleton and $p=1$ so that the condition $r_n = o(n^{1/p})$ reduces to sublinearity.
We explicitly highlight the fact that $\Qcal$-uniform and instance-log-optimality is obtained for free in the \emph{entire} alternative $\Qcal$ without needing to impose any uniform boundedness (or integrability) assumptions on any moments whatsoever. This is in contrast to uniform (weak and strong) laws of large numbers, central limit theorems, and so on which almost invariably invoke a uniform integrability condition on some higher moment (and indeed, such conditions have been shown to be necessary \citep{waudby2024distribution,waudby2025nonasymptotic}).
Let us now provide a distribution-uniform bound on the expected rejection time.
\begin{theorem}\label{theorem:uniform stopping time}
Let $\Qcalu \subseteq \Qcal$ be a subset of alternative distributions for which the $s^\tth$ moment of the log-wealth increments of $\eproc(\lambda_Q)$ is uniformly bounded for some $s > 2$:
\begin{equation}
    \sup_\Qinu \EE_Q \left \lvert \log \left ( \lambda_Q^\trans \bE_1 
    \right ) - \ell_Q^\star \right \rvert^s < \infty
\end{equation}
and the optimal growth rate is $\Qcalu$-uniformly lower-bounded, meaning there exists some $\ubar \ell(\Qcalu) > 0$ such that
\begin{equation}
  \inf_\Qinu \ell_Q^\star \geq \ubar \ell(\Qcalu).
\end{equation}
Then for any $\Pcal$-$e$-process $\eproc$ with sublinear portfolio regret, the expected rejection time of $\tau_\alpha  := \{ n \in \NN : \eproc_n \geq 1/\alpha \}$ can be uniformly upper bounded in the $\alphato$ regime as
\begin{equation}
  \limsup_\alphato \sup_\Qinu \frac{\EE_Q \left [ \tau_\alpha \right ]}{\log(1/\alpha)} \leq \frac{1}{\ubar \ell(\Qcalu)}.
\end{equation}
Finally, if the optimal growth rate is $\Qcalu$-uniformly upper-bounded by $\widebar \ell(\Qcalu)$, then any $\alpha$-correct stopping rule $\tau_\alpha'$ satisfies
\begin{equation}
  \frac{\EE_Q [\tau_\alpha']}{\log(1/\alpha)} \geq \frac{1}{\widebar \ell(\Qcalu)}.
\end{equation}
\end{theorem}
The proof of \cref{theorem:uniform stopping time} can be found in
\cref{proof:uniform stopping time}, and it follows from some simple arguments
once given access to a nonasymptotic bound on the expected stopping time. Finally, we note that the
discussions in \cref{section:corollaries} can all be stated in a
distribution-uniform fashion, and these seem to be the first of their kind in
the literature.

%% file: relatedwork.tex
\section{Further discussion of related work}\label{section:related work}
Here we outline some related work that have studied growth-rate-log-optimality properties and bounds on expected rejection times. We begin by discussing some classical foundations laid by Wald, Robbins, and others, but note that the rest of the discussion does not proceed chronologically.

As noted in the introduction, the field of sequential hypothesis testing goes back to the seminal work of \citet{wald1945sequential,wald1947sequential} and his sequential probability ratio test (SPRT) for simple (i.e., non-composite) and parametric hypotheses. The field saw a surge in interest in the 1960s and 1970s with the work of Robbins, Siegmund, Lai, and Darling, often taking an estimation (rather than testing) perspective through the construction of confidence sequences with explicit use of Ville's inequality \citep{ville1939etude} applied to nonnegative supermartingales \citep{darling1967confidence,robbins1968iterated,robbins1970statistical,robbins1972class,robbins1974expected,lai1976confidence}.\footnote{Confidence sequences are anytime-valid analogues of confidence intervals and are in a certain sense ``dual'' to sequential tests, but they lie outside the scope of this paper. See \citet[Section 5.4]{waudby2020estimating} for a description of this connection.} Moreover, many of these works considered nonparametric conditions but their test supermartingales do not resemble the ones considered here. Instead, one can find, for example, sub-Gaussian exponential supermartingales \citep{robbins1970statistical} which---with the power of hindsight---can be viewed as a strengthening and generalization of the Chernoff method for the derivation of concentration inequalities such as Hoeffding's \citep{hoeffding_probability_1963}. See \citet*{howard_exponential_2018,howard2018uniform} for a comprehensive study of this technique and its applications to a wide array of settings.

Particularly relevant to the present paper, \citet[\S 4]{wald1945sequential} considered bounds on the expected rejection time of the SPRT that scale inversely with the expected log-likelihood ratio. From this perspective, the maximum expected log-wealth increment $\ell_Q^\star$ can be thought of as a nonparametric and composite analogue of this expected ratio.

\citet{dixit2024anytime} show that in certain problems for which maximum likelihood estimators can be used (including nonparametric settings such as testing log-concavity; see \citep{wasserman2020universal,gangrade2023sequential,dunn2025universal}), $e$-processes built through predictive recursion enjoy almost sure exponential growth properties and they provide explicit rates of almost sure convergence under certain regularity conditions. However, the test supermartingales they consider do not include the general setting alluded to in \eqref{eq:intro-general-test-sm} nor the special cases considered in \cref{section:corollaries}.

However, several authors have considered the use of test (super)martingales that fall under the general representation of \eqref{eq:intro-general-test-sm} for problems with composite nulls and alternatives such as bounded mean testing \citep{hendriks2021test,waudby2020estimating,orabona2021tight,ryu2024confidence,cho2024peeking}, two-sample testing \citep{shekhar2023nonparametric,podkopaev2023sequentialTwoSample}, independence testing \citep{podkopaev2023sequentialTwoSample,podkopaev2023sequentialKernelized}, and difference-in-means testing of bounded tuples \citep{chugg2023auditing,chen2025optimistic}.

In the pursuit of bounded mean testing, \citet{waudby2020estimating} use precisely the test martingale of \cref{section:bounded two sided} and suggest (among other methods) the use of ONS \citep{hazan2007logarithmic,cutkosky2018black}, an online learning algorithm satisfying a certain notion of regret, but not portfolio regret. Moreover, they do not provide formal results relating to asymptotic almost sure log-optimality nor expected rejection times.
In different contexts of two-sample and independence testing, \citet*{shekhar2023nonparametric,podkopaev2023sequentialKernelized}, and \citet{podkopaev2023sequentialTwoSample} also consider the use of ONS and additionally provide analyses pertaining to almost sure growth rates, but these satisfy a weaker notion of optimality than that discussed in this paper. Moreover, they do not provide bounds on expected rejection times under this strategy. Similar analyses and discussions can be found in \citet*{chugg2023auditing} and \citet{chen2025optimistic} in the context of bounded mean testing, but with the addition of an upper bound on the expected rejection time. However, these bounds are not sharp in general and they can be derived from ours (at least qualitatively; see \cref{section:difference-in-means}). Notably, \citet{chen2025optimistic} do not use ONS but rely on optimistic interior point methods instead. This approach can improve on ONS under some special conditions that they discuss, but it does not satisfy the portfolio regret bound that we require, nor does it yield growth rate or rejection time bounds that are substantially sharper than those of \citet{chugg2023auditing} in general. None of the aforementioned works have exact matching lower and upper bounds for growth rates nor expected rejection times in the settings they consider. Moreover in the context of expected rejection time bounds, once taking $\alphato$, the aforementioned bounds do not match the lower bounds provided in \cref{section:stopping time,section:corollaries}.

The problem of best-arm identification in multi-armed bandits is intimately related to sequential hypothesis testing of the mean. More specifically, sequential tests are often used as tools in pursuit of a different end---namely to identify an optimal arm among several options as quickly as possible. Several prior works in that literature explicitly analyze the $\log(1/\alpha)$-rescaled expected rejection time $\EE_Q[\tau_\alpha] / \log (1/\alpha)$ of their underlying sequential tests in the $\alphato$ regime. See \citet*{garivier2016optimal}, \citet*{kaufmann2016complexity}, \citet*{kaufmann2021mixture}, \citet*{agrawal2020optimal} and \citet*{agrawal2021optimal} for some examples. The former three focus on exponential families and related settings such as sub-Gaussianity, while the latter considers heavy-tailed nonparametric distributions. All of these works obtain matching lower and upper bounds on $\EE_Q[\tau_\alpha] / \log(1/\alpha)$ as $\alphato$ with exact constants for their settings. The recent followup work of \citet{agrawal2025stopping} studies and provides lower bounds on expected stopping times for general null hypotheses against point alternatives. They similarly show that certain power one tests can attain these lower bounds in the $\alphato$ regime.

Returning to the context of testing the mean of bounded random variables, and as alluded to in the discussion surrounding \cref{corollary:adaptive optimality of UP and OJ}, \citet{orabona2021tight} were the first to suggest the use of Cover's universal portfolio algorithm as a betting strategy for test martingale construction and they show that this has excellent empirical performance when the tests are inverted to form confidence sequences. As a sequential testing-focused successor to regret-based concentration inequalities appearing in \citet[Section 7.2]{jun2019parameter} (see also \citet{rakhlin2017equivalence}), \citet{orabona2021tight} provide sharp bounds on the regret of the universal portfolio algorithm and use it to derive $e$-processes. They show that these $e$-processes retain much of the aforementioned empirical performance while being computationally preferable. Similar approaches can be found in \citet{ryu2024confidence,ryu2024gambling}, \citet*{ryu2025improved}, and \citet*{clerico2025confidence}. However, these works do not contain formal results on the asymptotic almost sure log-optimality of this approach in the sense of \cref{definition:log-optimality}, nor on the expected rejection time.
From an algorithmic perspective, \cref{corollary:adaptive optimality of UP and OJ} can be thought of as generalizing \citet{orabona2021tight} beyond bounded mean testing, but our main contributions are not algorithmic. Instead, we are focused on exactly these types of formal optimality results. From this theoretical perspective, our results can be viewed as showing that the algorithms proposed by \citet{orabona2021tight} and related methods \citep{ryu2024confidence,ryu2024gambling,ryu2025improved} (and any other methods satisfying sublinear portfolio regret more generally) enjoy strong optimality properties in the bounded means testing problem.

\citet{shekhar2023near} studied the properties of certain confidence intervals and sequences (including those suggested by \citet{orabona2021tight}) in terms of their effective widths and show that they are in a sense optimal, but they did not consider the log-optimality and rejection time desiderata specific to testing that we focus on here.
Our motivations have some overlap with those of \citet{shekhar2023near}, but from a perspective of testing rather than estimation; moreover, our results and proof techniques are substantially different from theirs.

%% file: conclusions.tex
\section{Summary and discussion}

We have considered a class of sequential hypothesis testing problems
that encompasses several problems commonly studied in the literature---including bounded mean testing, difference-in-mean testing, two-sample and independence testing with oracle witness functions, $e$-backtesting, among others. We showed that test supermartingales, $e$-processes, and stopping rules cannot have faster growth rates than $\ell_Q^\star$ nor expected rejection times smaller than $\log(1/\alpha) / \ell_Q^\star$. Conversely, we show that any $e$-process with sublinear portfolio regret attains both of these unimprovable bounds, and we highlight that Cover's universal portfolio algorithm can be used to construct one such $e$-process \citep{cover1991universal,cover1996universal}. 

Overall our results yield a constructive approach to obtaining optimal $e$-process-based sequential tests that are applicable to a wide range of settings. To the best of our knowledge, the matching lower and upper bounds on almost sure growth rates and expected rejection times are the first to appear in the sequential hypothesis testing literature for the null considered in \eqref{eq:null}. Finally, we also demonstrated that these bounds hold in a strictly stronger distribution-uniform sense which appears to open a new line of inquiry in the literature on log-optimality of $e$-processes and sequential tests. An interesting future direction would be to extend the almost sure growth rates and expected rejection times results to settings in which the random vectors of $e$-values $\infseqn{\bE_n}$ are not \iid{}, but can be drawn from arbitrary distributions.

%% file: acks.tex
\subsection*{Acknowledgments}
  The authors are grateful to Shubhada Agrawal and Aaditya Ramdas for comments on an early draft of this work. We also thank Liviu Aolaritei, Ben Chugg, Siva Balakrishnan, Noah Liniger, Sasha Rakhlin, and Ryan Tibshirani for helpful conversations. IW-S acknowledges support from the Miller Institute for Basic Research in Science.  MJ acknowledges support from the Chair ``Markets and Learning,'' supported by Air Liquide, BNP PARIBAS ASSET MANAGEMENT Europe, EDF, Orange and SNCF, sponsors of the Inria Foundation, and
  funding by the European Union (ERC-2022-SYG-OCEAN-101071601).
  Views and opinions expressed are those of the author(s) only and do not necessarily reflect those of the European Union or the European Research Council Executive Agency. Neither the European Union nor the granting authority can be held responsible for them.

%% file: appendix.tex
\section{Proofs of the main results}
In this section we present the proofs of the main results. They appear in the same order as their respective results in the main text.

\subsection{Proof of \cref{theorem:adaptive optimality}}
\label{proof:adaptive optimality-specific}

The proof proceeds in four steps. We begin by deriving a time-uniform
concentration inequality for the absolute difference between $\infseqn{W_n}$ and a
process with oracle access to $\lambda_Q$. We then prove properties
$(i)$-$(iii)$ where the aforementioned time-uniform concentration inequality
is used.

\begin{proof}[Proof of the time-uniform concentration inequality]
    Fix $\Qin$, $\eps > 0$, and define the process with oracle access to
    $\lambda_Q$ as: $W_n(\lambda_Q) \coloneqq \prod_{i=1}^n 
    \lambda_Q^\top \bE_i$. We first provide an upper bound on the
    following probability:
  \begin{equation}
    \Q \left ( \supkm \left ( a_k^{-1} \left \lvert \log W_k - \log
    W_k(\lambda_Q) \right \rvert \right ) \geq \eps \right
    ).
    \label{eq:proof-special-case-quantity-to-bound}
  \end{equation}
  Note that we can upper bound \eqref{eq:proof-special-case-quantity-to-bound} by
  \begin{align}
    \eqref{eq:proof-special-case-quantity-to-bound} \leq \underbrace{\Q \left
    ( \supkm \left (
    a_k^{-1}\log(W_k / W_k(\lambda_Q)) \right ) \geq \eps \right )}_{(+)}
    + \underbrace{\Q \left ( \supkm \left ( a_k^{-1}\log(W_k(\lambda_Q)
    / W_k) \right ) \geq \eps \right )}_{(-)},
  \end{align}
  where the terms $(+)$ and $(-)$ are probabilities in terms of events that $a_k^{-1}\log(W_k / W_k(\lambda_Q))$ are $\eps$-far apart for any $k \geq m$ from above or below, respectively. Looking first to $(+)$, we have
  \begin{align}
    (+) &= \Q \left ( \exists k \geq m : a_k^{-1} \log \left ( W_k / W_k(\lambda_Q) \right ) \geq \eps \right )\\
        &= \Q\left ( \exists k \geq m : W_k / W_k(\lambda_Q) \geq \exp \left \{ a_k \eps \right \} \right ) \\
        &\leq \sum_{k=m}^\infty \frac{\EE_Q (W_k / W_k(\lambda_Q))}{\exp \{a_k \eps \}} \leq \sum_{k=m}^\infty \exp \{-a_k \eps \},
  \end{align}
  where the first inequality follows from a union bound and Markov's inequality,
  and the second follows from the numeraire property of $\lambda_Q$ as
  stated  in \cref{lemma:kkt-supermartingale-doob}.
  Turning now to $(-)$, we have that
  \begin{align}
    (-) &= \Q \left ( \supkm \left( a_k^{-1}[\log(W_k(\lambda_Q)) - 
            \log(W_k)]) \right ) \geq \eps \right ) \\ 
    &\leq \Q \left ( \supkm \left (  a_k^{-1}
    \left [\log W_k^\Max - \log W_k \right ] \right ) \geq \eps \right ) \\
    &\leq \Q \left \{ \supkm a_k^{-1} r_k \geq
                                   \eps \right \} \\
    &= \1 \left \{ \supkm a_k^{-1} r_k \geq
                                   \eps \right \},
  \end{align}
  where the first inequality follows from
  $\eproc_k^\Max := \sup_{\lambda \in \simplex} \log \eproc_k (\lambda) \geq
  \log \eproc_k(\lambda_Q)$ by definition, the second inequality follows
  from the portfolio regret bound, and the final equality follows by noting 
  that all of the terms in the previous probability are deterministic.
  Putting the two inequalities for $(+)$ and $(-)$ together, we have that
  \begin{equation}
    \Q \left ( \supkm \left ( a_k^{-1} \left \lvert \log W_k - \log
    W_k(\lambda_Q) \right \rvert \right ) \geq \eps \right
    ) \leq \sum_{k=m}^\infty \exp \left \{ -a_k \eps
    \right \} + \1 \left \{ \supkm a_k^{-1}r_k \geq \eps \right \},
    \label{eq:proof-log-opt-time-unif-concentration}
  \end{equation}
  which completes the proof of the time-uniform concentration.
\end{proof}

With access to the inequality from \eqref{eq:proof-log-opt-time-unif-concentration}, 
properties $(i)$--$(iii)$ follow with just a few more arguments. We proceed by first
proving property $(ii)$, then $(iii)$, and then $(i)$. 

\begin{proof}[Proof of property $(ii)$.]
  By the inequality in \eqref{eq:proof-log-opt-time-unif-concentration}, it 
  follows that for any $\Qin$ and any $\eps > 0$,
  \begin{equation}
    \Q \left ( \supkm \left \lvert a_k^{-1} \log (\eproc_k /
            W_k(\lambda_Q))
    \right \rvert \geq \eps \right ) \leq \sum_{k=m}^\infty \exp \left \{ -a_k
    \eps \right \} + \1 \left \{ \supkm a_k^{-1} r_k \geq \eps \right \}.
  \end{equation}
  Now by the assumption that $\sum_{k=1}^\infty \exp \left \{ -a_k\eps \right \}
  < \infty$ and $a_k^{-1} r_k \to 0$, we have that the right-hand side of the
  above inequality vanishes as $\mto$. Since for any sequence of random
  variables $\infseqn{Z_n}$, it is the case that $Z_n \to 0$ as $\nto$ with
  $Q$-probability one if and only if $\supkm |Z_k| \to 0$ in $Q$-probability as
  $\mto$, we have that
  \begin{equation}
    a_n^{-1} \log \left ( \eproc_n / W_n(\lambda_Q) \right ) \to 0
    \quad\text{with $Q$-probability one,}
  \end{equation}
  which completes the proof of property $(ii)$.
\end{proof}

\begin{proof}[Proof of property $(iii)$]
      By the strong law of large numbers we have that with 
      $Q$-probability one,
      \begin{equation}
          \frac{1}{n} \log\left(W_n(\lambda_Q)\right) \to
          \EE_{Q}\left[\log\left((\lambda_Q)^\top \bE_1 \right)\right]   
          \equiv \ell_Q^\star \,. 
      \end{equation}
    Appealing to property $(ii)$ we have that
    \begin{equation}
      \frac{1}{n} \log (\eproc_n) = \frac{1}{n} \log(\eproc_n / 
      W_n(\lambda_Q)) + \frac{1}{n} \log(W_n(\lambda_Q)) \to \ell_Q^\star
    \end{equation}
    $Q$-almost surely for every $\Qin$.

    To demonstrate unimprovability of this limit for any other $e$-process
    $\eproc' \in \eprocclass$, we note that $\log (\eproc_n') =
    \log(\eproc_n' / \eproc_n(\lambda_Q)) + \log (\eproc_n(\lambda_Q))$ 
    and hence with $Q$-probability one,
    \begin{equation}
      \limsup_\nto \frac{1}{n} \log (\eproc_n') \leq \underbrace{\limsup_\nto
      \frac{1}{n} \log (\eproc_n' / \eproc_n(\lambda_Q))}_{\leq 0} +
      \underbrace{\limsup_\nto \frac{1}{n} \log(\eproc_n(\lambda_Q))}_{=
      \ell_Q^\star} 
      \label{eq:proof-log-opt-unimprov-decomp} \,.
    \end{equation}
    The first term in the inequality is nonpositive by the numeraire property
    (\cref{lemma:kkt-supermartingale-doob}) combined with the Borel-Cantelli
    lemma:
    \begin{align}
      \forall \delta>0,\quad\sum_{n=1}^\infty \Q \left ( \frac{1}{n}
      \log(\eproc_n' / \eproc_n(\lambda_Q)) \geq \delta \right ) 
      &= \sum_{n=1}^\infty \Q (\eproc_n' / \eproc_n(\lambda_Q) \geq 
      \exp \{ n\delta\} )\\
      &\leq \sum_{n=1}^\infty \EE_{Q}\left[\eproc_n' / \eproc_n(\lambda_Q)
      \right] \exp\left(-n\delta\right) \\
      &\leq \sum_{n=1}^\infty \exp\left(-n \delta \right)  < \infty \,,
    \end{align}
    where the first inequality follows from Markov's inequality and the second
    by \citet[Theorem 15.2.2]{cover1999elements} (see also a generalization in \cref{lemma:kkt-supermartingale-doob}). The second term in the upper bound
    from \eqref{eq:proof-log-opt-unimprov-decomp} is exactly $\ell_Q^\star$ by 
    the strong law of large numbers. Hence, putting these steps together we can
    conclude that
    \begin{equation}
        \limsup_{n \to \infty}\frac{1}{n} \log\left(\eproc_n'\right)
        \leq \ell_Q^\star.
    \end{equation}
    This completes the proof for property $(iii)$.
\end{proof}

\begin{proof}[Proof of property $(i)$]
      Let $\eproc' \in \eprocclass$ be any other $e$-process in $\eprocclass$. 
      Then we have for any $\Qin$,
      \begin{align}
        \log(\eproc_n / \eproc_n') &= \log \left (
        \frac{\eproc_n}{\eproc_n(\lambda_Q)} \cdot
        \frac{\eproc_n(\lambda_Q)}{\eproc_n'}\right ) \\
        &= \log ( \eproc_n / \eproc_n(\lambda_Q)) +
        \log(\eproc_n(\lambda_Q) / \eproc_n') \,.
      \end{align}
      Hence, we have that for any $\Qin$:
      \begin{align}
        \MoveEqLeft{\liminf_\nto \frac{1}{n} \left ( \log(\eproc_n /
        \eproc_n(\lambda_Q)) + \log (\eproc_n(\lambda_Q) / \eproc_n'
        \right )}\\
        \geq\ &\liminf_\nto \frac{1}{n} \log(\eproc_n /
        \eproc_n(\lambda_Q)) +
        \liminf_\nto \frac{1}{n} \log(\eproc_n(\lambda_Q) / \eproc_n') \geq 0,
      \end{align}
      where the final inequality follows from asymptotic equivalence of 
      $\eproc$ and $\eproc(\lambda_Q)$---as shown in property $(ii)$---and 
      from another application of the numeraire property with the Borel-Cantelli
      lemma. This completes the proof for property $(i)$ and hence the proof of
      \cref{theorem:adaptive optimality} altogether.
\end{proof}

 \begin{lemma}[The log-optimal strategy for \eqref{eq:intro-general-test-sm} is a numeraire portfolio (\citet{long1990numeraire})]\label{lemma:kkt-supermartingale-doob}
   Let $\eproc_n$ be any test $\Pcal$-supermartingale of the form \eqref{eq:intro-general-test-sm}. Let $\Qin$ be an arbitrary alternative distribution and let $\eproc_n(\lambda_Q) = \prod_{i=1}^n \lambda_\Q^\trans \bE_i$. Then the suboptimality wealth ratio $\infseqn{S_n^\brackQ}$ given by
   \begin{equation}
     S_n^\brackQ := \eproc_n / \eproc_n(\lambda_Q)
   \end{equation}
   is a nonnegative $Q$-supermartingale with $\EE_Q[S_1^\brackQ] \leq 1$. It follows from Doob's optional stopping theorem that for an arbitrary stopping time $\tau$,
   \begin{equation}\label{eq:lemma-numeraire}
     \EE_Q [S_\tau^\brackQ] \leq 1.
   \end{equation}
 \end{lemma}
 \cref{lemma:kkt-supermartingale-doob} is a strengthening of 
 \citet[Theorem 15.2.2]{cover1999elements} since the latter can be viewed as saying that
 $\EE_Q[S_n^\brackQ] \leq 1$ for all \emph{fixed and nonrandom} $n \in \NN$
 while \cref{lemma:kkt-supermartingale-doob} says that the same holds even when
 $n \equiv \tau$ is an arbitrary data-dependent stopping time $\tau$. %
 
 The proof of \cref{lemma:kkt-supermartingale-doob} uses the necessary and
 sufficient KKT conditions satisfied by $\lambda_Q$ but applies them in a
 conditional form; details are provided below. We remark that a technically accurate
 conclusion of \cref{lemma:kkt-supermartingale-doob} is that $S_n^\brackQ$ forms
 a test $Q$-supermartingale, but this is not necessarily a practically relevant
 interpretation since we are not interested in \emph{testing} the alternative
 $Q$, nor can we construct $S_n^\brackQ$ since it depends on the log-optimal
 strategy $\lambda_Q$ which is unknown for practical purposes. Instead,
 \cref{lemma:kkt-supermartingale-doob} serves as a technical device in the
 proof of \cref{theorem:adaptive optimality} and in many of the results to come.

\begin{proof}
  We aim to show that $S_n$ forms a nonnegative $\Qcal$-supermartingale with
  mean $\EE_Q[S_1] \leq 1$ for all $\Qin$. Indeed, nonnegativity is obtained by
  construction and $\EE_Q[S_1] \leq 1$ by \citet[Theorem
  15.2.2]{cover1999elements} so it remains to show that $S_n$ forms a
  $Q$-supermartingale for all $\Qin$. Writing out the conditional expectation of
  $S_n$, we have
  \begin{align}
    \EE_Q[S_n \mid \Fcal_{n-1}] &= \EE_Q \left [ \prod_{i=1}^n
    \frac{\lambda_i^\trans \bE_i}{\lambda_Q \bE_i} \Bigm \vert
    \Fcal_{n-1} \right ]\\
    &= \underbrace{\frac{\eproc_{n-1}}{\eproc_{n-1}^\star}}_{S_{n-1}}
    \underbrace{\EE_Q \left [ \frac{\lambda_n^\trans \bE_n}{\lambda_Q^\trans \bE_n} \Bigm \vert
    \Fcal_{n-1} \right ]}_{(\star)},
  \end{align}
  and it remains to show that $(\star) \leq 1$ $Q$-almost surely. Using the fact
  that the vectors of $e$-values $\infseqn{\bE_n}$ are
  identically distributed and by definition of $\lambda_Q$, we have
  \begin{align}
    \lambda_Q %
                    &= \argmax_{\lambda \in \simplex} \EE_Q \left [ \log
                    (\lambda^\trans \bE_n) \right ],
  \end{align}
  and thus $\lambda_Q$ satisfies the following KKT conditions (see
  \citet[Theorem 15.2.1]{cover1999elements}):
  \begin{equation}
    \EE_Q \left [ \frac{E_n^\brackj}{\lambda_Q^\trans \bE_n} \right ] \leq 1 \quad\text{for $j \in
    [0:d]$}.
  \end{equation}
  Appealing to the independence of the vectors $\infseqn{\bE_n}$, 
  we have that this inequality holds conditionally on $\Fcal_{n-1}$ for each $j \in [0:d]$,
  \begin{equation}
    \EE_Q \left [ \frac{E_n^\brackj}{\lambda_Q^\trans \bE_n} \Bigm \vert \Fcal_{n-1} \right ] = \EE_Q
    \left [ \frac{E_n^\brackj}{\lambda_Q^\trans \bE_n} \right ] \leq 1.
  \end{equation}
  Using the assumption that $\infseqn{\lambda_n}$ is a $\simplex$-valued predictable
  sequence---i.e., $\lambda_n$ is $\Fcal_{n-1}$-measurable---we have that
  \begin{align}
    \MoveEqLeft{\EE_Q \left [ \frac{\lambda_n^\trans
    \bE_n}{\lambda_Q^\trans \bE_n} \Bigm \vert \Fcal_{n-1} \right ]}
    = \sum_{j=0}^{d} \lambda_n^{(j)}
    \EE_Q\left[\frac{E_n^{(j)}}{\lambda_Q^\trans \bE_n} \Bigm \vert
    \Fcal_{n-1}\right]
    \leq \sum_{j=0}^d \lambda_n^{(j)} = 1\,.
  \end{align}
  Hence $S_n$ forms a nonnegative $Q$-supermartingale with mean $\EE_Q[S_1] \leq 1$ for every $\Qin$, completing the proof.
\end{proof}

\subsection{Proof of \cref{proposition:expected stopping time of constant betting strategy}}
\label{proof:expected stopping time of constant betting strategy}
The proof consists of two parts. In the first part we provide the proof of the
nonasymptotic lower bound on the expected rejection time of rebalanced
portfolios. In the second part we prove an upper bound for the expected
rejection time of rebalanced portfolios, and argue that upper and lower bounds
are match in the $\alpha \to 0+$.

\begin{proof}[Proof of the nonasymptotic lower bound.]
Let us begin by showing the lower bound on $\EE_Q[\tau_\alpha]$. Let $\lambda
\in \simplex$ be an arbitrary constant betting strategy 
and let $\tau \equiv \tau_\alpha := \inf \{ n \in \NN : \eproc_n(\lambda) \geq
1/\alpha \}$ be the stopping time of its associated test: $W_n(\lambda)
\coloneqq \prod_{i=1}^n \lambda^\trans \bE_i$.
We assume that $\EE_Q[\tau] < \infty$, since otherwise the lower bound is immediate.
Notice that by definition of $\tau$, we have that for any $Q \in \Qcal$,
  \begin{equation}
    \log \eproc_\tau(\lambda) \geq \log(1/\alpha)\quad\text{$Q$-almost surely.}
  \end{equation}
  Applying Wald's identity, we have that
  \begin{align}
    \EE_Q \left [ \log(W_\tau(\lambda)) \right ] &= \EE_Q \left [ \sum_{i=1}^\tau \log \left ( \lambda_\Q^\trans \bE_i \right ) \right ]\\
                                                 &= \EE_Q [\tau] \EE_Q \left [ \log \left ( \lambda_\Q^\trans \bE_1 \right ) \right ]\\
    &= \EE_Q [\tau] \ell_Q(\lambda),
  \end{align}
  and once combined with the inequality $\log \eproc_\tau(\lambda) \geq \log(1/\alpha)$, we have that
  \begin{equation}
    \EE_Q[\tau] \geq \frac{\log(1/\alpha)}{\ell_Q(\lambda)},
  \end{equation}
  which yields the lower bound of \cref{proposition:expected stopping time of constant betting strategy}.
\end{proof}

\begin{proof}[Proof of the nonasymptotic upper bound]
  We now derive an upper bound on the expected time to rejection.
  Fix a $\Qin$ and let $\delta \in (0, 1)$ be an arbitrary constant. Define
  the quantities $\eps$ and $m$ given by
  \begin{equation}
    \eps := \frac{\delta}{1+\delta}\cdot \ell_Q(\lambda)\quad\text{and}\quad
    m := \left \lceil \frac{\log(1/\alpha)}{\ell_Q(\lambda)-\eps} \right \rceil = \left \lceil \frac{(1+\delta)\log(1/\alpha)}{\ell_Q(\lambda)} \right \rceil. 
  \end{equation}
  Writing out the expected stopping time under $Q$, we have
  \begin{align}
    \EE_Q[\tau] &= \sum_{k=1}^\infty \Q (\tau  \geq k)\\
    &\leq  m + \sum_{k=m}^\infty \Q \left ( \tau \geq k \right ) \\
    &= m + \sum_{k=m-1}^\infty \Q \left ( \tau > k \right ) \\
    &\leq 1 + m + \sum_{k=m}^\infty \Q \left ( \tau > k \right ) \,,
    \label{eq:proof-constant-rebalanced-port-exp-time-to-rejection-upper-bound}
  \end{align}
  where the first inequality follows from upper bounding the first $m$
  probabilities in the series by one, and in the second inequality we upper
  bound the $m-1$th probability by one.

  Recall the definition of the stopping time $\tau \equiv \tau_{\alpha} 
  \coloneqq \left\{n \in \NN ~|~ W_n(\lambda) \geq 1/\alpha\right\}$. 
  This allows us to rewrite the event inside the
  probability from
  \eqref{eq:proof-constant-rebalanced-port-exp-time-to-rejection-upper-bound}
  and upper bound the resulting probability as:
  \begin{align}
      \eqref{eq:proof-constant-rebalanced-port-exp-time-to-rejection-upper-bound}
    &= 1 + m + \sum_{k=m}^\infty \Q \left ( W_k < 1/\alpha \right ) \\
    &= 1 + m + \sum_{k=m}^\infty \Q \left ( k^{-1} \log W_k < k^{-1}\log (1/\alpha) \right ) \\
    &\leq 1 + m + \sum_{k=m}^\infty \left [ \Q \left( \ell_Q(\lambda) - \eps <
    k^{-1} \log(1/\alpha) \right) + \Q \left ( |k^{-1} \log W_k -
    \ell_Q(\lambda) | \geq \eps \right ) \right ] \\
    &\leq 1 + m + \sum_{k=m}^\infty \left [ \1 \{ \ell_Q(\lambda) - \eps <
    m^{-1} \log(1/\alpha) \} + \Q \left ( |k^{-1} \log \eproc_k -
    \ell_Q(\lambda) | \geq \eps \right ) \right ] \,,
    \label{eq:proof-constant-rebalanced-port-exp-time-to-rejection-union-and-indicator}
    \end{align}
    where the first inequality follows from a union bound, and in the second
    inequality we have upper bounded the probability by an indicator and have
    used the fact that $k^{-1} \leq m^{-1}$ whenever $k \geq m$. Appealing to
    the definition of $m := \left \lceil \log(1/\alpha) / (\ell_Q(\lambda) -\eps
    ) \right \rceil $ allows us to simplify
    \eqref{eq:proof-constant-rebalanced-port-exp-time-to-rejection-union-and-indicator}
    as follows:
    \begin{equation}
        \eqref{eq:proof-constant-rebalanced-port-exp-time-to-rejection-union-and-indicator}
    = 1 + m + \sum_{k=m}^\infty  \Q \left ( |k^{-1} \log W_k -
    \ell_Q(\lambda)| \geq \eps \right ) \,.
    \end{equation}
    Now, appealing to
    \cref{lemma:concentration-via-nemirovski} and plugging in the values of
    $\eps$ and $m$ in terms of $\delta$, $\ell_Q(\lambda)$, and $\alpha$, we
    have the following upper bound:
    \begin{align}
    \EE_Q[\tau] &\leq 2 + m + \frac{1}{\eps^{s}}\sum_{k=m+1}^\infty \frac{\EE_Q
    \left \lvert \log(\lambda^\trans \bE_1) - \ell_\Q(\lambda) \right \rvert^s}{k^{s/2}}\\
                &\leq 2 + m + \frac{\rho_s(\ell_Q(\lambda))}{(s/2-1)} \frac{1}{\eps^{s}m^{s/2-1}} \\
                &\leq 3 + \frac{(1+\delta) \log (1/\alpha)}{\ell_Q(\lambda)} + \frac{\rho_s(\ell_Q(\lambda))}{(s/2 - 1)} \cdot \frac{(1+\delta)^s}{\delta^s \ell_Q(\lambda)^s} \cdot \frac{\ell_Q(\lambda)^{s/2 - 1}}{(1+\delta)^{s/2 - 1} \log (1/\alpha)^{s/2-1}}\\
      &= 3 + \frac{(1+\delta) \log (1/\alpha)}{\ell_Q(\lambda)} + \frac{\rho_s(\ell_Q(\lambda))}{(s/2 - 1)} \cdot \frac{(1+\delta)^{1 + s/2}}{\delta^s \ell_Q(\lambda)^{1 + s/2}} \cdot \frac{1}{\log (1/\alpha)^{s/2-1}}
  \end{align}
  Now, taking a limit supremum as $\alpha \to 0^+$---in the sense that $\limsup_{\alpha \to 0^+} f(\alpha)$ for a function $f(\cdot)$ is simply $\limsup_{\beta \to \infty} f(1/\beta)$---we have that
  \begin{equation}
    \limsup_{\alpha \to 0^+} \frac{\EE_Q[\tau]}{\log(1/\alpha)} \leq \frac{1+\delta}{\ell_Q(\lambda)}.
  \end{equation}
  Since this holds for any $\delta \in (0, 1)$, it holds with $\delta = 0$, and once combined with the lower bound that was derived previously, the limit supremum can be replaced with a limit and we have the following equality:
  \begin{equation}
    \lim_{\alpha \to 0^+}\frac{\EE_Q[\tau]}{\log(1/\alpha)} = \frac{1}{\ell_Q(\lambda)}, 
  \end{equation}
  which completes the proof.
\end{proof}

\begin{lemma}[A Chebyshev-Nemirovski concentration inequality]\label{lemma:concentration-via-nemirovski}
  Let $X_1, X_2, \dots, X_n$ be independent random variables with mean zero for which $\max_{i \in [n]}\EE_P |X_i|^{s} < \infty$ for some $s \geq 2$. Then for any $\eps > 0$,
  \begin{equation}\label{eq:nemirovski-lemma-main-equation}
    \P \left ( \left \lvert \frac{1}{n} \sum_{i=1}^n X_i  \right \rvert \geq \eps \right ) \leq \frac{1}{\eps^{s} n^{s/2}} \max_{i \in [n]}\EE_P |X_i|^s.
  \end{equation}
\end{lemma}
\begin{proof}
  Throughout, for a random variable $Z$, let $\| Z \|_{L_s(P)} = \left ( \EE_P \left \lvert Z \right \rvert^s \right )^{1/s}$ be the usual $L_s(P)$ norm. Then we can upper bound the probability in \eqref{eq:nemirovski-lemma-main-equation} directly as
  \begin{align}
    \P \left ( \left \lvert \frac{1}{n} \sum_{i=1}^n X_i  \right \rvert \geq \eps \right ) &= \P \left ( \left \lvert \frac{1}{n} \sum_{i=1}^n X_i  \right \rvert^{s} \geq \eps^{s} \right )\\
            &\leq \frac{1}{\eps^{s} n^{s}} \EE_P \left \lvert \sum_{i=1}^n X_i \right \rvert^{s} \label{eq:nemirovski-lemma-proof-markov}\\
            &= \frac{1}{\eps^{s} n^{s}} \left \| \sum_{i=1}^n X_i \right \|_{L_s(P)}^{s} \\
    \verbose{
            &= \frac{1}{\eps^{s} n^{s}} \left ( \left \| \sum_{i=1}^n X_i \right \|_{L_s(P)}^{2} \right )^{s / 2}\\
    }
            &\leq \frac{1}{\eps^{s} n^{s}} \left ( \sum_{i=1}^n \| X_i\|_{L_s(P)}^2\right )^{s / 2}\label{eq:nemirovski-lemma-proof-nemirovski}\\
    \verbose{
            &= \frac{1}{\eps^{s} n^{s}} \left ( n \| X_1\|_{L_s(P)}^2\right )^{s / 2}\\
            &= \frac{n^{s/2}}{\eps^{s} n^{s}} \left ( \| X_1\|_{L_s(P)}^2\right )^{s / 2}\\
    }
            &\leq \frac{1}{\eps^{s} n^{s/2}} \left ( \max_{i \in [n]}\| X_i\|_{L_s(P)}^2\right )^{s / 2}\\
            &= \frac{1}{\eps^{s} n^{s/2}} \max_{i \in[n]}\EE_P |X_i|^s,
  \end{align}
  where \eqref{eq:nemirovski-lemma-proof-markov} follows from Markov's inequality while \eqref{eq:nemirovski-lemma-proof-nemirovski} follows from an inequality due to \citet{nemirovski2000topics} (see \citet[Theorem 2.2]{dumbgen2010nemirovski}). This completes the proof.
\end{proof}

\subsection{Proof of \cref{proposition:expected stopping time oracle}}
\label{proof:expected stopping time oracle}

\begin{proof}
  Fix $\alpha\in(0, 1)$ and $\Qin$ and let $\tau_\alpha$ be any stopping rule that satisfies
  \begin{equation}
    \sup_\Pin \P \left ( \tau_\alpha < \infty \right ) \leq \alpha.
  \end{equation}
  Now fix $\Pin$. We begin by showing that
  \begin{equation}
    \frac{\EE_\Q[\tau_\alpha]}{\log(1/\alpha)} \geq \frac{1}{\KL(\Q \| \P)}.
  \end{equation}
  Assume that $\EE_\Q[\tau_\alpha] < \infty$ and $\KL(\Q \| \P) < \infty$, otherwise the inequality is trivial. The stopping rule $\tau_\alpha$ satisfies $\P(\tau_\alpha < \infty) \leq \alpha$, and hence by \citet[Eq.~(4.81)]{wald1945sequential}, it holds that
  \begin{equation}\label{eq:expected stopping time lower bound fixed P fixed Q}
    \EE_\Q [\tau_\alpha] \geq \EE_\Q[\tau_\alpha(\Q, \P)],
  \end{equation}
  where $\tau_\alpha(\Q, \P) := \inf\{n \in \NN : \prod_{i=1}^n \dQ(\bE_i) / \dP(\bE_i) \geq 1/\alpha \}$. Now, by Wald's identity,
  \begin{align}
    \EE_\Q[\tau_\alpha(\Q, \P)] \EE_\Q[\dQ(\bE_1) / \dP(\bE_1)] &= \EE_\Q\left [\sum_{i=1}^{\tau_\alpha(\Q, \P)}\log(\dQ(\bE_i) / \dP(\bE_i))\right ] \\
    &\geq \log(1/\alpha).
  \end{align}
  Rearranging and applying \eqref{eq:expected stopping time lower bound fixed P fixed Q}, we have that $\EE_\Q[\tau_\alpha] / \log(1/\alpha) \geq 1/ \KL(\Q \| \P)$. Since $P$ was arbitrary, we have that
  \begin{equation}
    \frac{\EE_\Q[\tau_\alpha] }{\log(1/\alpha)} \geq \frac{1}{\inf_\Pin \KL(\Q \| \P)}.
  \end{equation}
  Applying \cref{lemma:klinf identity}, the denominator is $\ell_Q^\star$, so we have that
  \begin{equation}
    \frac{\EE_\Q[\tau_\alpha]}{\log(1/\alpha)}\geq \frac{1}{\ell_Q^\star},
  \end{equation}
  completing the proof.
\end{proof}

The following lemma is well-known in the $e$-value literature but we provide it and its proof here for completeness.
\begin{lemma}\label{lemma:LR is log-opt}
  Let $\Pcal$ be a collection of probability measures and $Q$ a single probability measure. If $E$ is a $\Pcal$-$e$-value,
  \begin{align}
    \EE_Q[\log(E)] \leq \inf_{P \in \Pcal}\KL(Q \| P).
  \end{align}
\end{lemma}
\begin{proof}
  Fix $P\in \Pcal$. We will first demonstrate the inequality without the infimum. The result holds trivially if $\KL(Q \| P) = \infty$ so assume that $Q \ll P$. By Jensen's inequality, we have that
  \begin{align}
    \EE_Q[\log (E)] - \EE_Q[\log(\dd Q / \dd P)] &= \EE_Q \left [ \log \left (  \frac{E}{\dd Q / \dd P} \right ) \right ]\\
                                                 &\leq \log \left ( \EE_Q \left [ \frac{E}{\dd Q / \dd P} \right ] \right ) \\
                                                 &= \log \left ( \EE_P [E] \right ) \leq 0.
  \end{align}
  Since the inequality holds for any $\Pin$, we take infima over $\Pin$ to obtain
  \begin{equation}
    \inf_{\Pin} \EE_\Q \left [ \log(E) \right ] = \inf_\Pin \KL(\Q \| \P).
  \end{equation}
  Since the argument in the left-hand side does not depend on $\P$, the result follows.
\end{proof}

\begin{lemma}[Equality of optimal growth rates and an infimum of KL divergences]\label{lemma:klinf identity}
  Let $\Pcal$ be the collection of probability measures in \eqref{eq:null}. Let $Q$ be an alternative distribution. It holds that
  \begin{equation}
    \ell_Q^\star \equiv \max_{\lambda \in \simplex} \EE_\Q \left [ \log \left ( \lambda^\trans \bE_1 \right ) \right ] = \inf_\Pin \KL(\Q \| \P).
  \end{equation}
\end{lemma}
\begin{proof}
  Invoking \cref{lemma:LR is log-opt} and noting that $\lambda^\trans \bE_1$ is a $\Pcal$-$e$-value for every $\lambda \in \simplex$, we have
  \begin{equation}
    \max_{\lambda \in \simplex}\EE_Q[\log(\lambda^\trans \bE_1)] \leq \inf_{P \in \Pcal} \KL(\Q \| \P).
  \end{equation}
  We now focus on the reverse inequality. Let $\lambda_Q \in \argmax_{\lambda \in \simplex}\EE_Q[\log(\lambda^\trans \bE_1)]$. Note that $\lambda_Q^\trans \bE_1 > 0$ with $Q$-probability one (otherwise $\EE_Q[\log(\lambda_Q^\trans \bE_1)] = -\infty$ but this cannot be since $E_1 = 1$). Define the measure $P^\star$ for any Borel $A$ as
  \begin{equation}
    P^\star(A) = \EE_Q \left [ \frac{\1 \{ \bE_1 \in A \}}{\lambda_Q^\trans \bE_1} \right ] + \left ( 1 - \EE_\Q \left [ \frac{1}{\lambda_Q^\trans \bE_1} \right ]\right ) \1\{ x_1 \in A \}, 
  \end{equation}
  where $x_1$ is $(1, 0, \dots, 0)^\trans$. Since $\lambda_Q$ is the numeraire portfolio, $\EE_Q[(\lambda_Q^\trans \bE_1)^{-1}] \leq 1$ so $\P^\star$ is clearly a sub-probability measure. It is in fact a probability measure since $\P^\star(\{1 \} \times \RR_+^{d}) = \EE_Q[(\lambda_Q^\trans \bE_1)^{-1}] + (1-\EE_Q[(\lambda_Q^\trans \bE_1)^{-1}]) = 1$. Moreover, $\P^\star$ is a probability measure in $\Pcal$ since for any $j \in \{1,\dots, d\}$,
  \begin{align}
    \EE_{\P^\star} \left [ E_1^\brackj \right ]  &= \int \frac{E_1^\brackj}{\lambda_Q^\trans \bE_1} \dQ + \left ( 1-\EE_Q \left [ \frac{1}{\lambda_Q^\trans \bE_1} \right ] \right ) \underbrace{x_1^\brackj}_{=0} \\
    &= \EE_Q \left [ \frac{E_1^\brackj}{\lambda_Q^\trans \bE_1} \right ] \leq 1.
  \end{align}
  We now show that $\lambda_Q^\trans \bE_1 = \dQ / \dP^\star$. First, note that
  \begin{align}
    1 = \EE_\Q \left [ \frac{\lambda_\Q^\trans \bE_1}{\lambda_\Q^\trans \bE_1} \right ] &= \lambda_\Q^\brackzero\EE_\Q \left [ \frac{1}{ \lambda_\Q^\trans \bE_1} \right ] + \sum_{j=1}^d \lambda_\Q^\brackj\EE_\Q \left [ \frac{ E_1^\brackj}{\lambda_\Q^\trans \bE_1} \right ],
  \end{align}
  and hence if $\EE_\Q[(\lambda_\Q^\trans \bE_1)^{-1}] < 1$, then $\lambda_\Q^\brackzero = 0$.
  Therefore, for any Borel $A$,
  \begin{align}
    \int_A \lambda_\Q^\trans \bE_1 \dd \P^\star &= \int_A \frac{\lambda_\Q^\trans \bE_1}{\lambda_\Q^\trans \bE_1} \dQ + \left ( 1-\EE_\Q \left [ \frac{1}{\lambda_\Q^\trans \bE_1} \right ] \right ) \lambda_\Q^\trans x_1 \1 \{ x_1 \in A \} \\
                                                &= \int_A \dQ + \left ( 1-\EE_\Q \left [ \frac{1}{\lambda_\Q^\trans \bE_1} \right ] \right ) \lambda_\Q^\brackzero \1 \{ x_1 \in A \}\\
    &= \int_A \dQ.
  \end{align}
  It follows that $\lambda_\Q^\trans \bE_1 = \dQ / \dP^\star$. Since $\P^\star \in \Pcal$, it follows that
  \begin{equation}
    \inf_\Pin \KL(\Q \| \P) \leq \EE_\Q[\log (\lambda_\Q^\trans \bE_1)],
  \end{equation}
  completing the proof.
\end{proof}

\subsection{Proof of \cref{theorem:stopping time}}
\label{proof:stopping time}
\begin{proof}
    Our focus is to provide an upper bound on the expected rejection time for the
    $\Pe$-process $\infseqn{W_n}$ satisfying the regret bound $\Rcal(W_n) \leq r_n$,
    as the lower bound was provided in 
    \cref{proposition:expected stopping time oracle}.

    Fix a $\Qin$ and let $\delta \in (0,1)$ be an arbitrary constant. We define the
    quantities $\epsilon$ and $m$ given by
    \begin{equation}
        \eps := \frac{\delta}{1+\delta}\cdot \ell_Q^\star \quad\text{and}\quad
        m := \left \lceil \frac{\log(1/\alpha)}{\ell_Q^\star-\eps} \right \rceil =
        \left \lceil \frac{(1+\delta)\log(1/\alpha)}{\ell_Q^\star} \right \rceil. 
    \end{equation}
    Writing out the expected stopping time under $Q$ using the tail-sum formula
    we have:
    \begin{align}
        \EE_Q\left[\tau\right] 
        &= \sum_{k=1}^\infty \Q(\tau \geq k) \\
        &\leq m + \sum_{k=m}^\infty \Q(\tau \geq k) \\
        &= m + \sum_{k=m-1}^\infty \Q(\tau > k) \\
        &\leq 1 + m + \sum_{k=m}^\infty \Q(\tau > k) \,,
    \end{align}
    where in the first inequality we upper bounded the first $m$ probabilities
    by one, and in the second inequality we again upper bounded the $m-1$th 
    probability by one.

    We now substitute in the definition of the stopping time $\tau$ and obtain
    the following identity:
    \begin{align}
        &1 + m + \sum_{k=m}^\infty \Q(\tau > k) \\
        &\quad = 1 + m + \sum_{k=m}^\infty \Q(W_k < 1/\alpha) \\
        &\quad = 1 + m + \sum_{k=m}^\infty \Q\left(k^{-1} \log(W_k) < k^{-1}
            \log(1/\alpha)\right) \,.
    \end{align}
    We now relate the event inside the probability to the expected log-wealth
    increment under the log-optimal bet, $\ell_Q^\star$, and proceed to upper
    bound the series as follows:
    \begin{align}
        &1 + m + \sum_{k=m}^\infty \Q\left(k^{-1} \log(W_k) < k^{-1}
        \log(1/\alpha)\right) \\
        &\quad = 1 + m + \sum_{k=m}^\infty \Q\left(\ell_Q^\star - \epsilon +
        k^{-1}\log(W_k) - \ell_Q^\star + \epsilon < k^{-1}\log(1/\alpha)\right)
        \\
        &\quad \leq 1 + m + \sum_{k=m}^\infty \left[\Q\left( \ell_Q^\star -
        \epsilon < k^{-1} \log(1/\alpha)\right) +
        \Q\left(\ell_Q^\star - k^{-1}\log\left(W_k\right) \geq \epsilon
        \right)\right] \\
        &\quad \leq 1 + m + \sum_{k=m}^\infty \left[\Q\left( \ell_Q^\star -
        \epsilon < m^{-1} \log(1/\alpha)\right) +
        \Q\left(\ell_Q^\star - k^{-1}\log\left(W_k\right) \geq \epsilon
        \right)\right] \\
        &\quad = 1 + m + \sum_{k=m}^\infty \left[\1\left\{ \ell_Q^\star -
            \epsilon < m^{-1} \log(1/\alpha)\right\} +
        \Q\left(\ell_Q^\star - k^{-1}\log\left(W_k\right) \geq \epsilon
        \right)\right] \,,
        \label{eq:proof-stopping-time-upper-bound-sub-linear-regret}
    \end{align}
    where the first inequality follows from a union bound, the second inequality
    follows from the fact that $k^{-1} \leq m^{-1}$ whenever $k \geq m$, and in
    the last identity we have substituted the probability by an indicator as the
    event inside the probability is not random. Plugging in the definitions of
    $m$ and $\epsilon$ allows us to rewrite the series as:
    \begin{align}
        \eqref{eq:proof-stopping-time-upper-bound-sub-linear-regret} 
        &=
        1 + m + \sum_{k=m}^\infty
        \Q\left(\ell_Q^\star - k^{-1}\log\left(W_k\right) \geq \epsilon
        \right) \\
        &\leq
        1 + m + \sum_{k=m}^\infty
        \Q\left(\left\vert\ell_Q^\star -
            k^{-1}\log\left(W_k\right)\right\vert \geq \epsilon
        \right) \\
        &\leq
        1 + m + \underbrace{\sum_{k=m}^\infty
        \Q\left(\left\vert\ell_Q^\star -
            k^{-1}\log(W_k^\star) \right\vert \geq
        \epsilon / 2 \right)}_{(+)}
        + \underbrace{\sum_{k=m}^\infty \Q\left(\left\vert k^{-1}\log(W_k^\star) - 
        k^{-1}\log\left(W_k\right)\right\vert \geq \epsilon / 2 \right)}_{(-)}
        \,.
    \end{align}
    
    We will analyze the series $(+)$ and $(-)$ individually.
    We first focus on $(+)$ and recall the definition $\rho_s(\lambda_Q^\star)
    \coloneqq \EE_Q\left\vert\log\left(\lambda_Q^\trans \mathbf{E_1} -
    \ell_Q^\star \right) \right\vert^s$. Applying the Chebyshev-Nemirovski inequality
    (\cref{lemma:concentration-via-nemirovski}) combined with the assumption
    that the $s$th moment is finite for some $s > 2$ we obtain:
    \begin{align}
        (+) &\leq 1 + \frac{2^s \rho_s(\lambda_Q^\star)}{\epsilon^s}\sum_{k=m+1}^\infty
        \frac{1}{ k^{s/2}} \\
            &\leq 1 + \frac{2^s \rho_s(\lambda_Q^\star)}{\epsilon^s}
            \int_{m}^\infty \frac{1}{k^{s/2}}dk \\
            &= 1 + \frac{2^s \rho_s(\lambda_Q^\star)}{\epsilon^s
            m^{s/2-1}(s/2-1)} \,.
            \label{eq:proof-stopping-time-upper-bound-chebyshev-nemirovski}
    \end{align}
    Substituting in the definition of $\epsilon$ and $m$ we can see how
    \eqref{eq:proof-stopping-time-upper-bound-chebyshev-nemirovski} simplifies
    to:
    \begin{align}
        \eqref{eq:proof-stopping-time-upper-bound-chebyshev-nemirovski}   
        &= 1 + \frac{2^s \rho_s(\lambda_Q^\star)}{\left(\delta \ell_Q^\star / (1
        + \delta)\right)^s \left\lceil\log(1/\alpha)(1 + \delta) / \ell_Q^\star
        \right\rceil^{s/2 - 1}(s / 2 - 1)} \\
        &\leq 1 + \frac{2^s \rho_s(\lambda_Q^\star)}{\left(\delta \ell_Q^\star / (1
        + \delta)\right)^s \left(\log(1/\alpha)(1 + \delta) / \ell_Q^\star
        \right)^{s/2 - 1}(s / 2 - 1)} \\
        &= 1 + \frac{(1 + \delta)^{1 + s / 2}2^s \rho_s(\lambda_Q^\star)}
        {\delta^s (\ell_Q^\star)^{1 + s/2}
        \log^{s/2 - 1}(1/\alpha)(s/2-1)} \,.
    \end{align}

    We now turn our attention to $(-)$, and proceed to upper bound the $k$th
    summand as,
    \begin{align}
        &\Q\left(\left\vert k^{-1}\log(W_k^\star) - k^{-1}\log(W_k)
        \right\vert \geq \epsilon / 2 \right)\\
        &\quad\leq \Q\left(k^{-1}\log(W_k^\star) - k^{-1}\log(W_k) \geq
        \epsilon / 2 \right)
        + \Q\left(k^{-1}\log(W_k) - k^{-1}\log(W_k^\star) \geq \epsilon / 2 
        \right) \\
        &\quad \leq \Q\left(k^{-1} \log(W_k(\lambda^{\max})) - k^{-1}
        \log(W_k) \geq \epsilon / 2 \right) 
        + \Q\left(k^{-1}\log(W_k) - k^{-1}\log(W_k^\star) \geq \epsilon / 2 
        \right) \,,
        \label{eq:proof-stopping-time-upper-bound-numeraire-decomposition}
    \end{align}
    where the first inequality follows from a union bound, and the second
    inequality follows from the fact that with $Q$-probability one, 
    $\log(W_n(\lambda^{\max})) \geq  \log(W_n^\star)$ where $\lambda^{\max} 
    \coloneqq \argmax_{\lambda \in
    \simplex} \sum_{i=1}^n \log(\lambda^\trans \mathbf{E}_i)$. We now appeal to
    the numeraire property of $\lambda_Q$ to upper bound the second
    summand from 
    \eqref{eq:proof-stopping-time-upper-bound-numeraire-decomposition}, and we
    do so as follows:
    \begin{align}
        \Q\left(k^{-1} \log(W_k) - k^{-1} \log(W_k^\star) \geq \epsilon / 2\right)
        &= \Q\left(\log(W_k / W_k^{\star}) \geq k \epsilon / 2\right) \\
        &= \Q\left(W_k / W_k^\star \geq \exp\{k \epsilon / 2\}\right) \\
        &\leq \exp\{- k \epsilon / 2\} \EE_Q\left[W_k / W_k^\star\right] \\
        &\leq \exp\{- k \epsilon / 2\} \,,
    \end{align}
    where the first inequality follows from Markov's inequality and the second
    is due to the numeraire property of the log-optimal portfolio
    $\lambda_Q$ (\cref{lemma:kkt-supermartingale-doob}). 
    Putting the above steps together, and appealing to the
    regret bound $\Rcal(W_n) \leq r_n$ we have that,
    \begin{align}
        (-) &\leq 1 + \sum_{k=m + 1}^\infty \exp\left\{ -k \epsilon / 2\right\}
        + \sum_{k = m}^\infty \1\{r_k / k \geq \epsilon / 2\} \\
            &\leq 1 + \int_{m}^\infty \exp\left\{ -k \epsilon / 2\right\}
        + \sum_{k = m}^\infty \1\{r_k / k \geq \epsilon / 2\} \\
            &= 1 + \frac{2}{\epsilon}\exp\{-m\epsilon / 2\}
        + \sum_{k = m}^\infty \1\{r_k / k \geq \epsilon / 2\} \,.
    \end{align}
    Substituting in the definitions of $\epsilon$ and $m$ we have that $(-)$ can
    be further upper bounded by 
    \begin{align}
        (-) &\leq 1 + \frac{2 (1 + \delta)}{\delta \ell_Q^\star} \exp\left\{-\frac{\delta
        \log(1/\alpha)}{2}\right\} + \sum_{k=m}^\infty \1\left\{r_k / k \geq
    \epsilon / 2\right\}  \\
        &= 1 + \frac{2 (1 + \delta) \alpha^{\delta / 2}}{\delta
        \ell_Q^\star} + \sum_{k=m}^\infty \1\left\{r_k / k \geq \epsilon /
        2\right\} \,.
    \end{align}
    Hence, putting all of the previous steps together we obtain the following
    upper bound on the expected stopping time
    \begin{equation}
    \begin{aligned}
        \EE_Q[\tau] &\leq 4 + \frac{(1 + \delta) \log(1/\alpha)}{\ell_Q^\star}
        + \frac{2 (
        1 + \delta) \alpha^{\delta / 2}}{\delta \ell_Q^\star} \\
        &\qquad+ 
        \frac{(1 + \delta)^{1 + s / 2}2^s \rho_s(\lambda_Q)}
        {\delta^s (\ell_Q^\star)^{1 + s/2}
        \log^{s/2 - 1}(1/\alpha)(s/2-1)} + \sum_{k=m}^\infty \1\left\{r_k / k
        \geq \epsilon / 2\right\} \,.
    \end{aligned}
    \label{eq:proof-expected-stopping-time-nonasymptotic-bound}
    \end{equation}
    We note that the series is finite due to the sub-linearity of $r_k$.
    Dividing both sides by $\log(1/\alpha)$ and taking the limit superior as
    $\alpha \to 0^+$ we have
    \begin{equation}
        \limsup_{\alpha \to 0^+} \frac{\EE_Q[\tau]}{\log(1/\alpha)} \leq 
        \frac{(1 + \delta)}{\ell_Q^\star}
    \end{equation}
    Since $\delta \in (0,1)$ was arbitrary, the inequality holds with a one in
    the numerator. Combining this upper bound with the lower bound from 
    \cref{proposition:expected stopping time oracle} we can conclude that
    \begin{equation}
        \limsup_{\alpha \to 0^+} \frac{\EE_Q[\tau]}{\log(1/\alpha)} = 
        \frac{1}{\ell_Q^\star} \,,
    \end{equation}
    which completes the proof.
\end{proof}

\subsection{Proof of \cref{theorem:uniform log optimality}}\label{proof:uniform log optimality} 

\begin{proof}
    By \eqref{eq:proof-log-opt-time-unif-concentration}, we have that for any $Q \in \Qcal$ and any $\epsilon > 0$,
    \begin{equation}
        \Q \left ( \supkm \left ( a_k^{-1} \left \lvert \log W_k - \log
        W_k(\lambda_Q) \right \rvert \right ) \geq \eps \right
        ) \leq \sum_{k=m}^\infty \exp \left \{ -a_k \eps
        \right \} + \1 \left \{ \supkm a_k^{-1}r_k \geq \eps \right \}.
    \end{equation}
    Now by the assumption that $\sum_{k=1}^\infty \exp \left \{ -a_k\eps \right \} < \infty$ and $a_n^{-1} r_n \to 0$ and noticing that the right-hand side does not depend on any particular $\Q \in \Qcal$, we have that
    \begin{equation}
        \lim_\mto \sup_{Q \in \Qcal} \Q \left ( \supkm \left ( a_k^{-1} \left \lvert \log W_k - \log
        W_k(\lambda_Q) \right \rvert \right ) \geq \eps \right ) = 0,
    \end{equation}
    completing the proof for property $(ii)$.
\end{proof}

\begin{proof}[Proof of property $(iii)$]
    By uniform integrability of the $s^\tth$ moment and invoking the
    distribution-uniform strong laws of large numbers
    \citep{chung_strong_1951,waudby2024distribution}, we have that
      \begin{equation}\label{eq:proof-property-iii-slln}
        \forall \eps > 0,\quad  \lim_\mto \sup_\Qinu \Q \left ( \supkm
        \frac{1}{k^{1/s}} \left \lvert  \log \eproc_k(\lambda_Q) - k\ell_Q^\star
        \right \rvert \geq \eps \right ) = 0.
      \end{equation}
      
      By a triangle inequality, we have for any $\Qin$,
      \begin{align}
        \MoveEqLeft{\Q \left ( \supkm \frac{1}{k^{1/s}} \left \lvert \log
        \eproc_k - k\ell_Q^\star \right \rvert \geq \eps \right )}\\
        & \leq \Q \left ( \supkm \frac{1}{k^{1/s}} \left \lvert \log(
            \eproc_k/\eproc_k(\lambda_Q)) \right \rvert \geq \eps/2
                \right ) + \Q \left ( \supkm \frac{1}{k^{1/s}} \left \lvert
                \log \eproc_k(\lambda_Q) - k\ell_Q^\star \right \rvert \geq \eps/2
            \right ).
      \end{align}
      Using the assumption that $r_n = o(n^{1/s})$ and appealing to property
      $(ii)$, the first term vanishes $\Qcal$-uniformly and hence
      $\Qcalu$-uniformly. Combined with \eqref{eq:proof-property-iii-slln}, we
      have that for any $\eps > 0$,
      \begin{equation}
    \lim_\mto \sup_{Q \in \Qcalu}\Q \left ( \supkm \frac{1}{k^{1/s}} \left
\lvert \log \eproc_k - k\ell_Q^\star \right \rvert \geq \eps \right ) = 0,
      \end{equation}
      which completes the proof of property $(iii)$.
\end{proof}

\begin{proof}[Proof of property $(i)$]
  Let $\eproc' \in \eprocclass$ be any other $e$-process in $\eprocclass$. Then we have for any $\Qin$,
  \begin{align}
    \log(\eproc_n' / \eproc_n) &= \log \left ( \frac{\eproc_n'}{\eproc_n(\lambda_Q)} \cdot \frac{\eproc_n(\lambda_Q)}{\eproc_n} \right ) \\
    &= \log ( \eproc_n' / \eproc_n(\lambda_Q)) + \log(\eproc_n(\lambda_Q) / \eproc_n),
  \end{align}
  and hence we have that for any $Q \in \Qcal$ and any $\eps > 0$,
  \begin{align} 
    \Q \left ( \supkm k^{-1} \log( \eproc_k' / \eproc_k) \geq \eps \right ) \leq &\underbrace{\Q \left ( \supkm k^{-1}\log(\eproc_n'/\eproc_n(\lambda_Q)) \geq \eps / 2 \right )}_{(\star)} +\\
    &\underbrace{\Q \left ( \supkm k^{-1} \log(\eproc_n(\lambda_Q)/\eproc_n) \geq \eps / 2 \right )}_{(\dagger)}.
  \end{align}
  Now, by property $(ii)$ and sublinearity $r_n = o(n)$, we have that $(\dagger)$ vanishes $\Qcal$-uniformly as $\mto$ so we focus on $(\star)$. By a union bound and the numeraire property of $\lambda_Q$, we have
  \begin{align}
    (\star) &\leq \sum_{k=m}^\infty \Q \left ( k^{-1} \log (\eproc_k' / \eproc_k(\lambda_Q)) \geq \eps / 2 \right ) \\
    \verbose{
            &\leq \sum_{k=m}^\infty \Q \left ( \eproc_k' / \eproc_k(\lambda_Q) \geq \exp \{ k\eps / 2 \} \right ) \\
    }
            &\leq \sum_{k=m}^\infty \EE_Q [ \eproc_k' / \eproc_k^\Qstar] \exp \left \{ -k\eps / 2 \right \}\\
    &\leq \sum_{k=m}^\infty  \exp \left \{ -k \eps / 2 \right \}.
  \end{align}
  By summability of $\exp \left \{ -k \eps / 2 \right \}$ and noting that the final line no longer depends on the distribution $\Qin$, we have that $(\star)$ vanishes $\Qcal$-uniformly and hence
 \begin{equation}
   \lim_\mto \sup_{\Qin} \Q \left ( \supkm k^{-1} \log (\eproc_k' / \eproc_k) \geq \eps \right ) = 0,
 \end{equation} 
 which completes the proof for property $(i)$ and hence the proof of
 \cref{theorem:uniform log optimality} altogether.
\end{proof}

\subsection{Proof of \cref{theorem:uniform stopping time}}
\label{proof:uniform stopping time}
\begin{proof}
    The proof is an immediate consequence of the nonasymptotic bound given in 
    \eqref{eq:proof-log-opt-time-unif-concentration}.
    Indeed, to see why the theorem holds, notice how the upper bound on 
    $\EE_Q[\tau_{\alpha}]$ from \eqref{eq:proof-expected-stopping-time-nonasymptotic-bound} 
    is finite due to the sublinearity of the portfolio regret bound $r_n$.
    The bound remains finite when we take the supremum over $Q \in \Qcalu$.
    Moreover, from the assumption that $\inf_{Q \in \Qcalu} \ell_Q^\star \geq
    \underline{\ell}(\Qcalu) > 0$ we have that $1/\ell_Q^\star$ is uniformly
    upper-bounded by $1/\underline{\ell}(\Qcalu)$. Dividing both sides of the
    upper bound by $\log(1/\alpha)$ and taking the limit supremum as $\alpha \to
    0^+$ completes the proof of the first claim.

    The second claim is a consequence of 
    \cref{proposition:expected stopping time oracle} combined with the
    assumption that $\sup_{Q \in \Qcalu} \ell_Q^\star \leq
    \widebar{\ell}(\Qcalu)$. Hence, we have completed the proof of
    \cref{theorem:uniform stopping time}.
\end{proof}

\newpage
\section{Remarks on $L_1$ log-optimality}\label{section:l1}

The main text focused on the study of \emph{almost sure} log-optimality, and on demonstrating that sublinear portfolio regret suffices to attain it. As alluded to in the introduction, another notion of asymptotic log-optimality has appeared in the literature, but with respect to convergence in $L_1$. In particular, \citet[\S 5.2]{wang2022backtesting} provide a definition of $L_1$ log-optimality and show that a follow-the-leader-type strategy referred to as ``Growth Rate for Empirical E-statistics'' (GREE)---see also ``Growth Rate Adaptive to the Particular Alternative'' (GRAPA) in \citep{waudby2020estimating}---enjoys $L_1$ optimality within the class of $e$-processes described by \eqref{eq:intro-general-test-sm-special-case}. This raises two natural questions: (a) Does sublinear portfolio regret also suffice for $L_1$ log-optimality, and (b) Does GREE satisfy sublinear portfolio regret? We provide a positive answer to (a). As for (b), it is well-known that follow-the-leader approaches can suffer linear regret without further structure (such as bounded domains or gradients, neither of which we have assumed here). Finally, given the discussion in \cref{section:uniformity}, it is natural to wonder whether $L_1$-log-optimality results hold in a \emph{distribution-uniform} sense. We also provide a positive answer to that question, relying on a distribution-uniform $L_1$ law of large numbers (\cref{proposition:l1-lln}) which seems to be new to the literature.

\subsection{$L_1$ log-optimality under sublinear portfolio regret}
\label{sec:l1-log-optimality}
Let us recall what it means for an $e$-process to be (distribution-uniformly) $L_1$-log-optimal within a class $\Wcal$. The pointwise half of the following definition is essentially found in \citet{wang2022backtesting}.
\begin{definition}[Asymptotic $L_1$ instance-log-optimality]\label{definition:l1-logopt}
  Let $\Wcal$ be a class of $\Pe$-processes and fix $W \in \Wcal$. Fix a class of alternatives $\Qcal$. We say that $W$ is pointwise $\Qcal$-instance $L_1$-log-optimal within $\Wcal$ if for every $W' \in \Wcal$ and every $\Qin$, it holds that
  \begin{equation}
    \lim_\nto \EE_\Q \pospart{\frac{1}{n} \log (W_n') - \frac{1}{n} \log(W_n)} = 0.
  \end{equation}
  Furthermore, if $\widebar \Qcal \subseteq \Qcal$ is a collection of distributions, we say that $W$ is $\widebar \Qcal$-uniformly $L_1$-log-optimal if
  \begin{equation}
    \lim_\nto \sup_{\Q \in \widebar \Qcal} \EE_\Q \pospart{\frac{1}{n} \log (W_n') - \frac{1}{n} \log(W_n)} = 0.
  \end{equation}
\end{definition}
With \cref{definition:l1-logopt} in mind, we now provide a nonasymptotic upper bound on the expectation appearing in \cref{definition:l1-logopt} which then establishes $\Qcal$-uniform $L_1$-log-optimality under sublinear portfolio regret. The conclusion of \cref{theorem:l1-logopt} can be viewed as a nonasymptotic and distribution-uniform analogue of \citet[Theorem 7.22]{ramdas2024hypothesis} but for $e$-processes with sublinear portfolio regret.
\begin{theorem}\label{theorem:l1-logopt}
  Suppose that $W_n$ forms a $\Pe$-process in $\Wcal$ with portfolio regret $R_n$. Then for any other $W' \in \Wcal$, we have
  \begin{equation}\label{eq:l1logopt upper bound}
    \EE_\Q \pospart{\frac{1}{n}\log(W_n') - \frac{1}{n}\log(W_n)} \leq \frac{1 + R_n}{n}.
  \end{equation}
  In particular, if $R_n$ is sublinear then $W_n$ is $\Qcal$-uniformly asymptotically $L_1$-log-optimal. Furthermore, for any $\widebar \Qcal \subseteq \Qcal$ for which $\log(\lambda_\Q^\trans \bE_1)$ is $\widebar \Qcal$-uniformly integrable, it holds that
  \begin{equation}
    \lim_\nto \sup_{\Q \in \widebar \Qcal} \EE_\Q \left \lvert \frac{1}{n}\log(W_n) - \EE_\Q[\log (\lambda_\Q^\trans \bE_1)] \right \rvert = 0.
  \end{equation}
\end{theorem}

\begin{proof} Writing out $\EE_\Q \pospart{\frac{1}{n}\log(W_n') - \frac{1}{n}\log(W_n)}$, we have
  \begin{align}
    \EE_\Q \pospart{\frac{1}{n} \sum_{i=1}^n \log \left ( \lambda_i^\trans \bE_i \right ) - \frac{1}{n} \sum_{i=1}^n  \log \left ( W_n \right )} &\leq \EE_\Q \pospart{\frac{1}{n} \log \left ( \lambda_i^\trans \bE_i \right ) - \max_{\lambda \in \simplex} \frac{1}{n} \sum_{i=1}^n \log \left ( \lambda^\trans \bE_i \right ) + \frac{R_n}{n} } \\
                                                                                                                                   &\leq \EE_\Q \pospart{\frac{1}{n} \sum_{i=1}^n  \log \left ( \lambda_i^\trans \bE_i \right ) - \frac{1}{n}\sum_{i=1}^n \log \left ( \lambda_\Q^\trans \bE_i \right ) + \frac{R_n}{n} } \\
               &\leq \EE_\Q \pospart{\frac{1}{n} \sum_{i=1}^n  \log \left ( \lambda_i^\trans \bE_i \right ) - \frac{1}{n} \sum_{i=1}^n \log \left ( \lambda_\Q^\trans \bE_i \right ) } + \frac{R_n}{n}.
  \end{align}
  Writing the expectation as an integral of tail probabilities, we have
  \begin{align}
    \EE_\Q \pospart{\frac{1}{n} \log \left ( \lambda_i^\trans \bE_i \right ) - \frac{1}{n} \sum_{i=1}^n \log \left ( \lambda_\Q^\trans \bE_i \right ) } &= \int_{0}^\infty\Q \left ( \pospart{\frac{1}{n} \log \left ( \lambda_i^\trans \bE_i \right ) - \frac{1}{n} \sum_{i=1}^n \log \left ( \lambda_\Q^\trans \bE_i \right ) } > t \right ) \dd t \\
                                                                                                                                                        &= \int_{0}^\infty\Q \left ( \prod_{i=1}^n \frac{\lambda_i^\trans \bE_i}{\lambda_\Q^\trans \bE_i} > \exp \{ n t \} \right ) \dd t\\
                                                                                                                                                        &\leq \int_0^\infty \EE_Q \left [ \prod_{i=1}^n  \frac{\lambda_i^\trans\bE_i}{\lambda_\Q^\trans \bE_i} \right ] \exp \{- nt \} \dd t \\
                                                                                                                                                        &\leq \int_0^\infty \exp \{-nt\} \dd t = \frac{1}{n}.
  \end{align}
  Therefore, it follows that 
  \begin{equation}
    \EE_\Q \pospart{\frac{1}{n} \sum_{i=1}^n \log \left ( \lambda_i^\trans \bE_i \right ) - \frac{1}{n} \sum_{i=1}^n  \log \left ( W_n \right )} \leq \frac{1 + R_n}{n},
  \end{equation}
  which completes the claim of \eqref{eq:l1logopt upper bound}. As for the second claim, it is easy to deduce from the previous arguments that
  \begin{equation}
    \EE_\Q \left \lvert \frac{1}{n} \log( W_n ) - \frac{1}{n} \sum_{i=1}^n  \log (\lambda_\Q^\trans \bE_i) \right \rvert \leq \frac{2}{n}.
  \end{equation}
  Applying the triangle inequality and the distribution-uniform $L_1$ law of large numbers in \cref{proposition:l1-lln}, we have
  \begin{equation}
    \lim_{n\to\infty} \sup_{\Q \in \widebar \Qcal}\EE_\Q \left \lvert \frac{1}{n} \log (W_n) - \EE_\Q [\log (\lambda_\Q^\trans \bE_1)] \right \rvert = 0.
  \end{equation}
  This completes the proof of \cref{theorem:l1-logopt}.
\end{proof}

\begin{proposition}[A distribution-uniform $L_1$ law of large numbers]\label{proposition:l1-lln}
  Let $\Pcal$ be a collection of probability measures and let $X_1, \dots, X_n$ be \iid{} random variables with $\P$-mean zero for each $\Pin$. Suppose that $X_1$ is $\Pcal$-uniformly integrable, meaning that
  \begin{equation}
    \lim_{K \to \infty} \sup_\Pin \EE_\P \left [ \left \lvert X_1 \right \rvert \1 \{ |X_1| \geq K \} \right ] = 0.
  \end{equation}
  Then the $L_1$ law of large numbers holds uniformly in $\Pcal$:
  \begin{equation}
    \lim_{n\to\infty}\sup_\Pin \EE_\P \left \lvert \frac{1}{n} \sum_{i=1}^n X_i \right \rvert = 0.
  \end{equation}
\end{proposition}
\begin{proof}
  Put $S_n = \sum_{i=1}^n X_i$ and fix $\eps > 0$. Let $A \equiv A_n(\eps)$ be the event given by
  \begin{equation}
    A_n(\eps) := \left \{ |S_n / n| \geq \eps \right \}.
  \end{equation}
  Analyzing the expected absolute value of $S_n / n$, we have that for any $\Pin$,
  \begin{align}
    \EE_\P \left \lvert S_n / n \right \rvert &= \EE_\P \left [ |S_n / n| \1_{A_n(\eps)} \right ] + \EE_\P \left [ |S_n / n| \1_{A_n(\eps)} \right ] \\
                                                            &\leq \EE_\P \left[ \left \lvert X_1 \right \rvert \1_{A_n(\eps)} \right ] + \eps.
  \end{align}
  Letting $r_n(\eps) = \sqrt{\P (A_n(\eps))}$, we have
  \begin{align}
    \EE_\P \left \lvert S_n / n \right \rvert &\leq \EE_\P \left [ |X_1| \1 \left \{ |X_1| \geq \frac{1}{r_n(\eps)} \right \} \right ] +  \frac{\P(A_n(\eps))}{r_n(\eps)} + \eps.
  \end{align}
  Now, by \citet[Theorem 1]{waudby2024distribution}, we have that $r_n(\eps)$ vanishes to 0 uniformly in $\Pin$ as $n \to \infty$. Therefore, taking suprema over $\Pin$ on both sides of the inequality and limits as $n \to \infty$, we have that
  \begin{equation}
    \lim_{n\to\infty}\sup_\Pin \EE_\P \left \lvert S_n / n \right \rvert \leq \eps.
  \end{equation}
  Since $\eps > 0$ was arbitrary, we have the desired result.
\end{proof}

\newpage

\section{Remarks on admissibility of sequentially merged e-processes}
\label{section:admissibility}
In this section, we show that there is a sense in which the $e$-processes taking the form \eqref{eq:intro-general-test-sm} are the only admissible ones. Before stating the result, we first establish some definitions.
\begin{definition}[E-merging function]\label{definition:e-merging}
    We say that an increasing Borel function $g : [0, \infty]^{d+1} \to [0, \infty]$ is an \emph{$e$-merging function} if for any probability space $(\Omega, \Fcal, P)$ and $e$-variables $E^{(0)}m \dots, E^{(d)}$, it holds that $g\left(E^{(0)}, \dots E^{(d)}\right)$ is an $e$-variable. That is, it holds that $g\left(E^{(0)}, \dots E^{(d)}\right) \geq 0$ $P$-almost surely and $\EE_{P}\left[g\left(E^{(0)}, \dots E^{(d)}\right)\right] \leq 1$.
\end{definition}
In words, an \emph{e-merging} function is one that takes a tuple of e-variables and produces an e-variable. See also \citet{vovk2021values,vovk2024merging} and \citet{wang2025only}.
The following is a natural extension of \cref{definition:e-merging} to $e$-processes that are formed by taking sequential products of $e$-variables.

\begin{definition}[Sequentially merged $e$-process]
  Let $\infseqn{\bE_n}$ be \iid{} $(d+1)$-tuples of arbitrarily dependent $\Pcal$-e-variables. Define $\Fcal_n = \sigma(\bE_1,\dots, \bE_n)$ for each $n \in \NN$. We say that a $\Pcal$-$e$-process $\infseqn{W_n}$ is \emph{sequentially merged} if $W_n / W_{n-1} = g_n(E_n^\brackzero,\dots, E_n^\brackd)$ for some coordinate-wise increasing $\Fcal_{n-1}$-measurable $e$-merging function $g_n$.
\end{definition}

Following \citet{vovk2024merging}, \citet{wang2025only}, and \citet{clerico2025optimality}, an $e$-process $(W_n)_{n\in\NN}$ is said to be \emph{admissible} if it cannot be pointwise improved. That is, if there exists some other $e$-process $(W_n')_{n \in \NN}$ so that $W_n \leq W_n'$ for all $n \in \NN$, it must hold that $W_n = W_n'$ for all $n \in \NN$. With these definitions and concepts in mind, we now demonstrate that sequentially merged dependent tuples of $e$-variables must take the form in \eqref{eq:intro-general-test-sm}.

\begin{proposition}
\label{prop:admissibility}
  Any admissible sequentially merged $e$-process takes the form
  \begin{equation}
    \prod_{i=1}^n \lambda_i^\trans \bE_i
  \end{equation}
  for some predictable $\simplex$-valued sequence $\infseqn{\lambda_n}$
\end{proposition}
\begin{proof}
  We proceed by induction. The result for $n = 0$ and $n = 1$ follows from \citet{wang2025only}. For arbitrary $n \in \NN \setminus \{1\}$,
  we have that $W_{n-1} = \prod_{i=1}^{n-1} \lambda_{i}^\trans \bE_i$. Now, since $W_n$ is sequentially merged, it holds that $W_n / W_{n-1} = g_n(E_n^\brackzero, \dots, E_n^\brackd)$ for an $\Fcal_{n-1}$-measurable e-merging function $g_n$. Another application of \citet[Theorem 1]{wang2025only} yields $g_n(e_0, \dots, e_d) = \lambda^\brackzero e_0 + \dots + \lambda^\brackd e_d$ for some $\simplex$-valued and $\Fcal_{n-1}$-measurable $\lambda$. This completes the proof.
\end{proof}

    \newpage
\section{Remarks on two-sample and independence testing}\label{section:two-sample}
Suppose we have access to a sequence of \iid{} random variables $\infseqn{X_n, Y_n} \sim P \equiv P_X \times P_Y$ and we are tasked with testing whether $X_1$ and $Y_1$ have the same distribution. More formally, we posit the null $\Pcal := \{ P : P_X = P_Y \}$ versus the alternative $\Qcal := \{ P : P_X \neq P_Y \}$. The approach proposed by \citet{shekhar2023near} is based on considering a distance measure $d_\Gcal$ between $P_X$ and $P_Y$ that admits a variational representation for a class of functions $\Gcal$:
\begin{equation}\label{eq:ipm}
  d_\Gcal(P_X, P_Y) = \sup_{g \in \Gcal} \left \lvert \EE_{P}[ g(X_1)] - \EE_P [g(Y_1)] \right \rvert,
\end{equation}
and so that $g$ takes values in $[-1, 1]$ for any $g \in \Gcal$. The exact technical requirements for $\Gcal$ are  outside the scope of this paper. Importantly, however, \citet{shekhar2023nonparametric} show the value of focusing on a special class of functions $\Gcal$ with the property that
\begin{equation}
  d_\Gcal(P_X, P_Y) = 0  \quad\textit{if and only if} \quad P_X = P_Y,
\end{equation}
which allows them to reduce the problem of testing $\Pcal$ versus $\Qcal$ to that of whether a supremum of means $d_\Gcal(P_X, P_Y)$ is zero or non-zero. Of course, it is not yet immediately obvious if or how one can compute $d_\Gcal$ or an estimate thereof. The authors consider the so-called \emph{oracle witness function} $g^\star_{P} \in \Gcal$; i.e., the function that attains the supremum in \eqref{eq:ipm} (while it cannot be assumed in general that the supremum lies in $\Gcal$, the authors mainly consider cases where it does, and they analyze an approximate supremum in cases where it does not). Furthermore, they focus on classes $\Gcal$ that are symmetric so that $d_\Gcal(P_X, P_Y)$ can be written without the absolute value:
\begin{equation}\label{eq:ipm-final-form}
  d_\Gcal(P_X, P_Y) = \sup_{g \in \Gcal} \left ( \EE_P [g(X_1)] - \EE_P[g(Y_1)] \right )  = \EE_P[g^\star_P(X_1)] - \EE_P[g^\star_P(Y_1)].
\end{equation}
When written in the form of \eqref{eq:ipm-final-form}, we see that \citet{shekhar2023nonparametric} have effectively reduced the problem of two-sample testing to that of difference-in-means testing for bounded random tuples as in \cref{section:difference-in-means}, with the tuples $\infseqn{g_P^\star(X_1), g_P^\star(Y_1)}$ for some $g_P^\star$. Of course, we are glossing over numerous subtleties pertaining to the fact that $g_P^\star$ depends on the distribution $P \equiv P_X \times P_Y$, the fact that it is not known \emph{a priori}, and that it must be learned over time with data. These details are beyond the scope of this paper. Importantly for our purposes, though, the authors describe the oracle test supermartingale as one which uses the witness function $g_Q^\star$ under an alternative $Q$ along with the optimal constant betting strategy:
\begin{equation}
  \eproc_n^\TwoSTOracle := \prod_{i=1}^n \left ( 1 + \lambda^\star_Q \cdot (g_Q^\star(X_n) - g_Q^\star(Y_n)) \right ),
\end{equation}
where $\lambda_Q^\star := \argmax_{\lambda \in [0, 1]} \EE_Q [ \log(1 + \lambda \cdot (g_Q^\star(X_1) - g_Q^\star(Y_1) ) ) ]$. However, their tests do not attain the same rate of growth as $\eproc_n^\TwoSTOracle$ since their use of the ONS betting strategy effectively restricts the search space of $\lambda$ to $[0, 1/2]$ rather than $[0, 1]$. We conjecture that our results will provide the requisite foundations to develop two-sample tests that attain the same rate of growth as $\eproc^\TwoSTOracle$ (and bounds on expected rejection times) even when the witness function is unknown.

As a brief remark on sequential marginal independence testing, the high-level goal is akin to two-sample testing but the null is given by $\Pcal := \{ P : P_X \indep P_Y \}$ versus the alternative $\Qcal := \{ P : P_X \cancel{\indep} P_Y \}$. When analyzed in batches of two, independence testing can essentially be reduced to two-sample testing since $P_X \indep P_Y$ if and only if the following four tuples have the same distribution: $(X_1, Y_1)$, $(X_1, Y_2)$, $(X_2, Y_1)$, and $(X_2, Y_2)$---though there are additional subtleties that arise and these are again outside the scope of this paper. We direct the interested reader to \citet{podkopaev2023sequentialKernelized} for a detailed treatment of this problem. Again, importantly for our purposes, \citet{podkopaev2023sequentialKernelized} describe an oracle test that implements the maximizer over $\lambda \in [0, 1]$ of certain oracle log-wealth increments (see \citep[Remark 1]{podkopaev2023sequentialKernelized}) but similar to \citet{shekhar2023nonparametric} their use of the ONS betting strategy effectively restricts the search space to $[0, 1/2]$ so that the oracle growth rate cannot be achieved in general. 

\newpage
\section{Abstract portfolio regret and numeraire portfolios}\label{section:generalizations}
The main theorems of \cref{section:asymptotic-equivalence,section:stopping
time,section:uniformity} centrally relied on $e$-processes having sublinear portfolio regret as well as on certain properties of the log-optimal strategy $\lambda_Q^\star$. In the present section, we abstract away from the concrete
example of the test supermartingale given in \eqref{eq:intro-general-test-sm},
and show that the optimality results for that test supermartingale hold in a
more abstract setting. 

Concretely, throughout we fix a measurable set $\Theta$ and we denote specific elements of $\Theta$ by $\theta$. Here, $\theta\in \Theta$ will play the role of generalizing a particular constant rebalanced portfolio $\lambda \in \dsimp$, and $\Theta$-valued predictable sequences $\infseqn{\theta_n}$ will play the role of generalizing $\dsimp$-valued portfolios $\infseqn{\lambda_n}$. Letting the filtration $\Fcal$ be generated by some \iid{} sequence of random objects $\infseqn{X_n}$, let $E_n : \Theta \to [0, \infty]$ be a function so that $E_n(\theta)$ is $\sigma(X_n)$-measurable for each $\thetain$ and $n \in \NN$. The mapping $E_n(\theta)$ should be thought of as a generalization of the random variable $\lambda_n^\trans \bE_n$ that we have studied thus far, where the latter can be viewed as arising from letting $X_n = \bE_n$ and $\theta_n = \lambda_n \in \Theta = \dsimp$. Define the composite null hypothesis $\widebar \Pcal$ and alternative $\widebar \Qcal$ by
\begin{align}
  \widebar \Pcal := \{ P : E_1(\theta) \text{ is a $P$-$e$-value for \emph{every} $\thetain$} \}
\end{align}
\begin{center}and\end{center}
\begin{equation}
  \widebar \Qcal := \{ P : E_1(\theta) \text{ is not a $P$-$e$-value for \emph{some} $\thetain$}\}.
\end{equation}
Furthermore, define the collection $\eprocclass(\Theta)$ of $\widebar \Pcal$-$e$-processes by
\begin{equation}\label{eq:generalized class}
  \eprocclass(\Theta) := \left \{ \eproc : \eproc_n \leq \prod_{i=1}^n E_i(\theta_i) \text{ and $\infseqn{\theta_n}$ is a $\Theta$-valued $\Fcal$-predictable sequence} \right \}.
\end{equation}
With all of this setup in mind, we are ready to define both abstract portfolio regret and the $(\Theta, Q)$-numeraire portfolio.
\begin{definition}[Abstract portfolio regret and $(\Theta, Q)$-numeraire portfolios]
  \label{definition:abstract regret and numeraire}
  Fix an alternative distribution $\Qinbar$. We will say that $\eproc \in \eprocclass(\Theta)$ satisfies a \emph{abstract portfolio regret} bound of $r_n(Q)$ if
\begin{equation}\label{eq:abstract portfolio regret}
  \sup_\thetain \sum_{i=1}^n \log (E_i(\theta)) - \log (\eproc_n) \leq r_n(Q),
\end{equation}
with $Q$-probability one and we say that its regret is $\Qcalbar$-uniformly sublinear if there exists a deterministic sublinear sequence $\infseqn{r_n}$ so that $\supQbar r_n(Q) \leq r_n$ for each $n \in \NN$. Furthermore, we say that $\theta_Q^\star \in \Theta$ is the \emph{$(\Theta, Q)$-numeraire portfolio} if for all $\Theta$-valued predictable sequences $\infseqn{\theta_n}$, the process
\begin{equation}\label{eq:general-assumption-suboptimality-ratio}
  \prod_{i=1}^n \left (E_i(\theta_i) / E_i(\theta_Q^\star)\right )
\end{equation}
forms a nonnegative $Q$-supermartingale with $\EE_Q[E_1(\theta) / E_1(\theta_Q^\star)] \leq 1$.
\end{definition}
Consider the following two assumptions.

\begin{assumption}[$\Qcalbar$-uniformly sublinear regret bound]
\label{assumption:generalized-portfolio-regret}
    We assume there exist $\widebar \Pcal$-$e$-processes that satisfy the abstract portfolio regret bound from \eqref{eq:abstract portfolio regret} for some deterministic $r_n(Q)$ with $Q$-probability one for each $Q \in \Qcalbar$.
\end{assumption}

\begin{assumption}[Abstract numeraire portfolios]
\label{assumption:abstract-numeraire-portfolios}
    We assume there exists a $(\Theta, Q)$-numeraire $\theta_Q^\star \in \Theta$ for each $Q \in \Qcalbar$, meaning that $\theta_Q^\star$ satisfies \eqref{eq:general-assumption-suboptimality-ratio}.
\end{assumption}

The choice to refer to $\theta_Q^\star$ as the $(\Theta, Q)$-numeraire portfolio is directly inspired by both \citet{long1990numeraire} and \citet{larsson2024numeraire}. However, we wish to distinguish the above from the definition of a numeraire $e$-value as in \citet{larsson2024numeraire} since the authors study a stronger property than \eqref{eq:general-assumption-suboptimality-ratio} when $n=1$ where there is a unique $\widebar \Pcal$-$e$-value $E^\star_1$ for which $\EE_Q[E_1' / E_1^\star] \leq 1$ for \emph{any} other $\widebar \Pcal$-$e$-value $E_1'$, not just those indexed by $\thetain$. It is for this reason that we qualify the definition of $\theta_Q^\star$ with the set $\Theta$ rather than simply calling it ``the $Q$-numeraire portfolio.'' Moreover, we qualify the definition of $\theta_Q^\star$ with the particular distribution $Q$ in order to emphasize that we will be making use of a \emph{family} of portfolios $(\theta_Q^\star)_{\Qinbar}$ indexed by the composite alternative $\Qcalbar$. The numeraire $e$-value defined in \citet{larsson2024numeraire} is with respect to a point alternative. 

\subsection{Concentration and log-optimality for an abstract class of $e$-processes}\label{section:general-concentration-log-optimality}
A key result used to prove \cref{theorem:adaptive optimality} is a time-uniform exponential concentration inequality for $a_n^{-1} [\log(\eproc_n) - \log(\eproc_n(\lambda_Q^\star))]$ for any $\eproc \in \eprocclass$. In fact, the aforementioned concentration inequality holds for arbitrary test supermartingales satisfying certain properties introduced in \cref{definition:abstract regret and numeraire}. %
The details follow.

\begin{theorem}\label{lemma:adaptive optimality}
    Suppose that Assumptions~\ref{assumption:generalized-portfolio-regret} and~\ref{assumption:abstract-numeraire-portfolios} hold.
    Let $\eproc \in \eprocclass(\Theta)$ be any $\widebar \Pcal$-$e$-process satisfying the abstract portfolio regret bound \eqref{eq:abstract portfolio regret} for some deterministic $r_n(Q)$ with $Q$-probability one for each $\Qinbar$. Then defining the test $\widebar \Pcal$-supermartingale $\eproc_n(\theta_Q^\star) := \prod_{i=1}^n E_i(\theta_Q^\star)$, we have for any deterministic sequence $\infseqn{a_n}$ and $m \in \NN$,
    \begin{align}
      \forall \eps > 0,\quad &\Q \left ( \supkm | a_k^{-1} \log (\eproc_k / \eproc_k(\theta_Q^\star)) | \geq \eps \right ) \\
      &\leq \sum_{k=m}^\infty \exp \{ -a_k \eps \} + \1 \left \{ \supkm a_k^{-1} r_k(Q) \geq \eps \right \}.\label{eq:log-optimality inequality}
      \end{align}
      Consequently, if $r_n(Q)$ is $\widebar \Qcal$-uniformly sublinear, i.e.~$\sup_\Qinbar r_n(Q) = o(1)$, then Properties $(i)$, $(ii)$, and $(iii)$ hold from \cref{theorem:adaptive optimality} but with respect to the larger class $\eprocclass(\Theta)$ and with $\ell_Q^\star := \EE_Q[\log E_1(\theta_Q^\star)]$.

\end{theorem}

The first and second terms in the right-hand side of \eqref{eq:log-optimality inequality} can be thought of as controlling $a_k^{-1} \log (\eproc_k/ \eproc_k(\theta_Q^\star))$ with high probability from above and from below by appealing to the numeraire property \eqref{eq:general-assumption-suboptimality-ratio} and sublinear regret \eqref{eq:abstract portfolio regret}, respectively. The key ingredients of the proof can be further distilled in words as follows.
First, for any $\Qinbar$, the $(\Theta,Q)$-numeraire property in \eqref{eq:general-assumption-suboptimality-ratio} is used to derive an exponential concentration inequality implying that $\log(\eproc_n(\theta_Q^\star))$ will eventually forever exceed $\log(\eproc_n')$ for any other $\eproc' \in \eprocclass(\Theta)$ in finite time. This fact can be seen as a quantitative strengthening of \citet[Theorem 15.3.1]{cover1999elements} (and for more general processes). Second, by the assumption that $\eproc_n$ enjoys an abstract portfolio regret bound, we have that $\log (\eproc_n)$ is close to $\sup_\thetain \log (\eproc_n(\theta))$, which is itself always larger than $\log (\eproc_n(\theta_Q^\star))$ by definition. Combined with the first fact that $\log (\eproc_n(\theta_Q^\star))$ will eventually exceed $\log (\eproc_n)$, it must be the case that $\log(\eproc_n)$, $\log(\eproc_n(\theta_Q^\star))$, and $\sup_\thetain \log(\eproc_n(\theta))$ are all sandwiched within $r_n(Q)$ of each other in finite time with $Q$-probability one. The formal details are in \cref{proof:abstract adaptive optimality}.

We now study expected rejection times and show how they have matching lower and upper bounds within
$\eprocclass(\Theta)$ under sublinear abstract portfolio regret.

\subsection{Expected rejection times via regret and numeraire portfolios}\label{section:general-bounds-stopping-time}
  Similar to \cref{theorem:stopping time} in \cref{section:stopping time}, we will require that the $s^\tth$ moment of $\log E_1(\theta_Q^\star)$ \emph{under the $(\Theta, Q)$-numeraire} is bounded for some $s > 2$:
  \begin{equation}
    \widebar \rho_s(\theta_Q^\star) := \EE_Q \left \lvert \log E_1(\theta_Q^\star) - \ell_Q^\star \right \rvert^s < \infty.
  \end{equation}
  Similar to \cref{lemma:adaptive optimality}, once combined with the finite-moment condition above, the essential properties for deriving bounds on expected rejection times in the more general setting considered in the present section are (1) sublinear abstract portfolio regret and (2) the existence of a family of numeraire portfolios.
  The following result not only generalizes \cref{theorem:stopping time} to $e$-processes in $\eprocclass(\Theta)$ but the upper bound on the expected rejection time is written in a nonasymptotic fashion with explicit constants. While the expression is somewhat involved, it contains more information than when limits are taken with respect to $\alphato$.
\begin{theorem}\label{lemma:general nonasymptotic bound}
    Consider the collection of $\widebar \Pcal$-$e$-processes $\eprocclass(\Theta)$ given in \eqref{eq:generalized class} and let $\eproc \in \eprocclass(\Theta)$. Suppose that for every $\Qin$, $\eproc$ satisfies \cref{assumption:generalized-portfolio-regret}, \cref{assumption:abstract-numeraire-portfolios} holds, and that the $(\Theta, Q)$-numeraire portfolio $\theta^\star_Q$ has log increments with bounded $s^\tth$ moments: $\widebar \rho_s(\theta_Q^\star) < \infty$. Then the $Q$-expectation of $\widebar \tau_\alpha := \inf \{ n \in \NN : \eproc_n \geq 1/\alpha \}$ can be upper bounded for any $\delta \in (0, 1)$ as
  \begin{align}
    \EE_Q[\widebar \tau_\alpha] \leq &4 + \frac{(1+\delta)\log(1/\alpha)}{\ell_Q^\star} + \frac{2(1+\delta)\alpha^{\delta/2}}{\delta \ell_Q^\star} \\
    &+ \frac{(1+\delta)^{1+s/2}}{\delta^s (\ell_Q^\star)^{1+s/2} \log^{s/2 -
    1}(1/\alpha)} + \sum_{k=m}^\infty \1 \left \{ r_k(Q)/k \geq \frac{\delta
    \ell_Q^\star}{2(1+\delta)} \right \},\label{eq:stopping time upper bound}
  \end{align}
  where $m = \left \lceil (1+\delta)\log(1/\alpha) / \ell_Q^\star  \right \rceil$. 
  Furthermore, the following lower bound holds for any predictable strategy $\infseqn{\widetilde \theta_n}$ and corresponding rejection time $\widetilde \tau_\alpha := \inf \{ n \in \NN : \prod_{i=1}^n E_i(\widetilde \theta_i) \geq 1/\alpha  \}$:
  \begin{equation}\label{eq:general stopping time lower bound}
    \frac{\EE_Q[\widetilde \tau_\alpha]}{\log(1/\alpha)} \geq \frac{1}{\ell_Q^\star}.
  \end{equation}
\end{theorem}
The proof of \cref{lemma:general nonasymptotic bound} is provided in \cref{proof:general stopping time}.
Notice that the infinite series in \eqref{eq:stopping time upper bound} is always finite if $\supQbar r_k(Q)$ is sublinear since only finitely many of the summands can be nonzero. Moreover, as $\alpha$ becomes small, the integer $m \asymp \log(1/\alpha)$ becomes large, in which case the aforementioned series is not just finite but in fact zero. Finally, notice that the second term in the upper bound on $\EE_Q[\widebar \tau_\alpha]$ is the only one that diverges as $\alpha \to 0^+$, and hence in the small-$\alpha$ regime, the dominant term is $(1+ \delta) \log(1/\alpha ) / \ell_Q^\star$. This term matches the lower bound in \eqref{eq:general stopping time lower bound} up to a factor of $(1+\delta)$ which can be made arbitrarily close to 1 in the limit supremum, hence yielding the equality
\begin{equation}
    \lim_\alphato \frac{\EE_Q[\widebar \tau_\alpha]}{\log(1/\alpha)} = \frac{1}{\ell_Q^\star}.
\end{equation}
Let us now investigate some example testing problems that fall under the more abstract setting of the present section but that could not have been described by the class $\eprocclass$.

\begin{example}[Sub-Gaussian random variables]\label{example:sub-gaussian}
  Fix $\sigma > 0$. Let $\infseqn{X_n}$ be random variables and let $\Pcal$ consist of all measures for which $\infseqn{X_n}$ are \iid{} and $\sigma$-sub-Gaussian with mean $\mu_0$. Consider the class of $\Pe$-processes given by
  \begin{equation}
    W_n := \exp \left \{ \sum_{i=1}^n\lambda_i (X_i - \mu_0) - \sum_{i=1}^n \frac{\sigma^2 \lambda_i^2}{2} \right \}.
  \end{equation}
  By \citet{larsson2024numeraire}, there exists $\lambda_\Q$ for which
  \begin{equation}
    W_n(\lambda_Q) \coloneqq \exp \left \{ - \lambda_\Q \sum_{i=1}^n (X_i - \mu_0) + \frac{n\sigma^2 \lambda_\Q^2}{2}  \right \}
  \end{equation}
  forms a nonnegative $\Q$-supermartingale with initial value 1. Moreover, as noted in \citet{agrawal2025eventually}, it holds that Robbins' Gaussian mixture $\Pcal$-supermartingale
  \begin{equation}
      W_n^\mathrm{mix} := \exp \left \{ \frac{\sigma^2 \left(\sum_{i=1}^n (X_i - \mu_0)\right)^2}{2(n \sigma^4 + 1)} \right \}\left(n \sigma^4 + 1\right)^{-1/2}
  \end{equation}
  enjoys logarithmic regret with $\Q$-probability one. Consequently, for $\tau_\alpha := \1 \{ W_n^\mathrm{mix} \geq 1/\alpha \}$,
  \begin{equation}
    \Q \left ( \frac{1}{n} \log (W_n^\mathrm{mix}) \to \EE_\Q [(\mu_\Q - \mu_0)^2] \right ) = 1\quad\text{and}\quad \lim_\alphato\frac{\EE_\Q [\tau_\alpha]}{\log(1/\alpha)} = \frac{1}{(\mu_\Q - \mu_0)^2}.
  \end{equation}
\end{example}

\begin{example}[The Universal Inference $e$-process \citep{wasserman2020universal}]\label{example:universalinference}
  Let $\Pcal = \{ \P_\theta \}_{\theta \in \Theta}$ be a collection of distributions for a finite dimensional parameter space $\Theta$, and suppose that this collection has a common reference measure. Let $\Theta_0, \Theta_1 \subseteq \Theta$ be the null and alternative parameter spaces, respectively, so that $\Theta_0 \cap \Theta_1 = \emptyset$. The Universal Inference $e$-process of \citet{wasserman2020universal} is given by
  \begin{equation}
    L_n^\UI := \prod_{i=1}^n \frac{\dd \P_{\widehat \theta_{i-1}}(X_i)}{\max_{\theta \in \Theta_0}\dd\P_{\theta}(X_i)},
  \end{equation}
  where $\widehat \theta_{i-1}$ is some $\sigma(X_1, \dots, X_{i-1})$-measurable random variable taking values in $\Theta_1$. 
  Then for $\theta^\star \in \Theta$, it holds that
  \begin{equation}
    L_n^\star := \prod_{i=1}^n \frac{\dd \P_{\theta^\star}(X_i)}{\dd \P_{\theta_0}(X_i)}
  \end{equation}
  is a $(\Theta, \P_{\theta^\star})$-numeraire portfolio. Indeed, this is easy to see since for any $\Theta_1$-valued adapted sequence $\infseqn{\widehat \theta_n}$, we have
  \begin{equation}
    \frac{L_n^\UI}{L_n^\star}= \prod_{i=1}^n \frac{\dd \P_{\widehat \theta_{i-1}}(X_i)}{\dd \P_{\theta^\star}(X_i)}
  \end{equation}
  forms a $\P_{\theta^\star}$-supermartingale.

  First, consider the case where $\Theta_0 = \{\theta_0\}$ is a singleton.
  Suppose that $\infseqn{W_n}$ is a $\P_{\theta_0}$-$e$-process with an almost sure abstract portfolio regret given by $R_n(\theta)$ under $\P_\theta$. If $\sup_{\theta \in \Theta} R_n$ is almost surely sublinear, then by \cref{lemma:general nonasymptotic bound}, we have for any alternative $\theta^\star \in \Theta$,
  \begin{equation}
    \P_{\theta^\star} \left ( \frac{1}{n} \log W_n \to \KL(\P_{\theta^\star} \| \P_{\theta_0}) \right ) = 1 \quad\text{and}\quad\lim_\alphato \frac{\EE_{\P_{\theta^\star}}\left [ \tau_\alpha \right ]}{\log(1/\alpha)} = \frac{1}{\KL(\P_{\theta^\star} \| \P_{\theta_0})}.
  \end{equation}
  On the other hand, if $\Theta_0$ is not a singleton and it additionally holds that $\max_{\theta \in \Theta_0} \frac{1}{n} \sum_{i=1}^n  \log(\dd \P_\theta(X_i)) \to \max_{\theta \in \Theta_0}\EE_{\P_{\theta^\star}}[\log(\dd \P_\theta(X_1))]$ with $\P_{\theta^\star}$-probability one, then
  \begin{equation}
    \P_{\theta^\star} \left ( \frac{1}{n} \log W_n \to \min_{\theta \in \Theta_0}\KL(\P_{\theta^\star} \| \P_{\theta}) \right ) = 1 \quad\text{and}\quad\lim_\alphato \frac{\EE_{\P_{\theta^\star}}\left [ \tau_\alpha \right ]}{\log(1/\alpha)} = \frac{1}{\min_{\theta \in \Theta_0}\KL(\P_{\theta^\star} \| \P_{\theta})}.
  \end{equation}

  It is not obvious whether it is possible to derive $e$-processes with sublinear portfolio regret for given parametric families, but there are examples where such regret bounds are known. For example, logarithmic regret bounds for Bernoulli random variables were derived by \citet{krichevsky1981performance}.

  When testing the mean of a Gaussian random variable with variance $\sigma^2$, it can be shown that the likelihood ratio between some null $\mu_0$ and an alternative $\mu_1$ is given by and can be equivalently written as
 \begin{equation}
   \prod_{i=1}^n \frac{\exp \left \{ - \frac{(X_i-\mu_1)^2}{2\sigma^2} \right \}}{\exp \left \{ - \frac{(X_i-\mu_0)^2}{2\sigma^2} \right \}} = \exp \left \{ \lambda \sum_{i=1}^n (X_i - \mu_0) - \frac{n \sigma^2 \lambda^2}{2}\right \},
 \end{equation} 
 where $\lambda = (\mu_1 - \mu_0) / \sigma^2$. Clearly there is a one-to-one correspondence between $\lambda$ and $\mu_1$, and hence the logarithmic regret bound discussed in \cref{example:sub-gaussian} applies here as well.
\end{example}

\subsection{Proof of \cref{lemma:adaptive optimality}}\label{proof:abstract adaptive optimality}
\begin{proof}[Proof of \eqref{eq:log-optimality inequality}]
  We aim to provide an upper bound on the following probability for any $\Qinbar$ and any $\eps > 0$
  \begin{equation}
    \Q \left ( \supkm \left ( a_k^{-1} \left \lvert \log W_k - \log W_k(\theta_Q^\star) \right \rvert \right ) \geq \eps \right ).\label{eq:quantity-to-bound}
  \end{equation}
  Note that we can upper bound \eqref{eq:quantity-to-bound} by
  \begin{align}
    \eqref{eq:quantity-to-bound} \leq \underbrace{\Q \left ( \supkm \left ( a_k^{-1}\log(W_k / W_k(\theta_Q^\star)) \right ) \geq \eps \right )}_{(+)} + \underbrace{\Q \left ( \supkm \left ( a_k^{-1}\log(W_k(\theta_Q^\star) / W_k) \right ) \geq \eps \right )}_{(-)},
  \end{align}
  where the terms $(+)$ and $(-)$ are probabilities in terms of events that $a_k^{-1}\log(W_k / W_k(\theta_Q^\star))$ are $\eps$-far apart for any $k \geq m$ from above or below, respectively. Looking first to $(+)$, we have
  \begin{align}
    (+) &= \Q \left ( \exists k \geq m : a_k^{-1} \log \left ( W_k / W_k(\theta_Q^\star) \right ) \geq \eps \right )\\
        &= \Q\left ( \exists k \geq m : W_k / W_k(\theta_Q^\star) \geq \exp \left \{ a_k \eps \right \} \right ) \\
        &\leq \sum_{k=m}^\infty \frac{\EE_Q (W_k / W_k(\theta_Q^\star))}{\exp \{a_k \eps \}} \leq \sum_{k=m}^\infty \exp \{-a_k \eps \},
  \end{align}
  where the first inequality follows from a union bound and Markov's inequality, and the second follows from the numeraire property in \eqref{eq:general-assumption-suboptimality-ratio}.
  Turning now to $(-)$, we have that
  \begin{align}
    (-) &= \Q \left ( \supkm \left( a_k^{-1}\log(W_k(\theta_Q^\star) / W_k) \right ) \geq \eps \right ) \\
        &\leq \Q \left ( \supkm \left (  a_k^{-1} \left [\log W_k^\Max - \log W_k \right ] \right ) \geq \eps \right ) \\
    &\leq \1 \left \{ \supkm a_k^{-1} r_k(Q) \geq \eps \right \},
  \end{align}
  where $\eproc_k^\Max := \sup_{\theta \in \Theta} \log \eproc_k (\theta) \geq \log \eproc_k(\theta_Q^\star)$.
  Putting the two inequalities for $(+)$ and $(-)$ together, we have that
  \begin{equation}
    \eqref{eq:quantity-to-bound} \leq \sum_{k=m}^\infty \exp \left \{ -a_k \eps
    \right \} + \1 \left \{ \supkm a_k^{-1}r_k(Q)\geq \eps \right \},
  \end{equation}
  which completes the proof.
\end{proof}

With access to the inequality \eqref{eq:log-optimality inequality}, Properties $(i)$--$(iii)$ follow with just a few more arguments.
We proceed by first proving property $(ii)$, then $(iii)$, and then $(i)$. Throughout, let $\eproc_n^\Qstar := \eproc_n(\theta_Q^\star)$ and $r_n := \supQbar r_n(Q)$ to ease notation.
\begin{proof}[Proof of property $(ii)$.]
  By the inequality in \eqref{eq:log-optimality inequality}, it follows that for any $\Qinbar$ and any $\eps > 0$,
  \begin{equation}
    \Q \left ( \supkm \left \lvert a_k^{-1} \log (\eproc_k / \eproc_k^\Qstar) \right \rvert \geq \eps \right ) \leq \sum_{k=m}^\infty \exp \left \{ -a_k \eps \right \} + \1 \left \{ \supkm a_k^{-1} r_k \geq \eps \right \}.
  \end{equation}
  Now by the assumptions that $\sum_{k=1}^\infty \exp \left \{ -a_k\eps \right \} < \infty$ and $a_n^{-1} r_n \to 0$, we have that the right-hand side of the above inequality vanishes as $\mto$. Since for any sequence of random variables $\infseqn{Z_n}$, it is the case that $Z_n \to 0$ as $\nto$ with $Q$-probability one if and only if $\supkm |Z_k| \to 0$ in $Q$-probability as $\mto$, we have that
  \begin{equation}
    a_n^{-1} \log \left ( \eproc_n / \eproc_n^\Qstar \right ) \to 0 \quad\text{with $Q$-probability one,}
  \end{equation}
  which completes the proof of property $(ii)$.
\end{proof}
  \begin{proof}[Proof of property $(iii)$]
    By the strong law of large numbers, we have that with $Q$-probability one,
    \begin{align}
      \frac{1}{n} \log \eproc_n^\Qstar \to \EE_Q \left [ E_1(\theta_Q^\star) \right ] \equiv \ell_Q^\star.
    \end{align}
    Appealing to property $(ii)$ and sublinearity of $r_n$, we have that
    \begin{equation}
      \frac{1}{n} \log (\eproc_n) = \frac{1}{n} \log(\eproc_n / \eproc_n^\Qstar) + \frac{1}{n} \log(\eproc_n^\Qstar) \to \ell_Q^\star
    \end{equation}
    $Q$-almost surely for every $\Qinbar$.

    To demonstrate unimprovability of this limit for any other $e$-process
    $\eproc' \in \eprocclass(\Theta)$, we note that $\log (\eproc_n') =
    \log(\eproc_n' / \eproc_n^\Qstar) + \log (\eproc_n^\Qstar)$ and hence with
    $Q$-probability one,
    \begin{equation}
      \limsup_\nto \frac{1}{n} \log (\eproc_n') \leq \underbrace{\limsup_\nto
      \frac{1}{n} \log (\eproc_n' / \eproc_n^\Qstar)}_{\leq 0} +
      \underbrace{\limsup_\nto \frac{1}{n} \log(\eproc_n^\Qstar)}_{=
      \ell_Q^\star} 
      \label{proof-log-opt}
    \end{equation}
    where the first term is nonpositive by the numeraire property (\cref{assumption:abstract-numeraire-portfolios}) combined with the Borel-Cantelli lemma:
    \begin{align}
      \forall \delta>0,\quad\sum_{n=1}^\infty \Q \left ( \frac{1}{n} \log(\eproc_n' / \eproc_n^\Qstar) \geq \delta \right ) &= \sum_{n=1}^\infty \Q (\eproc_n' / \eproc_n^\Qstar \geq \exp \{ n\delta \} )\\
                                                                                                                               &\leq \sum_{n=1}^\infty \underbrace{\EE_Q[ \eproc_n' / \eproc_n^\Qstar ]}_{\leq 1} \exp \{ -n\delta \} < \infty
    \end{align}
    and the second term is exactly $\ell_Q^\star$ by the strong law of large numbers. This completes the proof for property $(iii)$.
  \end{proof}
    \begin{proof}[Proof of property $(i)$]
      Let $\eproc' \in \eprocclass(\Theta)$ be any other $e$-process in $\eprocclass(\Theta)$. Then we have for any $\Qinbar$,
      \begin{align}
        \log(\eproc_n / \eproc_n') &= \log \left ( \frac{\eproc_n}{\eproc_n^\Qstar} \cdot \frac{\eproc_n^\Qstar}{\eproc_n'} \right ) \\
        &= \log ( \eproc_n / \eproc_n^\Qstar) + \log(\eproc_n^\Qstar / \eproc_n'),
      \end{align}
      and hence we have that for any $\Qinbar$,
      \begin{align}
        \MoveEqLeft{\liminf_\nto \frac{1}{n} \left ( \log(\eproc_n / \eproc_n^\Qstar) + \log (\eproc_n^\Qstar / \eproc_n') \right )}\\
        \geq\ &\liminf_\nto \frac{1}{n} \log(\eproc_n / \eproc_n^\Qstar) + \liminf_\nto \frac{1}{n} \log(\eproc_n^\Qstar / \eproc_n') \geq 0,
      \end{align}
      where the final inequality follows from asymptotic equivalence of $\eproc$ and $\eproc^\Qstar$---as shown in property $(ii)$---and from another application of the numeraire property with the Borel-Cantelli lemma. This completes the proof for property $(i)$ and hence the proof of \cref{lemma:adaptive optimality} altogether.
    \end{proof}

\subsection{Proof of \cref{lemma:general nonasymptotic bound}}\label{proof:general stopping time}
\begin{proof}
  We start with the upper bound and then later derive the lower bound.
  \paragraph*{Upper bounding $\EE_Q[\widebar \tau_\alpha]$}
  Throughout, let $\tau \equiv \widebar \tau_\alpha := \inf \{ n \in \NN : \eproc_n \geq 1/\alpha \}$ and define $\ell_Q^\star := \EE_Q[\log E_1(\theta_Q^\star)]$. Let $\delta \in (0, 1)$ be arbitrary and define
  \begin{equation}
    \eps := \frac{\delta}{1+\delta}\cdot \ell^\star_Q\quad\text{and}\quad
    m := \left \lceil \frac{\log(1/\alpha)}{\ell_Q^\star-\eps} \right \rceil = \left \lceil \frac{(1+\delta)\log(1/\alpha)}{\ell_Q^\star} \right \rceil. 
  \end{equation}
  Writing out the expected stopping time under $Q$ as an infinite sum of tail probabilities, we have
\begin{align}
  \EE_Q \left [ \tau \right ] &= \sum_{k=0}^\infty \Q ( \tau \geq k ) \\
                              &\leq m + \sum_{k=m}^\infty \Q \left ( \tau \geq k \right )\\
                              &\leq 1 + m + \sum_{k=m}^\infty \Q \left ( \tau > k \right )\\
  &\leq 1 + m +  \underbrace{\sum_{k=m}^\infty \Q ( \eproc_k < 1/\alpha )}_{(\star)}.
\end{align}
where the first inequality upper bounds the first $m$ elements of the series by 1, the second does the same for one additional term and re-writes $\tau \geq k$ as $\tau > k-1$, and the third inequality follows from the definition of $\tau$ as the first $n$ such that $\eproc_n \geq 1/\alpha$.
Analyzing the sum in $(\star)$, we have
\begin{align}
  (\star) &= \sum_{k=m}^\infty \Q \left ( W_k < 1/\alpha \right )\\
  \verbose{
          &= \sum_{k=m}^\infty  \Q \left ( k^{-1}\log  W_k < k^{-1} \log(1/\alpha) \right )\\
                                                      &\leq \sum_{k=m}^\infty \left [ \Q \left ( \ell_Q^\star - \eps < k^{-1}\log(1/\alpha) \right ) + \Q \left ( \left \lvert  k^{-1} \log W_k - \ell_Q^\star \right \rvert \geq \eps \right ) \right ]\\
  }
                                                      &= \sum_{k=m}^\infty \Q (k^{-1} \log \eproc_k < k^{-1}\log(1/\alpha))\\
                                                      &\leq \sum_{k=m}^\infty \left [ \1 \left \{ \ell_Q^\star - \eps < m^{-1}\log(1/\alpha) \right \} + \Q \left ( \left \lvert  k^{-1} \log W_k - \ell_Q^\star \right \rvert \geq \eps \right ) \right ] \\
  &= \sum_{k=m}^\infty \Q \left ( \left \lvert  k^{-1} \log W_k - \ell_Q^\star \right \rvert \geq \eps \right ),
\end{align}
where the final equality follows from the definition of $m := \left \lceil \log(1/\alpha) / (\ell_Q^\star - \eps) \right \rceil $.
Now, letting $\eproc_n^\star := \prod_{i=1}^n E_i(\theta_Q^\star)$ and analyzing the sum in the last line, we have
\begin{align}
  (\star) &\leq \sum_{k=m}^\infty \Q \left ( \left \lvert  k^{-1} \log W_k - \ell_Q^\star \right \rvert \geq \eps \right ) \\
  \verbose{
  &\leq \sum_{k=m}^\infty \left [ \Q \left ( |k^{-1}\log W_k - k^{-1}\log W_k^\star | \geq \eps / 2 \right ) + \Q \left ( |k^{-1}\log W_k^\star - \ell_Q^\star | \geq \eps / 2 \right ) \right ] \\
  }
  &\leq \underbrace{\sum_{k=m}^\infty  \Q \left ( |k^{-1}\log W_k - k^{-1}\log W_k^\star | \geq \eps / 2 \right )}_{(\dagger)} + \underbrace{\sum_{k=m}^\infty \Q \left ( |k^{-1}\log W_k^\star - \ell_Q^\star | \geq \eps / 2 \right )}_{(\dagger\dagger)},
\end{align}
and we will provide upper bounds on $(\dagger)$ and $(\dagger\dagger)$ separately. First, notice that the $k^\tth$ summand of $(\dagger)$ can be bounded as follows for $\Delta_k := k^{-1} \log(W_k / W_k^\star)$:
\begin{align}
  \Q \left ( | \Delta_k | \geq \eps / 2 \right ) &\leq \Q \left ( \Delta_k \geq \eps / 2 \right ) + \Q \left ( -\Delta_k \geq \eps / 2 \right ) \\
          &\leq \Q \left (\log (W_k / W_k^\star) \geq k \eps/ 2 \right ) + \Q \left ( k^{-1} \sup_{\thetain} \log W_k(\theta) - k^{-1} \log W_k \geq \eps/2 \right ),
  \end{align}
  where the first inequality follows from a union bound and the second follows from the trivial fact that $\sup_\thetain \log W_k(\theta) \geq \log W_k^\star \equiv \log W_k(\theta_Q^\star)$ with $Q$-probability one. Appealing to the numeraire property to control the first term and the regret bound to control the second, we have that
  \begin{align}
          \Q(|\Delta_k| \geq \eps / 2) &\leq \Q \left (W_k / W_k^\star \geq \exp \{ k \eps/ 2 \} \right ) + \1 \left \{ r_k / k \geq \eps/2 \right \} \\
    \verbose{
          &\leq \EE_Q \left ( W_k / W_k^\star \right ) \exp \left \{ -k \eps / 2 \right \} + \1 \left \{ r_k / k \geq \eps/2 \right \} \\
    }
          &\leq \exp \left \{ -k \eps / 2 \right \} + \1 \left \{ r_k / k \geq \eps/2 \right \}.
\end{align}
Therefore, $(\dagger)$ can be bounded as
\begin{align}
  (\dagger) &\leq 1 + \sum_{k=m+1}^\infty \exp \left \{ -k\eps / 2 \right \} + \sum_{k=m}^\infty \1 \{ r_k / k \geq \eps / 2 \} \\
  \verbose{
  &\leq 1 + 1 + \int_{m}^\infty \exp \left \{ -u\eps / 2 \right \} \dd u + \sum_{k=m}^\infty \1 \{ r_k / k \geq \eps / 2 \} \\
  }
            &\leq 1 + \frac{2}{\eps}\exp \left \{ -m \eps / 2 \right \} + \sum_{k=m}^\infty \1 \{ r_k / k \geq \eps / 2 \} \\
  \verbose{
            &\leq 1 + \frac{2}{\eps}\exp \left \{ -\left [ \frac{\log(1/\alpha) \cancel{(1+\delta)}}{\cancel{\ell_Q^\star}} \right ] \frac{\delta\cancel{\ell_Q^\star}}{2\cancel{(1+\delta)}} \right \} + \sum_{k=m}^\infty \1 \{ r_k / k \geq \eps / 2 \} \\
  }
            &\leq 1 + \frac{2(1+\delta)}{\delta \ell_Q^\star}\exp \left \{ -\frac{\delta \log(1/\alpha) }{2}  \right \} + \sum_{k=m}^\infty \1 \{ r_k / k \geq \eps / 2 \} \\
            &= 1 + \frac{2(1+\delta)\alpha^{\delta /2 }}{\delta \ell_Q^\star}  + \sum_{k=m}^\infty \1 \{ r_k / k \geq \eps / 2 \}, 
\end{align}
where we observe that the third term is finite if and only if $r_k$ is sublinear.
Turning now to $(\dagger\dagger)$, we have by the Chebyshev-Nemirovski inequality derived in \cref{lemma:concentration-via-nemirovski} combined with the assumption that the $s^\tth$ moment is finite for some $s > 2$,
\begin{align}
  (\dagger \dagger) &\leq \sum_{k=m}^\infty \Q (k^{-1} \log W_k^\star - \ell_Q^\star \geq \eps / 2) \\
  \verbose{
  &\leq 1 + \sum_{k=m-1}^\infty\Q (k^{-1} \log W_k^\star - \ell_Q^\star \geq \eps / 2)\\
  }
                    &\leq 1 + \frac{2^s}{\eps^{s}} \sum_{k=m+1}^\infty \frac{\rho_s(\theta_Q^\star)}{k^{s/2}}\\
                    &\leq 1 + \frac{2^s\rho_s(\theta_Q^\star)}{\eps^{s}} \int_{m}^\infty \frac{1}{u^{s/2}}\dd u\\
  &= 1 + \frac{2^s\rho_s(\theta_Q^\star)}{\eps^{s}m^{s/2 - 1}(s/2 - 1)}.
\end{align}
Plugging in $\eps := \delta \ell_Q^\star / (1+\delta)$ and $m = \left \lceil \log(1/\alpha) (1+\delta) / \ell_Q^\star \right \rceil$, we have
\begin{align}
  (\dagger \dagger) &\leq 1 + \frac{2^s \rho_s(\theta_Q^\star)}{\left ( \delta \ell_Q^\star / (1+\delta) \right )^s \left \lceil \log(1/\alpha) (1 + \delta) / \ell_Q^\star \right \rceil^{s/2 - 1} (s/2 - 1) } \\
                    &\leq 1 + \frac{(1+\delta)^s2^s \rho_s(\theta_Q^\star)}{\left ( \delta \ell_Q^\star \right )^s (\log(1/\alpha) (1 + \delta) / \ell_Q^\star )^{s/2 - 1} (s/2 - 1) } \\
    &= 1 + \frac{(1+\delta)^{1 + s/2} 2^s \rho_s(\theta_Q^\star)}{\delta^s(\ell_Q^\star)^{1+s/2} \log^{s/2-1}(1/\alpha)(s/2 - 1)} \,.
\end{align}
Putting everything together, we have that for any $\delta > 0$,
\begin{align}
\label{eq:proof-general-expected-rejection-time-upper-bound}
  &\EE_Q \left [ \tau \right ] \leq 1 + m + (\dagger) + (\dagger\dagger) \\
  &\leq 4 + \frac{(1+\delta)\log(1/\alpha)}{\ell_Q^\star} + \frac{2(1+\delta)\alpha^{\delta/2}}{\delta \ell_Q^\star} \\
  &\qquad+ \frac{(1+\delta)^{1+s/2} 2^s \rho_s(\theta_Q^\star)}{\delta^s (\ell_Q^\star)^{1+s/2} \log^{s/2 - 1}(1/\alpha)(s/2-1)} + \sum_{k=m}^\infty \1 \left \{ \frac{r_k}{k} \geq \frac{\delta \ell_Q^\star}{2(1+\delta)} \right \} \,.
\end{align}

  \paragraph*{Lower bounding $\EE_Q[\widetilde \tau_\alpha]$}
  Let $\infseqn{\theta_n}$ be an arbitrary predictable sequence and let $\tau \equiv \tau_\alpha := \inf \{ n \in \NN : \eproc_n(\theta_1,\dots, \theta_n) \geq 1/\alpha \}$ be the stopping time of its associated test. Assume that $\EE_Q[\tau]< \infty$ since otherwise the result is immediate. Notice that by definition of $\tau_\alpha$, we have that for any $Q \in \Qcal$,
  \begin{equation}
    \eproc_\tau \geq 1/\alpha
  \end{equation}
  with $Q$-probability one.
  Multiplying both sides by $\eproc_\tau^\star / \eproc_\tau^\star$ where $\eproc_n^\star = \eproc_n(\theta_Q^\star)$ under the $(\Theta,Q)$-numeraire portfolio, and taking logarithms and expectations under $Q$, we have the following inequality:
  \begin{align}
    \MoveEqLeft{\EE_Q \left [ \log ( \eproc_\tau / \eproc_\tau^\star) \right ] + \EE_Q \left [ \log \eproc_\tau^\star \right ]}\\
    & = \EE_Q \left [ \log ( \eproc_\tau / \eproc_\tau^\star) \right ] + \EE_Q \left [ \sum_{i=1}^\tau \log\left ( E_i(\theta_Q^\star) \right ) \right ]  \geq \log(1/\alpha).
  \end{align}
  Applying Wald's identity to the second term on the left-hand side and rearranging, we have that
  \begin{align}
    \EE_Q(\tau) \ell_Q^\star &\equiv \EE_Q (\tau) \cdot \EE_Q \left [ \log \left ( E_1(\theta_Q^\star) \right ) \right ] \\
                                          &\geq \log(1/\alpha) - \EE_Q \left [ \log(\eproc_\tau / \eproc_\tau^\star) \right ]\\
    &\geq \log(1/\alpha) - \underbrace{\log \EE_Q \left [ \eproc_\tau / \eproc_\tau^\star \right ]}_{(\star)},\label{eq:stopping time lower bound proof after applying wald}
  \end{align}
  where the second inequality follows from an application of Jensen's inequality. Moreover, we have that $(\star) \leq 0$ by \cref{assumption:abstract-numeraire-portfolios}. Dividing both sides by $\ell_Q^\star$ yields the desired result:
  \begin{equation}
    \EE_Q (\tau) \geq \frac{\log(1/\alpha)}{\ell_Q^\star},
  \end{equation}
  which completes the proof of the lower bound and hence of \cref{lemma:general nonasymptotic bound}.
\end{proof}

\subsection{Proofs for \cref{example:sub-gaussian}}
\label{proof:sub-gaussian}
In this subsection we first present the derivation of the test supermartingale whose portfolios are mixed using a Gaussian prior. We then prove the regret result for the mixed test supermartingale.
\begin{proof}[Derivation of the mixed test supermartingale]
    We now present the derivation of the mixed test supermartingale. Throughout we let $\lambda = \lambda_i$ for each $i \in [n]$ and for some $\lambda \in \RR$. Defining the shorthand $S_n \coloneqq \sum_{i=1}^n (X_i - \mu_0)$, we have that the test supermartingale takes the following form
    \begin{equation}
        W_n = \exp\left\{\lambda S_n - \frac{n \sigma^2 \lambda^2}{2}\right\} \,.
    \end{equation}
    We use as our mixing prior over the portfolios the centered normal distribution with variance $\sigma^2$. We now proceed to characterize the resulting, mixed, test supermartingale. Indeed, integrating over the portfolio $\lambda$ using centered normal prior we can see how
    \begin{align}
        W_n^{\mathrm{mix}} 
        &\coloneqq \frac{1}{\sqrt{2 \pi \sigma^2}}\int_{\RR} \exp\left\{\lambda S_n - \frac{n \sigma^2 \lambda^2}{2}\right\} \exp\left\{-\frac{\lambda^2}{2 \sigma^2}\right\} \dd\lambda \\
        &= \frac{1}{\sqrt{2 \pi \sigma^2}} \int_{\RR} \exp\left\{\lambda S_n - \lambda^2 \left(\frac{n \sigma^4 + 1}{2 \sigma^2}\right)\right\} \dd\lambda \\
        &= \frac{1}{\sqrt{2 \pi \sigma^2}} \int_{\RR} \exp\left\{\frac{2\lambda \sigma^2 S_n - \lambda^2 (n \sigma^4 + 1))}{2 \sigma^2}\right\} \dd\lambda \\
        &= \frac{1}{\sqrt{2 \pi \sigma^2}} \int_{\RR} \exp\left\{\frac{-a (\lambda^2 - 2\lambda b/a)}{2 \sigma^2}\right\} \dd\lambda \,,
        \label{eq:proof-mixed-supermartingale}
    \end{align}
    where in the last step we defined the shorthands $a \coloneqq (n \sigma^4 + 1)$ and $b \coloneqq \sigma^2 S_n$.

    We continue by completing the square inside the exponential function. Namely, we have that
    \begin{align}
        \lambda^2 - 2 \lambda \frac{b}{a}
        &= \lambda^2 - 2 \lambda \frac{b}{a} + \left(\frac{b}{a}\right)^2 - \left(\frac{b}{a}\right)^2 \\
        &= \left(\lambda - \frac{b}{a} \right)^2 - \left(\frac{b}{a}\right)^2 \,.
    \end{align}
    Plugging the identity into \eqref{eq:proof-mixed-supermartingale} results in the following expression
    \begin{align}
        \eqref{eq:proof-mixed-supermartingale}
        &= \frac{1}{\sqrt{2 \pi \sigma^2}} \int_{\RR} \exp\left\{\frac{-a \left[(\lambda - b/a)^2 - (b/a)^2\right]}{2 \sigma^2}\right\} \dd\lambda \\
        &= \exp\left\{\frac{b^2 / a}{2 \sigma^2}\right\} \frac{1}{\sqrt{2 \pi \sigma^2}} \int_{\RR} \exp \left\{-\frac{(\lambda - b/a)^2}{2 \sigma^2 / a}\right\} \dd\lambda \\
        &= \frac{\exp\left\{\frac{b^2 / a}{2 \sigma^2}\right\}}{\sqrt{a}} \frac{1}{\sqrt{2 \pi \sigma^2/a}} \int_{\RR} \exp \left\{-\frac{(\lambda - b/a)^2}{2 \sigma^2 / a}\right\} \dd\lambda \,,
        \label{eq:proof-mixed-supermartingale-transformed-integral}
    \end{align}
    where in the last step we multiplied the expression by $\sqrt{a} / \sqrt{a}$ and rearranged the terms. Notice how \eqref{eq:proof-mixed-supermartingale-transformed-integral} consists of a term multiplied by the density of a normal distribution with mean $b/a$ and variance $\sigma^2 / a$. Since the density in \eqref{eq:proof-mixed-supermartingale-transformed-integral} integrates to one, we can see that when we substitute in the definitions of $a$ and $b$ the expression simplifies to
    \begin{align}
        \eqref{eq:proof-mixed-supermartingale-transformed-integral}
        = \frac{\exp\left\{\frac{\sigma^2 S_n^2}{2(n \sigma^4 + 1)}\right\}}{\sqrt{n \sigma^4 + 1}} \,,
    \end{align}
    which completes the derivation of our mixed test supermartingale.
\end{proof}

\begin{proof}
    We now prove that $W_n^{\mathrm{mix}}$ enjoys logarithmic regret with $Q$-probability one. To do so, we first characterize the $\Pe$-process that uses the best in hindsight portfolio and subsequently prove the logarithmic regret.
    
    The best in hindsight portfolio is the solution to the following maximization problem
    \begin{equation}
        \lambda^{\mathrm{max}} \coloneqq \argmax_{\lambda \in \RR} \left\{\lambda S_n - \lambda^2 \frac{n \sigma^2}{2}\right\} \,,
    \end{equation}
    where we recall that $S_n \coloneqq \sum_{i=1}^n (X_i - \mu_0)$. Through some direct computations we can see how $\lambda^{\mathrm{max}} = S_n / (n\sigma^2)$, which means that the $\Pe$-process that uses this best in hidnsight portfolio takes the form
    \begin{equation}
        \log\left(W_n^{\mathrm{max}}\right) = \frac{S_n^2}{2n\sigma^2} \,.
    \end{equation}

    Hence, we can now directly compute the regret of $W_n^{\mathrm{mix}}$ as:
    \begin{align}
        \log\left(W_n^{\mathrm{max}}\right) - \log\left(W_n^{\mathrm{mix}}\right)
        &= \frac{S_n^2}{2n\sigma^2} - \frac{\sigma^2 S_n^2}{2(n\sigma^4 + 1)} + \frac{1}{2}\log\left(n\sigma^4 + 1\right) \\
        &= \frac{S_n^2}{2n\sigma^2(n\sigma^4 + 1)} + \frac{1}{2}\log\left(n \sigma^4 + 1\right) \,.
    \end{align}
    We can therefore see that $W_n^{\mathrm{mix}}$ attains logarithmic regret with $Q$-probability one, completing the proof.
\end{proof}